\documentclass{amsart}
\usepackage[english]{babel}
\usepackage{graphicx}
\pdfoutput=1
\usepackage[hidelinks]{hyperref}
\usepackage{amssymb}
\usepackage{amsmath}
\usepackage[foot]{amsaddr}
\usepackage{amsfonts}
\usepackage{bbm}
\usepackage[dvipsnames, table]{xcolor}
\usepackage[top=1in, bottom=1.25in, left=1.25in, right=1.25in]{geometry}
\usepackage{caption}
\usepackage{subcaption}
\usepackage{multirow}

\graphicspath{{images/}}

\input{packages.sty}

\definecolor{maria}{HTML}{0090A0}

\colorlet{BLUE}{blue}
\colorlet{BLACK}{black}

\newcommand{\A}[1]{{\color{black}#1}}

\newcommand{\C}[1]{{\color{black}#1}}

\begin{document}

\title[Driving bifurcating nonlinear PDE($\bmu$)\lowercase{s} by \ocp \lowercase{s}]{Driving bifurcating parametrized nonlinear PDE\lowercase{s} by optimal control strategies: application to Navier-Stokes equations with model order reduction}

%\author{Federico Pichi$^{1,*}$, Maria Strazzullo$^{1,*}$, Francesco Ballarin$^1$ and Gianluigi Rozza$^1$}
%\address{$^1$ \textup{mathLab, Mathematics Area, SISSA, via Bonomea 265, I-34136 Trieste, Italy.}}
%\address{$^*$ \textup{Both authors contributed equally to this work.}}

\author{Federico Pichi$^{1, 2, *}$}
\author{Maria Strazzullo$^{1, *}$}
\author{Francesco Ballarin$^3$}% \address{}
\author{Gianluigi Rozza}
\address{ $^1$ \textup{mathLab, Mathematics Area, SISSA, via Bonomea 265, I-34136 Trieste, Italy.}}
\address{$^2$ \textup{Chair of Computational Mathematics and Simulation Science, \'Ecole Polytechnique F\'ed\'erale de Lausanne, Lausanne, Switzerland}}
\address{$^3$ \textup{Department of Mathematics and Physics, Catholic University of the Sacred Heart, Brescia, Italy}}
\address{$^*$ \textup{Both authors contributed equally to this work.}}

\begin{abstract}
This work deals with optimal control problems as a strategy to drive bifurcating solution of nonlinear parametrized partial differential equations towards a desired branch. Indeed, for these governing equations, multiple solution configurations can arise from the same parametric instance. We thus aim at describing how optimal control allows to change the solution profile and the stability of state solution branches. First of all, a general framework for nonlinear optimal control problem is presented in order to reconstruct each branch of optimal solutions, discussing in detail the stability properties of the obtained controlled solutions. Then, we apply the proposed framework to several optimal control problems governed by bifurcating Navier-Stokes equations in a sudden-expansion channel, describing the qualitative and quantitative effect of the control over a pitchfork bifurcation, and commenting in detail the stability eigenvalue analysis of the controlled state. Finally, we propose reduced order modeling as a tool to efficiently and reliably solve parametric stability analysis of such optimal control systems, which can be challenging to perform with standard discretization techniques such as Finite Element Method.
\end{abstract}

\maketitle

\section*{Introduction}
\label{intro}
Nonlinear parametrized partial differential equations (PDE($\boldsymbol{\mu}$)s) play a fundamental role in several fields, ranging from Continuum Mechanics to Quantum Mechanics, passing through Fluid Dynamics.
When compared to linear PDE($\boldsymbol{\mu}$)s with smooth parametric dependence, the most striking difference is that in the linear case a stable solution evolves in a continuous (and thus, unique) manner when the parameter changes slightly. In contrast, in the nonlinear case the solution for a given parameter $\bmu$ of a given PDE($\boldsymbol{\mu})$ may not be unique. Indeed, the model can suddenly change its behavior, together with the stability properties of its solutions. Models with such a feature are called bifurcation problems \cite{Prodi,Caloz,seydel2009practical}. Examples of models which are characterized by the non-uniqueness of the solution are the Von K\'arm\'an plate model for buckling \cite{vonka,bauerreiss,berger, pichirozza}, the Gross-Pitaevskii equation for Bose-Einstein condensates \cite{Middelkamp_et_al2011,doi:10.1137/1.9781611973945,Charalampidis_et_al2018,pichiquaini} and Navier-Stokes equations in a channel \cite{AQpreprint, cardio, pintore2019efficient}.

The critical point at which the system loses its original features is called bifurcation point, and it is usually denoted by $\bmu^*$. Many different bifurcating configurations can emerge from these points, and an intuitive way to visualize them is through a plot of a scalar output of the solution against the parameter value for which the solution has been computed, e.g.\ the so called \textit{bifurcation diagram}. As an example, a branch of solutions can give rise at $\bmu^*$ to two further symmetric branches of solutions which coexist with the pre-existing one, and switch their stability properties with it. When such situation occurs, the model is said to undergo a \textit{pitchfork bifurcation} phenomenon \cite{seydel2009practical}.
In particular, we are interested in an application in Fluid Dynamics, where we consider sudden-expansion channel flows, which are motivated by many practical scenarios. Flow profiles are described by Navier-Stokes equations, parametrized by the viscosity value $\mu$.
Here we consider a simplified version for a model of a cardiac disease, called mitral valve regurgitation, that may cause either a symmetric regurgitant jet or a wall-hugging, non-symmetric regurgitant jet.
The latter phenomenon, which can be clinically detected through echocardiography, is called the Coanda effect \cite{tritton2012physical}, and expresses the tendency of a fluid jet to be attracted to a nearby surface. This represents an issue from the medical point of view, because the wall-hugging jet might lead to inaccurate echocardiography measurements. It is therefore of the utmost practical interest to try to drive the system towards the branch of symmetric solutions, which are more favorable for the measurement process.

Towards this goal, parametrized optimal control problems (\ocp s) governed by PDE($\boldsymbol{\mu}$)s might be employed. Indeed, \ocp s can be interpreted as an input-output system which achieves an observable configuration \cite{bochev2009least, gunzburger2003perspectives, hinze2008optimization, lions1971, troltzsch2010optimal}. They have been exploited in several applications in different scientific fields, see e.g.\ \cite{leugering2014trends} for an overview.
The main goal of this work is to use \ocp s to drive bifurcating state profiles towards a different desired state, which might possibly belong to another state solution branch.

In such \ocp s, a parametric study of the state solution is necessary in order to understand the behavior of the system.
The complexity of this task can be tackled using a combination of existing methodologies, which allow the complete reconstruction of the aforementioned bifurcation diagram.
However, the discretization through standard techniques can lead to a possibly huge system to be solved for many values of the parameters. In this work we exploit Galerkin Finite Element (FE) approach, which can be challenging when several instances of parametrized problem have to be studied, most of all because the optimal control setting requires additional equations to be solved on top of the original PDE($\boldsymbol{\mu}$)s. For this reason we also propose Reduced Order Modeling (ROM) as a tool to overcome this issue, see \cite{hesthaven2015certified} as an introductory reference. Namely, we build a low-dimensional space from FE \emph{high-fidelity} optimal solutions and we perform a Galerkin projection in a \emph{reduced space}, usually much smaller than the considered high-fidelity one. The main ingredients are the reduced space construction exploiting Proper Orthogonal Decomposition (POD) algorithm \cite{ballarin2015supremizer, burkardt2006pod, Chapelle2013} and then, using a Galerkin projection, solving a reduced problem in the lower dimensional reduced space, for every parameter instance.

The main novelty of this work is the formulation, numerical simulation and subsequent model reduction of \ocp s governed by bifurcating parametrized Navier-Stokes equations, together with the analysis of the properties of their eigenvalues. We develop several test cases corresponding to different control actions in order to understand and validate our findings. At the best of our knowledge, a thorough description of the mathematical analysis and widespread exploitation of \ocp s for bifurcating nonlinear PDE($\boldsymbol{\mu}$)s have not been extensively addressed in literature. This work could then pave the way towards application of \ocp s as a tool to drive the behavior of complex bifurcation problems to a desired branch of solutions.

The work is outlined as follows: in Section \ref{general_problem_sec} we describe the structure of \ocp s for general nonlinear bifurcating systems both at the continuous and discretized level, exploiting a branch-wise approach through a FE approximation. Furthermore, we introduce the basic notions of stability analysis and its connection to eigenvalue problems for the uncontrolled and the controlled systems. In Section \ref{sec:state} we introduce the uncontrolled sudden-expansion channel problem, while in Section \ref{NS_ocp} we build several \ocp s assigning different roles to the control variable. There we also perform a global eigenvalue analysis in order to understand the main features of the achieved optimal solutions. We present two different boundary control problems and two distributed ones, which give very different solution configurations. The ROM strategy is described in Section \ref{sec_ROM}: the presented approach is tested for all the numerical simulations performed at the FE level. Conclusions follow in Section \ref{conclusions}.

\section{Nonlinear Parametrized Optimal Control Problems and Bifurcating systems}
\label{general_problem_sec}
In this Section we introduce a generic nonlinear \ocp s. We will focus on the minimization of quadratic cost functional under nonlinear PDE($\bmu$) constraint for Hilbert spaces, following the Lagrangian approach \cite{gunzburger2003perspectives, hinze2008optimization}. In Section \ref{general_problem} and \ref{FE} we provide existence results and optimality conditions for nonlinear \ocp s in their continuous and discretized version, respectively. Then, Section \ref{bif} will describe the spectral properties of the optimization system at hand.

\subsection{Problem Formulation}
\label{general_problem}
Optimal control is a mathematical tool which aims at modifying the natural behavior of a system. Let us suppose to have a \emph{state} PDE($\bmu$)
\begin{equation}
\label{eq:state}
G(y; \bmu) = f,
\end{equation}
with \emph{state variable} $y \eqdot y(\bmu) \in \mathbb Y$, i.e.
$G: \mathbb Y \times \mathcal P \rightarrow \mathbb Y\dual$ where $\mathbb Y$ is a Hilbert space, $f \in \state\dual$ is a forcing term,
$\mathcal P \subset \mathbb R^P$
is a parameter space of dimension $P \geq 1$, while $G(y; \bmu) =
E_{\textit{n}\ell}(y; \bmu) + E_{\ell}(y; \bmu)$ is the \emph{state operator}, with $E_{\ell} \in \Cal L(\state, \state \dual)$ and $E_{\textit{n}\ell}$ representing the linear and nonlinear contributions, respectively.
Here, we call $\Cal L \cd$ the space of linear continuous functions between two spaces.
We now want $y$ to be the most similar to a known solution profile
$y_\text{d} \eqdot y_\text{d}(\bmu) \in \mathbb Y_{\text{obs}} \supseteq \mathbb Y$.
To this end, a new variable is introduced in the equation, the \emph{control variable} $u \eqdot u(\bmu) \in \mathbb U$, with $\mathbb U$ another possibly different Hilbert space. Let us define the \emph{controlled equation} 
%$\mathcal E(y,u; \bmu)$, where
$\mathcal E: \mathbb Y \times \mathbb U \times \mathcal P \rightarrow \mathbb Y\dual$ as 
%Then, the controlled equation will be of the following form:
\begin{equation*}
\mathcal E(y,u; \bmu) \eqdot\; G(y; \bmu) - C(u) - f = 0,
\end{equation*}
where $C \in \Cal L(\control, \state \dual)$ is the \emph{control operator} describing the action of the variable $u$ on the system\footnote{Parametrized control operators are also possible, with a straightforward extension of the methodology presented herein.}. In other words, we are trying to change the behavior of the state PDE($\bmu$) through $C(u)$.
The \ocp $\;$ reads: given a $\bmu \in \mathcal P$, find the pair $(y,u) \in \mathbb Y \times \mathbb U$ which solves
\begin{equation}
\label{min_problem}
\min_{y \in \mathbb Y, u \in \mathbb U} J(y,u; y_\text{d}) \text{ subject to } \mathcal E(y,u; \bmu) = 0,
\end{equation}
where $J: \mathbb Y \times \mathbb U \times \mathbb Y_{\text{obs}} \rightarrow \mathbb R$ is the \emph{objective functional} defined by
\begin{equation}
J(y,u; y_\text{d}) \eqdot \half \norm{y - y_\text{d}}_\mathbb {Y_{\text{obs}}}^2 + \alf \norm{u}_{\mathbb U}^2,
\end{equation}
and $\alpha \in (0, 1]$ is a \emph{penalization parameter}. The role of $\alpha$ is of great interest: indeed, a large value of $\alpha$ translates in a poor capability of the system to be controlled, while $\alpha \ll 1$ allows the functional to be minimized with larger values of the variable $u$. Problem \eqref{min_problem} admits a solution if \cite[Section 1.5.2]{hinze2008optimization}:
\begin{enumerate}[(i)]
\item $\mathbb U$ is convex, bounded and closed;
\item $\mathbb Y$ is convex and closed;
\item for every $\bmu \in \Cal P$, the controlled system $\mathcal E (y,u; \bmu) = 0$ has a bounded solution map
$u \in \mathbb U \mapsto y(u) \in \mathbb Y$;
\item for a given $\bmu \in \Cal P$, the map $(y,u, \bmu) \in \mathbb Y \times \mathbb U \times \mathcal P \rightarrow \mathcal E (y,u; \bmu) \in \mathbb Y\dual$ is weakly continuous with respect to (w.r.t.) the first two arguments;
\item for a given $y_{\text{d}} \in \mathbb Y_{\text{obs}}$, the objective functional $J(y,u; y_\text{d})$ is weakly lower semicontinuous w.r.t.\ $y$ and $u$.
\end{enumerate}
We now discuss the Lagrangian structure and the necessary first order optimality conditions. First of all, let $z \eqdot z(\bmu) \in {\mathbb Y^{\ast \ast}} = \mathbb Y$ be an arbitrary variable called \emph{adjoint variable}. Let us call $X = (y,u,z) \in \mathbb X \eqdot \mathbb Y \times \mathbb U \times \mathbb Y$ and let us build the \emph{Lagrangian functional}
$\Lg:\mathbb X \times \mathbb Y_{\text{obs}} \times \mathcal P \rightarrow \mathbb R$ as
\begin{equation}
\label{lagrangian_functional}
\Lg(X; y_\text{d}, \bmu) \eqdot J(y, u; y_\text{d}) + \la z, \Cal E(y,u; \bmu) \ra_{\mathbb Y \mathbb Y\dual},
\end{equation}
where $\la \cdot, \cdot \ra_{\mathbb Y \mathbb Y\dual}$ is the duality pairing of $\mathbb Y$ and $\mathbb Y\dual$. The introduction of the adjoint variable allows to treat problem \eqref{min_problem} in an unconstrained fashion finding the stationary point of \eqref{lagrangian_functional}. We remark that we consider $z$ in the same space of the state variable for a proper definition of the discretized problem: we will clarify the reason in Section \ref{FE}. Moreover, the variable $X$ will inherit the parameter dependence by definition, i.e $X \eqdot X(\bmu)$.
Furthermore, let us assume that the following hold:
\begin{enumerate}[resume*]
\item $\mathbb U$ is nonempty;
\item $J : \mathbb Y \times \mathbb U \times \mathbb Y_{\text{obs}} \rightarrow \mathbb R$ and $\mathcal E : \mathbb Y \times \mathbb U \times \mathcal P \rightarrow \mathbb Y\dual$ are continuously Fr\'echet differentiable w.r.t.\ $y$ and $u$;
\item given $\bmu \in \mathcal P$, the controlled system $ \mathcal E(y, u; \bmu) = 0$ has a unique solution $y = y(u) \in \mathbb Y$ for all $u \in \mathbb U$;
\item given $\bmu \in \mathcal P$, $D_y \mathcal E (y, u; \bmu) \in \mathcal L(\mathbb Y, \mathbb Y\dual)$ has a bounded inverse for all control variables $u$.
\end{enumerate}
The Fr\'echet derivative w.r.t.\ a variable $\star$ will be indicated as $D_{\star}$, as already done in (ix). Assuming
$({y}, u) \in \state \times \control$ to be a solution to \eqref{min_problem} for a given $\bmu \in \Cal P$, thanks to hypotheses
(vi) - (ix), there exists an adjoint variable $z \in \state$ such that the following variational system is satisfied \cite{hinze2008optimization}:
\begin{equation}
\label{KKT}
\begin{cases}
D_{y}\Lg(X; y_\text{d}, \bmu) [\omega] = 0 & \forall \omega \in \state,\\
D_u\Lg(X; y_\text{d}, \bmu) [\kappa] = 0 & \forall \kappa \in \control,\\
D_z\Lg(X; y_\text{d}, \bmu) [\zeta] = 0 & \forall \zeta \in \state,\\
\end{cases}
% \end{equation}
\qquad \text{or, in strong form,} \qquad
% \begin{equation}
% \label{KKT_explicit}
\begin{cases}
y + D_y \mathcal E (y, u; \bmu)\dual (z) = y_\text{d}, &\\
\alpha u - C\dual (z )= 0, &\\
\mathcal E(y,u;\bmu) = 0, &\\
\end{cases}
\end{equation}
where $D_y \mathcal E (y, u; \bmu) \dual \in \Cal L(\state, \state \dual)$ is the adjoint operator of the Fr\'echet linearization of $\mathcal E (y,u; \bmu)$ w.r.t.\ the state variable, while $C\dual \in \Cal L(\mathbb Y, \mathbb U\dual) $ is the adjoint of the control operator.
We will refer to problem \eqref{KKT} as \emph{optimality system}, in weak or strong form, respectively.
Moreover, writing the latter in compact form, it reads: given $\bmu \in \Cal P$, find $X \in \mathbb X$ such that
\begin{equation}
\label{ocp}
\mathcal G(X; \bmu) = {\Cal F},
\end{equation}
with
\begin{equation*}
\mathcal G(X; \bmu) \eqdot \begin{bmatrix} y + D_y \mathcal E (y, u; \bmu)\dual (z) \\ \alpha u - C\dual (z) \\ G(y, \bmu) - C(u) \end{bmatrix} \quad
\text{and} \quad \Cal F \eqdot \begin{bmatrix} y_{\text{d}} \\ 0 \\ f \end{bmatrix}.
\end{equation*}

In the nonlinear case, even considering only the state equation, the solution for a given parameter $\bmu$ may not be unique. Therefore, the local well-posedness of the problem \eqref{ocp}, which strongly relies on its local invertibility assumptions (viii) and (ix), can fail due to the singularity of the state equation.
Therefore, we can talk about \emph{solution branches}, i.e.\ multiple solution behaviors for a given value of the parameter.
We denote by $k$ the number of branches, and by $\Cal X_i$, $i = 1, \hdots, k$ the set of each solution on the $i$-th branch. We call \emph{solution ensemble} the set of all the solution branches $\Cal X_i$:
\begin{equation}
\label{ensemble}
\Cal X \eqdot \bigcup_{i = 1}^k \{ X(\bmu) \in \Cal X_i \; | \; \bmu \in \mathcal P \}.
\end{equation}
In the next Section, we will discuss the FE approximation of a solution for a fixed value of the parameter to the nonlinear \ocp , restricting ourselves to the well-posed setting.

\subsection{The FE Approximation}
\label{FE}
We are interested in the numerical approximation of the solution ensemble \eqref{ensemble} of the nonlinear \ocp \, in \eqref{ocp}, defined over an open and bounded regular domain $\Omega \subset \mathbb R^d$.
Our aim is to discretize the problem at hand, in order to investigate its qualitative changes w.r.t.\ the values of the parameter. % As we said, an overall viewpoint can be obtained through the bifurcation diagram, thus, in this Section, we will review the building blocks to approximate a specific solution $X(\bmu)$ for a chosen branch.
We remark that, even considering a single branch, the fulfillment of the well-posedness conditions can fail at some critical point $\bmu^{*}$.
Hence, in the following, we assume $\bmu \neq \bmu^*$ and $X(\bmu) \in \Cal X_i$ for some $i \in \{1, \dots, k\}$, thus we call $\Cal X_i$ a \textit{non-singular branch}.
Furthermore, we assume the nonlinearity to be at most quadratic in the state variable, \A{guided by the numerical results we are going to present in the following Sections. However, the structure and the methodology do not change when dealing with nonlinearities of higher order.}

To approximate the system in \eqref{ocp},
first of all we define the triangulation $\Cal{T^{N_T}}$ of $\Omega$, where $K$ is an element of $\Cal{T^{N_T}}$ and $\Cal N_{\Cal T}$ is the number of cells. Then, let us consider the discrete spaces $\state \discy = \state \cap\mathbb K_{r_y}$
and $\control \discu = \control \cap \mathbb K_{r_u}$, where
$
\mathbb K_r = \{ v \in C^0(\overline \Omega) \; : \; v |_{K} \in \mbb P^r, \; \; \forall \, K \in \Cal{T^{N_T}} \},
$
and $\mbb P^r$ is the space of all the polynomials of degree at most equal to $r$. Let us consider the FE function space $\mbb X \disc = \state \discy \times \control \discu \times \state \discy \subset \mbb X$, of dimension
$\Cal N = 2\Cal N_y + \Cal N_u$. The FE approximation of the parametrized problem \eqref{ocp} reads: given $\bmu \in \Cal P$, find $X\disc \eqdot X \disc (\bmu) \in \mathbb X\disc$ such that
\begin{equation}
\label{FE_ocp}
\mathcal G(X\disc; \bmu) = {\mathcal F}.
\end{equation}
We now want to make the algebraic structure of the system \eqref{FE_ocp} explicit. After the FE discretization, we can define $\mathsf y$, $\mathsf u$ and $\mathsf z$ as the column vectors which entries are given by the FE coefficients of the state, control and adjoint variables in their approximated spaces, respectively. In the same fashion, we call $\mathsf y_\text{d}$ the column vector of the FE coefficients representing the desired state profile. Let us focus on the structure of the optimality problem. At the FE level, applying our controlled state to the FE basis, we can derive the matrices
$\mathsf {E}_{\textit{n}\ell} + \mathsf {E}_{\ell} - \mathsf C$ and the forcing term vector $\mathsf f$. Moreover, we define the mass matrices $\mathsf M_y$ and $\mathsf M_u$ for state/adjoint variables and control, respectively. We still need to understand the algebraic structure of $D_y \Cal E(y,u; \bmu)$. %Assuming the operator $G$ to be quadratic in the variable $y$, \A{toglie frase prima} 
The Fr\'echet derivative of the controlled state equation w.r.t.\ the state $y$ will be
$\mathsf {E}_{\textit{n}\ell}'[\mathsf y] + \mathsf E_{\ell}$. In other words, the linear state structure is preserved, the nonlinear operator is linearized in $\mathsf {E}_{\textit{n}\ell}'[\mathsf y] $ and the control contribution disappears. Then, the global matrix formulation of the optimization system \eqref{FE_ocp} is
\begin{equation}
\label{algebra_ocp}
\overbrace{
\begin{bmatrix}
\mathsf M_y & 0 & \mathsf {E}_{\textit{n}\ell}'[\mathsf y]^T + \mathsf E_{\ell}^T\\
0 & \alpha \mathsf M_u & - \mathsf C^T \\
\mathsf {E}_{\textit{n}\ell} + \mathsf {E}_{\ell} & - \mathsf C & 0 \\
\end{bmatrix}\underbrace{
\begin{bmatrix}
\mathsf y \\
\mathsf u \\
\mathsf z \\
\end{bmatrix}}_{\mathsf X}}^{\mathsf G(\mathsf X; \boldsymbol \mu)}
=
\overbrace{
\begin{bmatrix}
\mathsf M_y \mathsf y_\text{d} \\
0 \\
\mathsf f \\
\end{bmatrix}}^{\mathsf F},
\end{equation}
which in compact form reads:
\begin{equation}
\label{G_compact_FE}
\mathsf R(\mathsf X; \bmu)\eqdot \mathsf G(\mathsf X; \boldsymbol \mu) - \mathsf F = 0,
\end{equation}
where $\mathsf R(\mathsf X; \bmu)$ represents the \emph{global residual} of the optimality system.
To solve system \eqref{G_compact_FE}, we rely on Newton's method and we solve
\begin{equation}
\mathsf {X}^{j + 1} = \mathsf {X}^j+ \mathsf{Jac}(\mathsf X^{j}; \boldsymbol \mu)^{-1}(\mathsf F - \mathsf G(\mathsf X^j; \boldsymbol \mu)), \spazio j \in \mathbb N,
\end{equation}
until a residual based convergence criterion is satisfied.
Since the matrix $\mathsf {E}_{\textit{n}\ell}'[\mathsf y]^T$ still depends on $\mathsf y$, the Jacobian matrix will be of the following nature:
\begin{equation}
\label{J_ocp}
\mathsf{Jac}(\mathsf X^j; \bmu) = \begin{bmatrix}
\mathsf M_y + \mathsf D_{\mathsf y}( \mathsf {E}_{\textit{n}\ell}'[\mathsf y^j]^T)[\mathsf z^j] & 0 & \mathsf {E}_{\textit{n}\ell}'[\mathsf y^j]^T + \mathsf E_\ell^T\\
0 & \alpha \mathsf M_u & - \mathsf C^T \\
\mathsf {E}_{\textit{n}\ell}'[\mathsf y^j] + \mathsf E_\ell & - \mathsf C & 0 \\
\end{bmatrix},
\end{equation}
where the matrix $\mathsf D_{\mathsf y}( \mathsf {E}_{\textit{n}\ell}'[\mathsf y^j]^T)[\mathsf z^j]$ does not depend anymore on the state, but only on the $j-$th value of the adjoint variable. %This effort is useful to understand the \emph{saddle point} structure of the Jacobian matrix. Indeed, 
Now, one can write
\begin{equation}
\label{J_saddle}
\mathsf{Jac}(\mathsf {X}^j; \bmu) =
\begin{bmatrix}
\mathsf A & \mathsf B^T \\
\mathsf B & 0 \\
\end{bmatrix},
\end{equation}
where
\begin{equation}
\mathsf A =
\begin{bmatrix}
\mathsf M_y + \mathsf D_{\mathsf y}( \mathsf {E}_{\textit{n}\ell}'[\mathsf y^j]^T)[\mathsf z^j] & 0\\
0 & \alpha \mathsf M_u & \\
\end{bmatrix}
\spazio \text{and} \spazio \mathsf B =
\begin{bmatrix}
\mathsf {E}_{\textit{n}\ell}'[\mathsf y^j] + \mathsf E_\ell & - \mathsf C
\end{bmatrix}.
\end{equation}
We remark that, % the assumption of at most quadratic nonlinearity in the state variable allows 
in the considered numerical settings, $\mathsf A$ is symmetric and, thus, we will always refer to \eqref{J_saddle} as a saddle point structure. To guarantee the solvability of the system we consider $\mathsf A$ an invertible matrix. Furthermore, we need the following \emph{Brezzi inf-sup condition} to be verified:
\begin{equation}
\label{FE_lbb}
\beta_\textit{Br} \disc(\bmu) \eqdot \adjustlimits \inf_{0 \neq \mathsf z} \sup_{0 \neq \mathsf x} \frac{\mathsf z^T\mathsf B \mathsf x}{\norm{\mathsf x}_{\state \times \control}\norm{\mathsf z}_{ \state}} \geq \hat{\beta}_\textit{Br} \disc > 0,
\end{equation}
where $\mathsf x = [{\mathsf y}^T, {\mathsf u}^T]^T$. %In the FE context, the inequality \eqref{FE_lbb} holds when the function spaces for state and adjoint coincide \cite{negri2015reduced,negri2013reduced}. 
The assumption $z \in \state$, will guarantee the fulfillment of the inf-sup stability condition in the FE approximation \cite{negri2015reduced,negri2013reduced}.
For general nonlinear problems, $\mathsf A$ may possibly be different from $\mathsf A^T$, however the well-posedness results can be extended, and the interested reader may refer to \cite{Benzi, GeneralizedSaddlePoint}. \\
%In the next Section we will present the peculiarities of the FE discretization of a bifurcation problem. Moreover, we will also introduce the spectral analysis to describe the branch-wise procedure we implemented in our algorithm for the reconstruction of the bifurcation diagram.
In the next Section, we will describe the branch-wise procedure implemented to reconstruct the bifurcation diagram.

\subsection{Bifurcation and stability analysis}
\label{bif}
%The study of the solution to a general nonlinear PDE($\bmu$) is usually very complicated, most of all of the form \eqref{ocp}, both from the theoretical and the numerical standpoint.
For nonlinear PDE($\bmu$)s, the assumptions (vii)-(ix) ensure the applicability of the well-known Implicit Function Theorem \cite{ciarlet2013linear, Prodi}.
Indeed, under those assumptions, one expects that when the parameter changes slightly, a stable solution evolves continuously in a unique manner.
When such conditions fail to be fulfilled, the model may undergo \textit{bifurcation} phenomena.
In particular, we consider the case in which the state equation has a singularity at some parameter value $\bmu^*$, and for the sake of clarity, we will assume $f = 0$ (or, equivalently, include the forcing term $f$ in the expression of $G$).
Indeed, from the mathematical perspective, we can give a precise definition of such points \cite{Prodi}.

\begin{definition}
\label{de:bifurcation_points}
A parameter value $\bmu^* \in \Cal P$ is a \textit{bifurcation point} for \eqref{eq:state} from the solution $y^* \eqdot y(\bmu^*)$, if there exists a sequence $(y_n, \bmu_n) \in \mbb Y \times \Cal P$, with $y_n \neq y^*$, such that (i) $G(y_n; \bmu_n) = 0$, and (ii) $(y_n, \bmu_n) \to (y^*, \bmu^*)$.
% \begin{itemize}
% \item[(i)] $G(y_n; \bmu_n) = 0$;
% \item[(ii)] $(y_n, \bmu_n) \to (y^*, \bmu^*)$.
% \end{itemize}
\end{definition}

Thus, the bifurcation phenomena is a paradigm for non-uniqueness in nonlinear analysis, and a necessary condition is the failure of the Implicit Function Theorem.
\begin{proposition}
\label{pr:ift}
A necessary condition for $\bmu^*$ to be a bifurcation point for $G$ is that the partial derivative $D_yG(y^*; \bmu^*)$ is not invertible.
\end{proposition}

Moreover, given the existence of many possible configurations for the system, a natural question is to understand which one inherits the stability of the unique solution, when it exists.
To perform a stability analysis, one of the most common and widely-studied method is the spectral analysis of the problem, which consists in the investigation of the eigenvalues of the system.
This technique, which has its roots in the ordinary differential equation (ODEs) theory \cite{seydel2009practical,kuznetsov2004elements,kielhofer2006bifurcation}, allows to understand the stability property of a solution to \eqref{eq:state} by means of the sign of the spectrum of the model operator.
In particular, for a general nonlinear uncontrolled problem, one linearizes the equation
around the solution under investigation, $\hat{y} = y(\hat{\bmu})$, and then solves the eigenvalue problem given by
\begin{equation}
\label{eq:eigen_state}
D_yG(\hat{y}; \hat{\bmu}) y_e = \rho_{ \hat{\bmu}} y_e ,
\end{equation}
where the pair $(\rho_{\hat \bmu}, y_e)$ represents respectively the eigenvalues and the eigenvector of $D_y G$ at $\hat{y}$ for each $\hat \bmu$ fixed.

This analysis provides us information about the physical stability of the problem. Indeed, the stability analysis is strongly bound to the investigation of the solution features after a small perturbation. If the perturbation is small enough and the dynamics of the system remains in a neighborhood of the solution, then it will be called a \textit{stable solution}. Thus, it is fundamental to observe that, in connection with ODEs stability theory, a positive eigenvalue gives an exponentially divergent behavior, while a negative one produces only small oscillations around the solution.
Therefore, it is clear that in order to have a stable solution, all eigenvalues must have negative real parts.%, in such a way that the divergent trend is prevented.

Dealing with the controlled problem \eqref{ocp} makes the analysis more involved. Indeed, the adjoint variable does not have physical meaning, making a straightforward application of the considerations above not applicable. Due to the high indefiniteness of the saddle point matrix \eqref{J_saddle}, a standard sign-analysis is no longer possible, see \cite{Benzi, BenziSimoncini, BenziWathen} as references on the topic.
Nevertheless, we can consider the eigenvalue problem for the system of the optimality conditions, in order to investigate \textit{a posteriori} the spectral property of \eqref{ocp}, as:
\begin{equation}
\label{eq:eigen_ocp}
D_X\Cal G(\hat{X}; \hat{\bmu}) X_e = \sigma_{\hat \bmu} X_e ,
\end{equation}
where $\hat{X} = X(\hat{\bmu})$ is the solution of which we are investigating the stability property and $(\sigma_{\hat \bmu}, X_e)$ is the eigenpair formed by the $\hat \bmu$-dependent eigenvalues $\sigma_{\hat \bmu}$ and eigenvectors $X_e$.
We will refer to \eqref{eq:eigen_state} as the \textit{state eigenvalue problem} and to \eqref{eq:eigen_ocp} as the \textit{global eigenvalue problem}.

As we said, the main issue with bifurcating system is the lack of the invertibility for the Fr\'echet derivative of the operator $\Cal G$ due to Proposition \ref{pr:ift}.
The latter, being equivalent to the injectivity and surjectivity of $\Cal G$, can be rewritten in terms of the \textit{continuous Babu{\v s}ka inf-sup stability}: there exists an inf-sup constant $\hat{\beta}_{\textit{Ba}} > 0$ such that
\begin{equation}
\label{eq:inf-sup_1}
\beta_{\textit{Ba}}(\bmu) = \adjustlimits \inf_{X \in \mbb X} \sup_{Y \in \mbb X} \frac{\langle D_X\Cal G[\hat{X}](X; \bmu), Y \rangle_{\mathbb X \mathbb X\dual}}{\norm{X}_{\mbb X}\norm{Y}_{\mbb X}} \geq \hat{\beta}_{\textit{Ba}} \qquad \forall \, \bmu \in \Cal P ,
\end{equation}
\begin{equation}
\label{eq:inf-sup_2}
\adjustlimits \inf_{Y \in \mbb X} \sup_{X \in \mbb X} \frac{\langle D_X\Cal G[\hat{X}](X; \bmu), Y \rangle_{\mathbb X \mathbb X\dual}}{\norm{X}_{\mbb X}\norm{Y}_{\mbb X}} > 0 \qquad \forall \, \bmu \in \Cal P .
\end{equation}
It is clear that the inclusion property $\mathbb{X}^\Cal N \subset \mathbb{X}$
is only a necessary but not sufficient condition for \eqref{eq:inf-sup_1} and \eqref{eq:inf-sup_2} to hold at the discrete level.
Hence, an additional assumption has to be required for the \textit{discrete Babu{\v s}ka inf-sup stability} of $\Cal G$: there exists a constant $\hat{\beta}_\textit{Ba} \disc > 0$ such that
\begin{equation}
\label{eq:inf_sup_disc}
\beta_\textit{Ba} \disc (\bmu) = \adjustlimits \inf_{\mathsf{X} \neq 0} \sup_{\mathsf{Y} \neq 0} \frac{\mathsf{Y}^T \mathsf{Jac}\ \mathsf{X}}{\norm{\mathsf{X}}_{\mbb X^{\Cal N}}\norm{\mathsf{Y}}_{\mbb X^{\Cal N}}} \geq \hat{\beta}_\textit{Ba} \disc \qquad \forall \, \bmu \in \Cal P .
\end{equation}

We remark that the continuous condition on the surjectivity \eqref{eq:inf-sup_2} is no longer needed for the discrete inf-sup stability \eqref{eq:inf_sup_disc}. In fact, while the assumption \eqref{eq:inf_sup_disc} corresponds to the non singularity of the matrix $\mathsf{Jac}$, a discrete counterpart of the assumption \eqref{eq:inf-sup_2} would require its surjectivity, or equivalently the injectivity of the transpose matrix $\mathsf{Jac}^T$, which being square would be the same as requiring \eqref{eq:inf_sup_disc}.
We also highlight that both continuous and discrete inf-sup conditions are satisfied as long as $\bmu \neq \bmu^*$.

We can finally present the branch-wise procedure we developed in order to deal with the numerical computation of multiple branches of solutions. Such an approach requires the combination of different methodologies for the approximation of the bifurcation diagram and the analysis of its stability properties.
In order to keep the presentation simple, we consider that the first component $\mu$ of the parameter $\bmu \in \Cal P \subset \mathbb R^P$ is the one responsible for the bifurcation behavior of the model, and in order to follow a single branch we consider that all the $P-1$ remaining parameters (thus the global physical/geometrical configuration) are unchanged on the branch.
\A{Even though in this work we will only deal with bifurcation phenomena with co-dimension one, the following methodology can be adapted to the multi-parameter and/or co-dimension $>1$ case, see e.g.\ \cite{pichi2021artificial}, by carefully choosing the technique to follow the branching behavior.}

Algorithm \ref{alg:01} summarizes how to reconstruct each branch $\mathcal{X}_i$ of solutions. More precisely, we combine, respectively:
\begin{itemize}
\item[{$\small{\circ}$}] Newton's method, as the nonlinear solver,
\item[{$\small{\circ}$}] Galerkin FE method, as the discretization phase,
\item[{$\small{\circ}$}] simple continuation method, as the bifurcation path tracer,
\item[{$\small{\circ}$}]generalized eigenvalue problems, as the stability detectors.
\end{itemize}
At the very beginning, one has to chose the branch to approximate, and the most preferable way to ``guide" the nonlinear solver to the desired configuration is through the initial guess.
Thus, in order to reconstruct a branch we consider the discrete version of the parameter space $\mathcal{P}_K = [\bmu_1, \dots, \bmu_K] \subset \mathcal{P}$ of cardinality $K$. We can take $\mathcal{P}_K$ as an ordered set, with the natural ordering induced by the first parameter component. Such ordering serves to assign the solution obtained for a given parameter $\bmu_{j-1}$ as the initial guess for the nonlinear solver at next iteration for $\bmu_{j}$. This allows us to follow the bifurcation behavior of the model. We choose the simplest variant of the continuation methods \cite{allogwer}, where the parametric set is fixed a priori, since it works well with pitchfork like bifurcation \cite{pichirozza,pichiquaini}. A more involved methodology have to be implemented when dealing with e.g.\ turning points or secondary bifurcations.

The next step is the actual discretization of the problem by means of the Newton-Kantorovich method \cite{ciarlet2013linear} combined with the Galerkin FE method. The initial guess for the former is set to the solution for the previous parameter value, while the latter projects the problem into a finite dimensional space, obtaining a linear system that we repeatedly solve until a convergence criterion is satisfied (here we chose a threshold tolerance $\epsilon$ for the norm of the global residual \eqref{G_compact_FE} at the $i-$th iteration of the Newton's method).
In Algorithm \ref{alg:01}, we call $\mathsf{Jac_y(\hat y, \hat \bmu)}$ the Jacobian matrix referred to the state equation \eqref{eq:state} and with $\mathsf V$ and $\mathsf {V_y}$ the scalar product matrices of the global optimization variable and of the state variable, respectively.
Finally, having computed a solution $\mathsf{X}_j$ of the problem \eqref{ocp} for the parameter $\bmu_j$, we can investigate its stability properties solving the two generalized eigenproblems, to recover the physical stability and the spectral properties, respectively for the state and global eigenvalue problems.

\begin{remark}
\label{re:branch}
Moreover, we highlight that the choice of the initial guess is fundamental, but it is not always sufficient to recover the full bifurcation diagram. In such cases, different techniques were proposed, e.g.\ manipulating the set $\mathcal{P}_K$ through predictor-corrector continuation methods, which involves pseudo-arclength strategy and homotopy \cite{allogwer,seydel2009practical}. Finally, when it is difficult to choose a proper guess, one can rely on:
\begin{itemize}
\item[{$\small{\circ}$}] the discretized version of analytic expressions, resembling the main properties of the sought solution \cite{pichiquaini};
\item[{$\small{\circ}$}] a deflation method, which requires only one initial guess, and discovers the full diagram preventing the convergence to already discovered solutions, helping the solver to find new branches \cite{pintore2019efficient, Charalampidis_et_al2018};
\item[{$\small{\circ}$}] the eigenvectors of the global eigenvalue problem, that have been used in \cite{pichirozza} to obtain the direction of the bifurcation branch in a neighborhood of the bifurcation points.
\end{itemize}
\end{remark}

\begin{algorithm}[H]
\caption{A pseudo-code for the reconstruction of a branch}\label{alg:01}
\begin{algorithmic}[1]
\State{$\mathsf{X}_0=\mathsf{X}_{guess}$}\Comment{Initial guess}
\For{$\bmu_j \in \mathcal{P}_K$}\Comment{Continuation loop}
\State{{$\mathsf{X}_j^{(0)} = \mathsf{X}_{j-1}$}} \Comment{{Continuation guess}}
\While{$|| \mathsf R(\mathsf{X}_j^{(i)}; \bmu_j)|| > \epsilon$}\Comment{Newton's method}
\State{$\mathsf{Jac}(\mathsf{X}_j^{(i)}; \bmu_j)\delta \mathsf{X} = \mathsf{R}(\mathsf{X}_j^{(i)}; \bmu_j)$}\Comment{Galerkin FE method}
\State{$\mathsf{X}_j^{(i+1)} = \mathsf{X}_j^{(i)} - \delta \mathsf{X}$}
\EndWhile
\State{$\mathsf{Jac_{y}}\mathsf{(y_j; \bmu_j) y_e = \rho_{\bmu_j} V_yy_e}$}\Comment{State eigenproblem}
\State{$\mathsf{Jac(X_j; \bmu_j) X_e = \sigma_{\bmu_j} \mathsf VX_e}$}\Comment{Global eigenproblem}
\EndFor
\end{algorithmic}
\end{algorithm}
\section{Bifurcations for Navier-Stokes Equations: the Coanda Effect}
\label{sec:state}
In this Section we analyze a bifurcating phenomenon deriving from Navier-Stokes equations in a sudden-expansion channel flow problem. Indeed, consider the channel geometry depicted in Figure \ref{fig:channel}. A fluid characterized by a high viscosity presents a jet which is symmetric w.r.t.\ the horizontal axis. Furthermore, a pair of vortices, called Moffatt eddies \cite{moffatt_1964}, form downstream of the expansion. Lowering the viscosity, the inertial effects of the fluid become more important and the two symmetric recirculation regions break the symmetry. Indeed, as the length of the recirculation zones increases, one can observe a non-uniform decrease of the pressure along the vertical axis. Thus, when we reach the aforementioned critical value, one recirculation zone expands whereas the other shrinks, giving rise to an asymmetric jet. This phenomenon is called the \textit{Coanda effect} and has been extensively studied in literature within different contexts \cite{tritton2012physical,AQpreprint,khamlich2021model,cardio,HESS2019379,pintore2019efficient,pichi2021artificial}.

From the mathematical point of view, this translates to a PDE($\bmu$) which, decreasing the viscosity $\mu$ below a certain critical value $\mu^{\ast}$, admits the existence of more solutions for the same value of $\mu \in \Cal P$.
During the study of the solution for different viscosity values, we expect the system to show two qualitatively different configurations:
\begin{itemize}
\item[{$\small{\circ}$}] a physically unstable configuration with a symmetric jet flow, the \emph{symmetric solution},
\item[{$\small{\circ}$}] a physically stable configuration with a wall-hugging jet, the \emph{asymmetric solution}.
\end{itemize}
These solutions, depicted in Figure \ref{fig:NS_sol_hf_bif}, coexist for parameter values below the critical one $\mu^{\ast}$ and belong to different branches that intersect in the bifurcation point, forming the so-called pitchfork bifurcation.

In the next Sections we introduce the mathematical formulation of Navier-Stokes equations describing the flow in a channel. This will serve us to highlight the bifurcated behavior of the system and its stability properties and it will be fundamental to understand how different controls affect the original system.

\subsection{Navier-Stokes problem as the state equation}
\label{sec:NS}
\begin{figure}[H]
\centering
\includegraphics[scale=0.4]{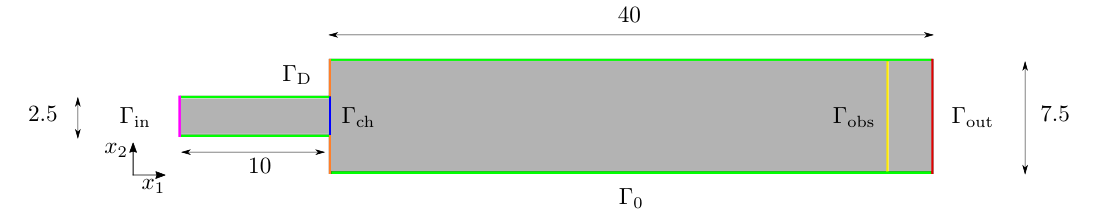}
\caption{\emph{Uncontrolled system}: domain $\Omega$ which represents a straight channel with a narrow inlet. }
\label{fig:channel}
\end{figure}

Here we consider a simplified setting with a two-dimensional planar straight channel with a narrow inlet and a sudden expansion, depicted in Figure \ref{fig:channel}, which represents a simplification of the left atrium and the mitral valve, respectively.
We define $\Gamma_{\text{in}} = \{0\}\times[2.5, 5]$ and $\Gamma_{\text{out}} = \{50\}\times[0, 7.5]$, where inflow and outflow boundary conditions are imposed, respectively. We indicate with $\Gamma_{\text{wall}}$, the boundaries representing the walls, in this case $\Gamma_{\text{wall}}= \Gamma_{\text{D}} \cup \Gamma_{0}$, where $\Gamma_{\text{D}} = \{\{10\}\times[0, 2.5]\}\cup \{\{10\}\times[5, 7.5]\}$ and
$ \Gamma_{0} = \partial \Omega \setminus \{\Gamma_{\text{in}} \cup \Gamma_{\text{D}} \cup \Gamma_{\text{out}}\}$.

The steady and incompressible Navier-Stokes equations for a viscous flow in $\Omega$ read as:
\begin{equation}
\label{eq:NS_eq}
\begin{cases}
-\mu \Delta v + v\cdot\nabla v + \nabla p=0 \quad &\text{in} \ \Omega, \\
\nabla \cdot v = 0 \quad &\text{in} \ \Omega, \\
v = v_{\text{in}} \quad &\text{on} \ \Gamma_{\text{in}}, \\
v = 0 \quad &\text{on} \ \Gamma_{\text{wall}}, \\
- pn + (\mu \nabla v) n = 0 \quad &\text{on} \ \Gamma_{\text{out}},

\end{cases}
\end{equation}
where $v = (v_{x_1}, v_{x_2})$ is the velocity of the fluid, $p$ is its pressure normalized over a constant density and $\mu$ represents the kinematic viscosity.
We supplement the system \eqref{eq:NS_eq} with proper boundary conditions: a stress free boundary condition on the velocity at the outlet, $\Gamma_{\text{out}}$ with outer normal $n$, a no-slip (homogeneous) Dirichlet boundary condition on $\Gamma_{\text{wall}}$ and a non-homogeneous Dirichlet boundary conditions $v_{\text{in}}$ at the inlet $\Gamma_{\text{in}}$ given by $v_{\text{in}}(x_2) = [20(5-x_2)(x_2 -2.5), 0]^T$.
% \begin{equation*}
% v_{\text{in}}(x_2) = \begin{bmatrix} 20(5-x_2)(x_2 -2.5) \\ 0 \end{bmatrix}.
% \end{equation*}
For later convenience, we introduce the dimensionless Reynolds number, which represents the ratio between inertial and viscous forces, and is given by $\text{Re} = Uh / \mu$, where $U$ and $h$ are characteristic velocity (i.e., maximum inlet velocity, $U = 31.25$) and characteristic length of the domain (i.e., length of the inlet section, $h = 2.5$), respectively. In the following we will consider $\mu$ as the parameter. %However, especially while commenting results, we will often refer to the associated Reynolds number.

Fixed the domain $\Omega$, the flow regime varies as we consider different values for the viscosity $\mu$ in $\mathcal{P} \subset \mathbb{R}$.
As we said in the introduction to this Section, this model exhibits a bifurcating behavior. Indeed we have the existence and uniqueness of the solution only above a certain critical value for the viscosity, that for this test case corresponds to $\mu^* \approx 0.96$. {Such value has been found in different works and numerical contexts for this benchmark \cite{pintore2019efficient,khamlich2021model,pichi2021artificial}. It can be obtained either ``a posteriori" by looking at the behaviour of the flow while varying the viscosity, or ``a priori" by investigating the change of sign of the leading eigenvalue w.r.t.\ the parameter. } %Thus, the problem for higher values of the Reynolds number loses the well-posedness and we have to refer to solution branches of a pitchfork bifurcation. 
To investigate the loss of uniqueness in a neighborhood of this pitchfork bifurcation, we set the parameter space as $\mathcal{P} = [0.5, 2.0]$, such that the first critical point $\mu^*$ is included.
These values for the viscosity correspond to Re in the interval [39.0, 156.0].

Let $\V=\left(H^1(\Omega)\right)^2$, $\V_{\text{in}}=\{v \in \V \mid v=v_{\text{in}} \text{ on }\Gamma_{\text{in}}, v=0 \text{ on }\Gamma_{\text{wall}}\}$, $\V_0=\{v \in \V \mid v=0 \text{ on }\Gamma_{\text{in}} \cup \Gamma_{\text{wall}}\}$ be the function spaces for velocity. Furthermore, let $\Q=L^2(\Omega)$ be the function space for pressure. The weak formulation of \eqref{eq:NS_eq} reads as: given $\mu \in \mathcal{P}$, find $v \in \V_{\text{in}}$ and $p \in \Q$ such that
\begin{equation}
\label{eq:gal_ns}
\left\{
\begin{aligned}
\mu\int_\Omega\nabla v\cdot\nabla \psi \, d\Omega + \int_\Omega \left(v\cdot\nabla v\right)\psi \, d\Omega - \int_\Omega p\nabla\cdot \psi \, d\Omega = 0 \quad\quad &\forall \, \psi \in \V_0, \\
\int_\Omega \pi\nabla\cdot v \, d\Omega = 0\quad\quad &\forall \, \pi \in \Q.
\end{aligned}
\right.
\end{equation}
We can rewrite the formulation of \eqref{eq:gal_ns} in an equivalent way
as: given $\mu \in \mathcal{P}$, find $v \in \V_{\text{in}}$ and $p \in \Q$ such that
\begin{equation}
\label{eq:gal_ns2}
\begin{cases}
a(v,\psi; \mu) +s(v,v,\psi) +b(\psi,p) = 0 \quad &\forall \, \psi \in \V_0, \\
b(v,\pi) = 0\quad &\forall \, \pi \in \Q ,
\end{cases}
\end{equation}
having introduced the following bilinear and trilinear forms for all $v$, $\bar v$, $\psi \in \V$ and $p \in \Q$,
\begin{equation}
\label{eq:forms}
a(v, \psi; \mu) =\mu\int_\Omega\nabla v\cdot\nabla \psi \, d\Omega, \qquad
b(v, p) = -\int_\Omega(\nabla\cdot v) \hspace{.05cm}p \, d\Omega, \qquad
s(v, \bar v, \psi)=\int_\Omega \left(v\cdot\nabla \bar v\right) \psi \, d\Omega. 
\end{equation}

\subsection{Numerical approximation of the problem}

We can now discuss the numerical approximation of the Navier-Stokes equation, that will be the state equation of the control problem in the next Sections. We consider a mesh on the domain $\Omega$ with $\Cal{N_T} =2785$ cells and $\Cal N_{y} = 24301$ degrees of freedom associated to a Taylor-Hood $\mathbb{P}^2$-$\mathbb{P}^1$ discretization of $\V \times \Q$. This choice is motivated by the well-known stability results of the Taylor-Hood Finite Element pair \cite{quarteroni2008numerical}.

\begin{figure}
\centering
\includegraphics[width=9cm]{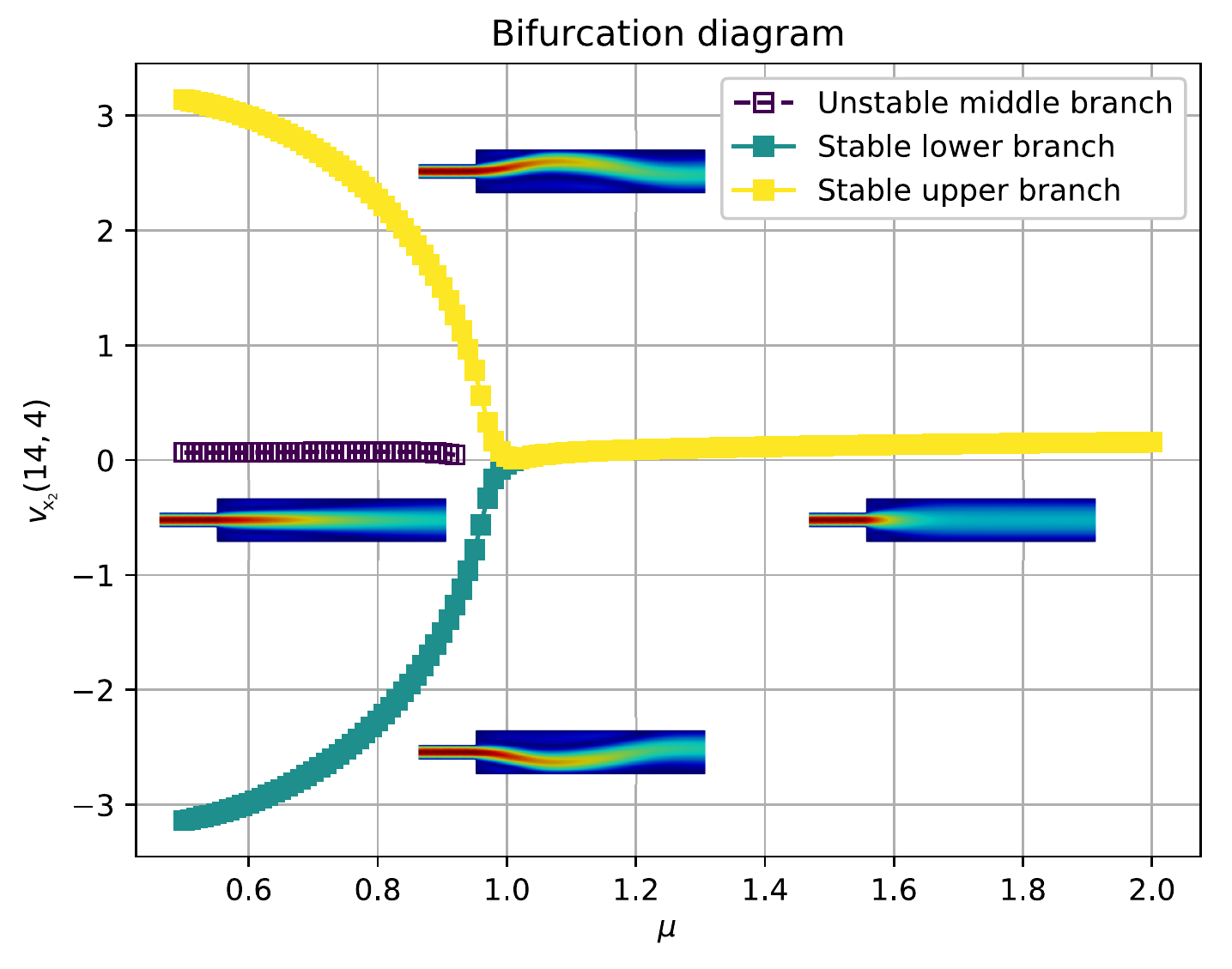}
\caption{\emph{Uncontrolled Navier-Stokes system}: Bifurcation diagram.}
\label{fig:bifurcation}
\end{figure}

In order to plot the bifurcation diagram, we choose an output value that results in a symmetry indicator of the approximated solution. \A{This function is given by the value of the vertical component of the velocity in a point of the channel, i.e.\ $v_{x_2}$ evaluated for $(x_1, x_2) = (14, 4)$, or the nearest node to that value since the mesh is unstructured. We remark that this output provides a merely graphical intuition about the symmetry breaking, and can be chosen to be efficiently computed by means of the analysis in Section \ref{sec_ROM}.}
In Figure \ref{fig:bifurcation} we plot the bifurcation diagram with all the solution branches found for the system \eqref{eq:gal_ns2} in the viscosity range chosen.
The numerical approximation clearly shows that a supercritical pitchfork bifurcation occurs around the critical viscosity value $\mu^* \approx 0.96$.
It is evident that we have a unique solution for all $\mu > \mu^*$, thus when the fluid behaves like a Stokes one, while we find three qualitatively different solutions increasing the Reynolds number. The bifurcation point $\mu^*$ is also the one responsible for the change in stability properties of the model. Indeed, the unique symmetric solution remains stable until it encounters the critical value $\mu^*$, where it becomes unstable. Moreover, this feature is inherited by the bifurcating solutions, which evolve as a physically stable branch.

\begin{figure}
\centering
\includegraphics[width=7cm]{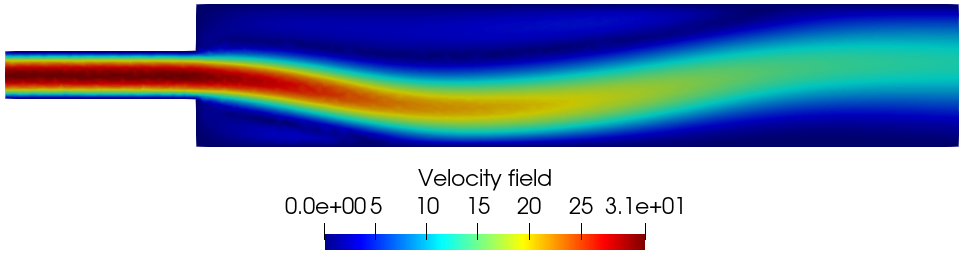}\qquad
\includegraphics[width=7cm]{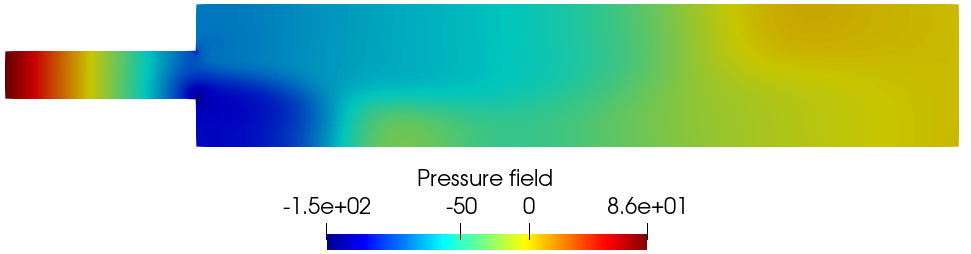}
\\
\vspace{1em}

\includegraphics[width=7cm]{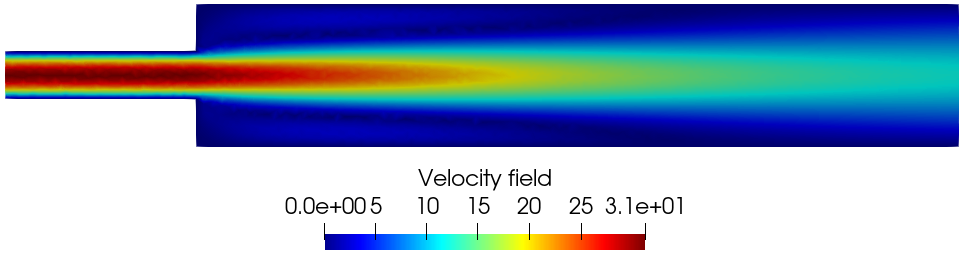}\qquad
\includegraphics[width=7cm]{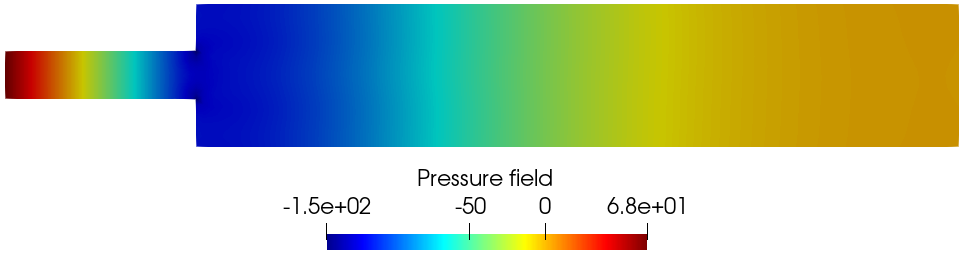}
\caption{\emph{Uncontrolled Navier-Stokes system}: representative solutions for $\mu = 0.5$, velocity and pressure fields, lower and middle branch, top and bottom, respectively.}
\label{fig:NS_sol_hf_bif}
\end{figure}

\begin{figure}
\centering
\includegraphics[width=7cm]{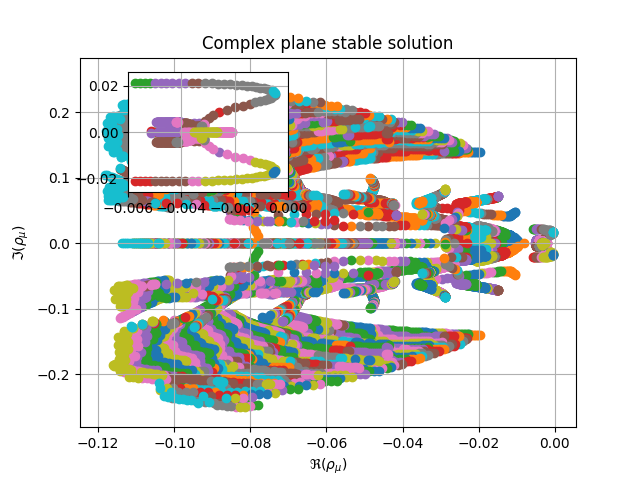}\qquad
\includegraphics[width=7cm]{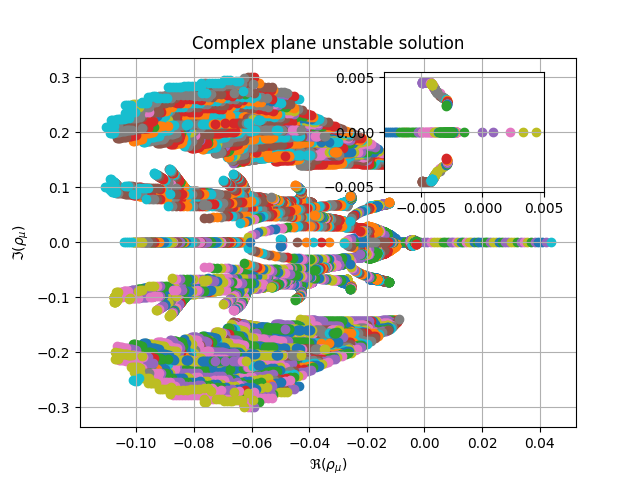}
\caption{\emph{Uncontrolled Navier-Stokes system}: eigenvalues of the state eigenproblem in the complex plane: stable and unstable solutions, left and right panels respectively.}
\label{fig:plot_eig_rc_nc}
\end{figure}

Some representative solutions for the lower and middle branch are presented in Figure \ref{fig:NS_sol_hf_bif}. Velocity and pressure fields belonging to different branches, for the same viscosity value $\mu = 0.5$, present qualitatively dissimilar behavior.
Indeed, the pressure for the lower branch decreases near the bottom-left corner of the expansion, causing the velocity to deflect, hugging the lower wall. Finally, thanks to the no-slip boundary condition the flux goes back to mid line, ending with a non-axis-symmetric outflow.

The stability analysis is performed through the eigenvalue analysis depicted in Section \ref{bif}, where Algorithm \ref{alg:01} has been applied exclusively to the state equation \eqref{eq:state}.
In particular, we analyzed the behavior of the first $N_{eig} = 100$ eigenvalues of \eqref{eq:eigen_state}, by means of the Krylov-Schur algorithm, varying the viscosity of the fluid.
Such eigenvalues are plotted in Figure \ref{fig:plot_eig_rc_nc} for the stable lower branch (left panel) and unstable middle branch (right panel). Note that since the Navier-Stokes operator is not symmetric, we have the presence of both real and complex eigenvalues. Here we are just interested in pitchfork bifurcation, thus in the zoom we follow the behavior of the biggest real eigenvalue and its sign. When investigating the stability of the lower branch, all the eigenvalues of the Navier-Stokes system, linearized around this stable solution, have negative real part. From the consideration of Section \ref{bif}, we can assert the stability of the wall-hugging branch. Indeed, the zoom in the left plot of Figure \ref{fig:plot_eig_rc_nc} shows no crossing of negative-real part eigenvalues. On the contrary, the close up in the right plot, which corresponds to the symmetric flow, shows the sign change of the biggest eigenvalue, thus characterizing a physically unstable solution.

\section{Steering bifurcating governing equations towards desired branches by means of Optimal Control}
\label{NS_ocp}
\begin{figure}
\centering
\begin{subfigure}[b]{0.49\textwidth}
\centering
\includegraphics[width=\textwidth]{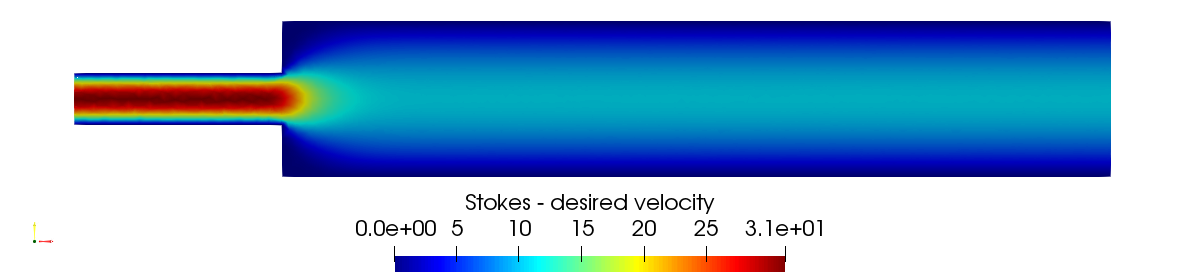}
\caption{}
\label{fig:vd_S}
\end{subfigure}
\hfill
\begin{subfigure}[b]{0.49\textwidth}
\centering
\includegraphics[width=\textwidth]{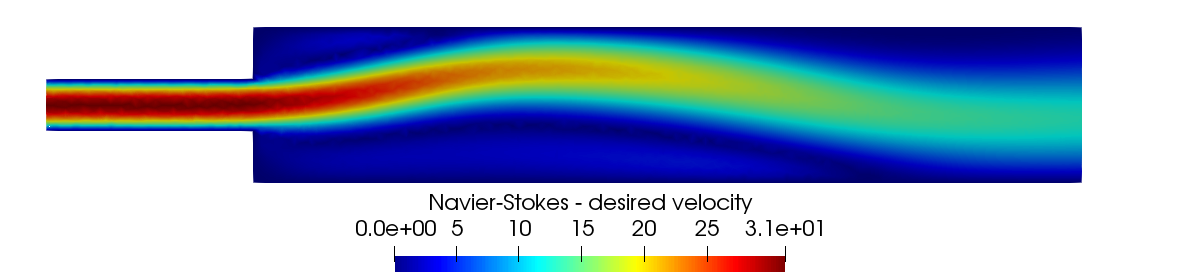}
\caption{}
\label{fig:vd_NS}
\end{subfigure}
\caption{\emph{Desired velocity profiles}: (a) symmetric profile obtained as Stokes solution for $\mu = 1$; (b) asymmetric profile given by the physically stable Navier-Stokes solution for $\mu = 0.49$.}
\label{fig:vd}
\end{figure}

In this Section we focus on several \ocp s governed by Navier-Stokes equations \eqref{eq:NS_eq} in the geometrical configuration of a contraction-expansion channel. We aim at understanding how different control problems can affect the solution behavior discussed in Section \ref{sec:state} for the uncontrolled case, especially when bifurcation phenomena are taken into account. %Indeed, the Coanda effect occurs and below the value $\mu^{*}$ three different solutions coexist. The multiplicity of configurations can be problematic in several applications: in the specific context of the mitral valve regurgitation, for example, a wall-hugging solution can represent an issue for accurate experimental measurements by echocardiography. 
This leads us to analyze the controlled systems, trying to reach state profiles which are different from the expected uncontrolled solution. Our goal is to investigate and better understand the role that optimal control plays as an attractor towards a desired configuration.\\
We thus follow the general procedure described in Section \ref{general_problem}, and discuss the specific case of optimal control for the Coanda effect. Nonetheless, the procedure adopted here is general and can be used in wide variety of applications.\\
For all the applications, we will simulate the physical phenomenon over the domain $\Omega$ shown in Figure \ref{fig:channel}.
Moreover, for the OCPs structure, we will require
the velocity solution $v \in \mathbb V$ to be the most similar to a desired profile $v_\text{d} \in \mathbb V_{\text{obs}} \eqdot (L^2(\Gamma_{\text{obs}}))^2$. The \emph{observation domain} $\Gamma_{\text{obs}}= \{47\}\times [0, 7.5]$ is a line near the end of the channel. This structure allows the control to change the solution at the outflow following a prescribed convenient configuration. During the rest of this work, we will employ two velocity solution profiles, which are showed in Figure \ref{fig:vd}: we will denote them as the \emph{symmetric desired profile (or target)} for Figure \ref{fig:vd_S} and the \emph{asymmetric desired profile (or target)} for Figure \ref{fig:vd_NS}. The first is the result of a Stokes system over $\Omega$ for $\mu = 1$ with the same boundary conditions of the Navier-Stokes uncontrolled equations \eqref{eq:NS_eq}. On the contrary, the latter is the physically stable solution of \eqref{eq:NS_eq} for $\mu = 0.49$. While the former choice aims at controlling the system towards a globally symmetric configuration with a weaker outgoing flux, the latter is set to achieve the opposite goal.

The purpose of steering the bifurcating behavior is summarized in the minimization of the functional
\begin{equation}
\label{eq:J_NS}
J_{\text{NS}}(v,u; v_\text{d}) = \half \norm{v - v_\text{d}}_ {\mathbb V_{\text{obs}} }^2+ \alf \norm{u}_{\control}^2,
\end{equation}
where $\control \eqdot (L^2(\Omega_u))^2$ with $\Omega_{u} \subset \overline \Omega$: indeed, the control action can be performed even over a portion of the boundary $\partial \Omega$. We will refer to $\Omega_u$ as the \emph{control domain}.
Within this study we will analyze how the choice of $\Omega_{u}$, combined with different values of the penalization parameter $\alpha$, affects the solution behavior of the system, compared to the uncontrolled Navier-Stokes state equation.
% In the next Section, we provide a general description of the \ocp s governed by Navier-Stokes equations. 
\C{We remark that hypothesis (i)-(ix) are verified in this specific \ocp\, setting, see e.g.\ \cite{hinze2008optimization}.}
Then, after a general introduction in Section \ref{sec:41}, we will provide the analysis of different optimal control systems, as follows.
\begin{enumerate}
\item[Section \ref{neumann}.]
\begin{enumerate}
A \emph{weak control} is built by controlling a Neumann boundary and the optimality system slightly affects the usual bifurcating nature of the uncontrolled Navier-Stokes equations.
\end{enumerate}
\item[Section \ref{distributed}.]
\begin{enumerate}
A \emph{strong control} effect can be observed over the classical bifurcating behavior of the uncontrolled solution by acting on the forcing term.
\end{enumerate}
\item [Section \ref{inlet}.]
\begin{enumerate}
The \emph{penalization parameter} $\alpha$ is analyzed while acting at the end of the inlet channel, and we discuss how changing $\alpha$ results in different orders of magnitude for the optimal control.
\end{enumerate}
\item [Section \ref{dirichlet}.]
\begin{enumerate}
We show how imposing different \emph{boundary flux} conditions completely changes the known behavior of the starting system.
\end{enumerate}
\end{enumerate}

Finally, in Section \ref{comparison}, remarks and comparisons on the spectral analysis of the four test cases are presented.

\subsection{\ocp s governed by Navier-Stokes equations}
\label{sec:41}
We recast \ocp s constrained to Navier-Stokes equations in the algebraic formulation presented in Section \ref{general_problem}. \\
The steady and incompressible controlled Navier-Stokes equations in a given domain $\Omega$ are:
\begin{equation}
\label{eq:OCP_NS_eq}
\begin{cases}
-\mu \Delta v + v\cdot\nabla v + \nabla p=C(u) \quad &\text{in} \ \Omega, \\
\nabla \cdot v = 0 \quad &\text{in} \ \Omega, \\
\end{cases}
\end{equation}
accompanied by some boundary conditions. The control operator $C : \mathbb U \rightarrow \V\dual$ can represent an external forcing term or a boundary term. If $C$ is defined in the whole domain we will say that the control is \emph{distributed}, while if it is defined in a portion of the internal domain, we will deal with \emph{localized control}. Furthermore, we will refer to \emph{Neumann control} and \emph{Dirichlet control}, if the control acts as Neumann or Dirichlet boundary conditions, respectively.
%Despite the variety of ways in which the control can act, the optimization system preserves a common structure that we are going to describe in the following. \\
The weak formulation of \eqref{eq:OCP_NS_eq} reads: given $\mu \in \mathcal{P}$, find $v \in \V$, $p \in \Q$ and $u \in \control$ such that
\begin{equation}
\label{eq:gal_ocp_ns2}
\begin{cases}
a(v,\psi; \mu) +s(v,v,\psi) +b(\psi,p) = c(u, \psi) \quad &\forall \, \psi \in \V, \\
b(v,\pi) = 0\quad &\forall \, \pi \in \Q ,
\end{cases}
\end{equation}
where $a (\cdot, \cdot; \mu)$, $b \cd $ and $s(\cdot, \cdot, \cdot)$
have been already defined in \eqref{eq:forms} while $c : \control \times \state \rightarrow \mathbb R$ is a bilinear form associated to the operator $C$.
First of all, to derive the optimality conditions, we need the adjoint variables $w \in \mathbb V$ and $q \in \mathbb Q$ for velocity and pressure, respectively. Let $X = ((v,p),u,(w,q)) \in \mathbb X \eqdot \mathbb Y \times \mathbb U \times \mathbb Y$ be an optimal solution, where $\state \eqdot \mathbb V \times \mathbb Q$. The Lagrangian functional for this specific problem is
\begin{equation}
\begin{aligned}
\label{lg_ocp_ns}
\Lg_{\text{NS}}(X; v_\text{d}, \bmu) = J_{\text{NS}}(v, u; v_\text{d}) \, + \,
\mu \int_\Omega\nabla v\cdot\nabla w \, & d\Omega \, +\,
\int_\Omega \left(v\cdot\nabla v\right)w \, d\Omega \\
& \, -\, \int_\Omega p\nabla\cdot w \, d\Omega
\,+\,
\int_\Omega q\nabla\cdot v \, d\Omega \,-\, c(u, w).
\end{aligned}
\end{equation}
The optimality system built through Fr\'echet differentiation is given by:
\begin{equation}
\label{KKT_NS}
\begin{cases}
D_{v}\Lg_{\text{NS}}(X; v_\text{d}, \bmu)[\varphi] = 0 & \forall \, \varphi \in \mathbb V,\\
D_{p}\Lg_{\text{NS}}(X; v_\text{d}, \bmu)[\xi] = 0 & \forall \, \xi \in \mathbb Q,\\
D_u\Lg_{\text{NS}}(X; v_\text{d}, \bmu)[\tau] = 0 & \forall \, \tau \in \control,\\
D_{w}\Lg_{\text{NS}}(X; v_\text{d}, \bmu)[\psi] = 0 & \forall \, \psi \in \mathbb V,\\
D_{q}\Lg_{\text{NS}}(X; v_\text{d}, \bmu)[\pi] = 0 & \forall \, \pi \in \mathbb Q,\\
\end{cases}
\end{equation}
where the first two equations form the \emph{adjoint system}, the last two form the \emph{state system} \eqref{eq:gal_ocp_ns2}, while differentiating w.r.t.\ the variable $u$ leads to the \emph{optimality equation}.
In particular, the adjoint system combined with the optimality equation has the following form:
\begin{equation}
\label{eq:gal_adj_ocp_ns2}
\begin{cases}
m(v, \varphi) + a(w,\varphi; \mu) +s(\varphi, v, w) +s(v,\varphi, w) +b(\varphi,q) = m(v_\text{d}, \varphi) \quad &\forall \, \varphi \in \V, \\
b(w,\xi) = 0\quad &\forall \, \xi \in \Q ,\\
\alpha r(u, \tau) = c(\tau, w) \quad & \forall \, \tau \in \control,
\end{cases}
\end{equation}
%while the optimality equation is given by
%\begin{equation}
%\label{eq:gal_opt_ocp_ns2}
%\alpha r(u, \tau) = c(\tau, w) \quad \forall \, \tau \in \control,
%\end{equation}
where $m\goesto{\V}{\V}{\mathbb R}$ and $r\goesto{\control}{\control}{\mathbb R}$ terms come from
the Fr\'echet derivative of \eqref{eq:J_NS} w.r.t.\ the velocity and control, respectively. They represent the $L^2$ scalar product in
$\Gamma_{\text{obs}}$ and $\Omega_{u}$, respectively. Furthermore, we remark that $s(\varphi,v, w) + s(v,\varphi, w)$ is the linearization around $v$ of the trilinear form $s(v,v,\varphi)$, by definition. Therefore, the strong formulation for \eqref{eq:gal_adj_ocp_ns2}
%and \eqref{eq:gal_opt_ocp_ns2} 
reads:
\begin{equation}
\label{eq:OCP_ADJ_OPT_NS_eq}
\begin{cases}
v\mathbb{I}_{\Omega_{\text{obs}}} -\mu \Delta w - v\cdot\nabla w + (\nabla v)^T w + \nabla q= v_\text{d} \mathbb{I}_{\Omega_{\text{obs}}} \quad &\text{in} \ \Omega, \\
\nabla \cdot w = 0 \quad &\text{in} \ \Omega, \\
\alpha u \mathbb{I}_{\Omega_u} = C^* w \quad &\text{in} \ \Omega, \\
\end{cases}
\end{equation}
where $\mathbb{I}_{\Omega_u}$ and $\mathbb{I}_{\Omega_{\text{obs}}}$ are the indicator functions of the control and observation domains, respectively. \A{The detailed derivation of the optimality system is addressed in several works, see e.g.\ \cite{hinze2008optimization,Fursikov1998852,Gunzburger2000249}}. The global optimization problem reads: given $\mu \in \Cal P$, find $X = ( (v, p),u, (w, q)) \in \mathbb X$ such that \eqref{eq:OCP_NS_eq} and \eqref{eq:OCP_ADJ_OPT_NS_eq} are verified.
\\We remark that, if we call $y \eqdot (v, p)$ and $z \eqdot (w, q)$, we recover the global algebraic formulation presented in Section \ref{general_problem} and the saddle point structure is preserved. Indeed,
let us suppose that we apply the Taylor-Hood approximation $\mathbb{P}^2$-$\mathbb{P}^1$ for state $y$ and adjoint variable $z$. Furthermore, we discretize the space $\control$ with FE using $\mathbb P^2$ polynomials.
Recalling the notation of Section \ref{FE}, we define the quantities
\begin{equation}
\mathsf y =
\begin{bmatrix}
\mathsf v \\
\mathsf p
\end{bmatrix} , \qquad
\mathsf z =
\begin{bmatrix}
\mathsf w \\
\mathsf q
\end{bmatrix}, \qquad
\mathsf M_{y} =
\begin{bmatrix}
\mathsf M_v & 0 \\
0 & 0
\end{bmatrix}, \quad \text{and} \quad
\mathsf C =
\begin{bmatrix}
\mathsf C_v \\
0
\end{bmatrix},
\end{equation}
where $\mathsf v, \mathsf p, \mathsf w, \mathsf q$ are the column vectors of FE coefficients for state and adjoint, velocities and pressures respectively, while
$\mathsf M_{v}$ is the mass velocity matrix and $\mathsf C_v$ derives by the bilinear form $c\cd$.
Furthermore, the linearized state equation structure can be now expressed as
\begin{equation}
\label{eq:NS_matrix}
\mathsf {E}_{\textit{n}\ell}'[\mathsf y^j] + \mathsf E_\ell =
\begin{bmatrix}
\mathsf S[\mathsf v^j] & 0 \\
0 & 0
\end{bmatrix} +
\begin{bmatrix}
\mathsf K & \mathsf D^T \\
\mathsf D & 0
\end{bmatrix} =
\begin{bmatrix}
\mathsf K + \mathsf S[\mathsf v^j]& \mathsf D^T \\
\mathsf D & 0
\end{bmatrix}
,
\end{equation}
where $\mathsf K$ is the stiffness matrix associated to the bilinear form $a(\cdot, \cdot; \mu)$, $\mathsf D$ is the continuity equation matrix coming from $b \cd$ and $ \mathsf S[\mathsf v^j]$ is the algebraic formulation of $s(v, \cdot, \cdot) + s(\cdot, v, \cdot)$ evaluated at the FE velocity basis functions. It remains to understand the specific structure of $\mathsf D_{\mathsf y}( \mathsf {E}_{\textit{n}\ell}'[\mathsf y]^T)[\mathsf z^j]$ defined in \eqref{J_ocp}. To this end, we define $s_{\text{ad}}(v, w, \varphi)$ as the \emph{adjoint operator} of the linearized trilinear form $s(v, v, \cdot)$ around the state velocity $v$. Applying $s_{\text{ad}}(v, \cdot, \cdot)$ to the basis functions of $\mathbb V^{\Cal N_v}$
will result in $ \mathsf S[\mathsf v^j]^T$. In the Jacobian matrix evaluation, a linearization of $s_{\text{ad}}(w,v, \varphi)$ is performed not only in $w$, but also w.r.t.\ the variable $v$. This process will lead to
\begin{equation}
\mathsf D_{\mathsf y^j}( \mathsf {E}_{\textit{n}\ell}'[\mathsf y^j]^T)[\mathsf z^j] =
\begin{bmatrix}
\mathsf D_{\mathsf v}(\mathsf S[\mathsf v^j]^T)([\mathsf w^j]) & 0 \\
0 & 0 \\
\end{bmatrix},
\end{equation}
where $\mathsf D_{\mathsf v}(\mathsf S[\mathsf v^j]^T)([\mathsf w^j]) $ is given by the form $ s_{\text{ad}} (w, \cdot, \cdot)$ applied to the FE velocity basis. Then, the Jacobian reads
\begin{equation}
\label{eq:J_ocp_NS}
\mathsf{Jac}_{\text{NS}}(\mathsf X^j; \bmu) =
\begin{bmatrix}
\mathsf M_v + \mathsf D_{\mathsf v}(\mathsf S[\mathsf v^j]^T)[\mathsf w^j] & 0 & 0 & \mathsf K^T + \mathsf S[\mathsf v^j]^T & \mathsf D^T \\
0 & 0 & 0 & \mathsf D & 0 \\
0 & 0 & \alpha \mathsf M_u & - \mathsf C^T_v & 0 \\
\mathsf K + \mathsf S[\mathsf v^j] & \mathsf D^T & - \mathsf C_v & 0 & 0 \\
\mathsf D & 0 & 0& 0 & 0 \\
\end{bmatrix}
=
\begin{bmatrix}
\mathsf A & \mathsf B^T \\
\mathsf B & 0 \\
\end{bmatrix}
\end{equation}
where $\mathsf X$ is the FE coefficient vector of the optimal solution and
\begin{equation}
\mathsf A =
\begin{bmatrix}
\mathsf M_v + \mathsf D_{\mathsf v}(\mathsf S[\mathsf v^j]^T)([\mathsf w^j]) & 0 & 0 \\
0 & 0 & 0 \\
0 & 0 & \alpha \mathsf M_u \\
\end{bmatrix}
\quad \text{and} \quad \mathsf B =
\begin{bmatrix}
\mathsf K + \mathsf S[\mathsf v^j] & \mathsf D^T & - \mathsf C_v \\
\mathsf D & 0 & 0 \\
\end{bmatrix}.
\end{equation}
As already specified in Section \ref{FE}, we assume that for $\mu \neq \mu^*$ the saddle point \eqref{eq:J_ocp_NS} is well-posed. Moreover, we highlight that we are dealing with a \emph{nested saddle point} structure: indeed, for the state equation \eqref{eq:NS_matrix} we require that, for a given $\mu \neq \mu^*$ and fixed $\mathsf v^j$,
the matrix $\mathsf K + \mathsf S[\mathsf v^j] $ is invertible and that Brezzi inf-sup condition holds, i.e.
\begin{equation}
\label{NS_FE_lbb}
\beta_{\text{Br}, \text{NS}}\disc \eqdot \adjustlimits\inf_{\mathsf p \neq 0} \sup_{\mathsf v \neq 0} \frac{\mathsf p^T\mathsf D \mathsf v}{\norm{\mathsf v}_{\V}\norm{\mathsf p}_{ \Q}} \geq \overline \beta_{\text{Br}, \text{NS}} \disc > 0.
\end{equation}
This is indeed the case for the Taylor-Hood discretization introduced in Section \ref{sec:state}.

In the next subsections we analyze how the controlled problem behaves, comparing its properties with the ones of the uncontrolled system presented in Section \ref{sec:state}.
We will use the word \emph{natural optimal branch} to describe the branch that is obtained by running Algorithm \ref{alg:01} with a trivial initial guess. This branch may consist of either symmetric or asymmetric configurations, depending on the test case. Further branches may exist, but are much harder to compute in practice and require very tailored initial guesses that can be provided by running Algorithm \ref{alg:01} in a neighborhood of $\mu^*$, and will be named \emph{non-natural optimal branches}\footnote{For the sake of exposition, each branch is extended to $\mu > \mu^*$ with the unique solution.}. We interpret the concept of natural optimality as a \emph{numerical stability} property of the optimal control system. Indeed, as already specified in Section \ref{general_problem}, for \ocp s it makes no sense to talk about the physical stability of the global optimal solution. In fact, the system is ``artificially" built by adding non-physical adjoint variables, with the aim of changing the system behavior.

\subsection{Neumann Control: weak steering}
\label{neumann}
\begin{figure}
\centering
\begin{subfigure}[b]{0.49\textwidth}
\centering
\hspace*{-2mm}\includegraphics[width=\textwidth]{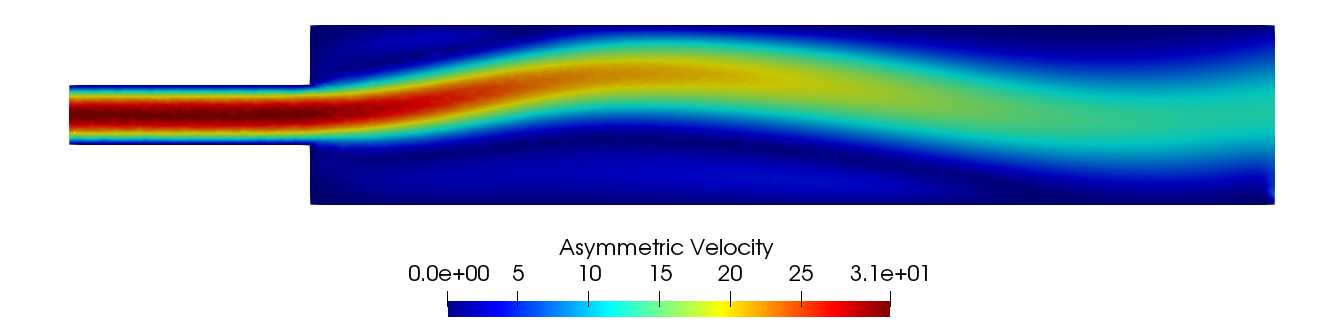}
\caption{}
\label{fig:v_as}
\end{subfigure}
\hfill
\begin{subfigure}[b]{0.49\textwidth}
\centering
\hspace*{-6mm}\includegraphics[width=\textwidth]{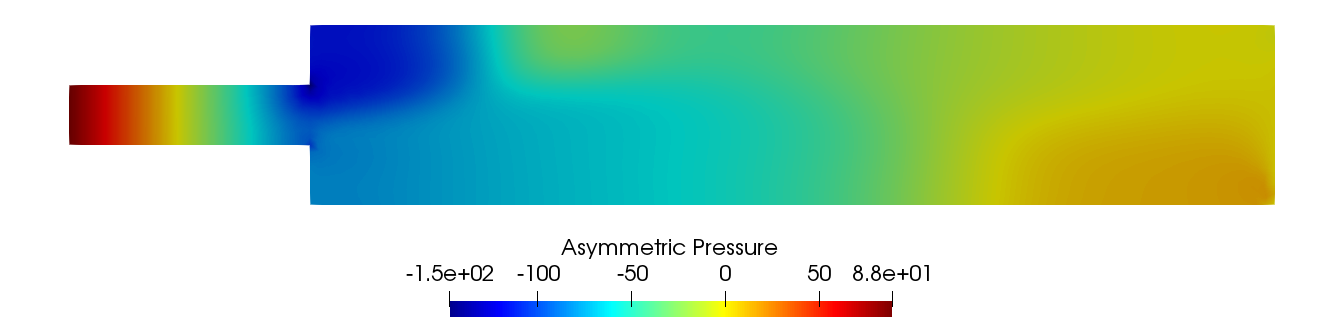}
\caption{}
\label{fig:p_as}
\end{subfigure}
\begin{subfigure}[b]{0.49\textwidth}
\centering
\includegraphics[width=0.98\textwidth]{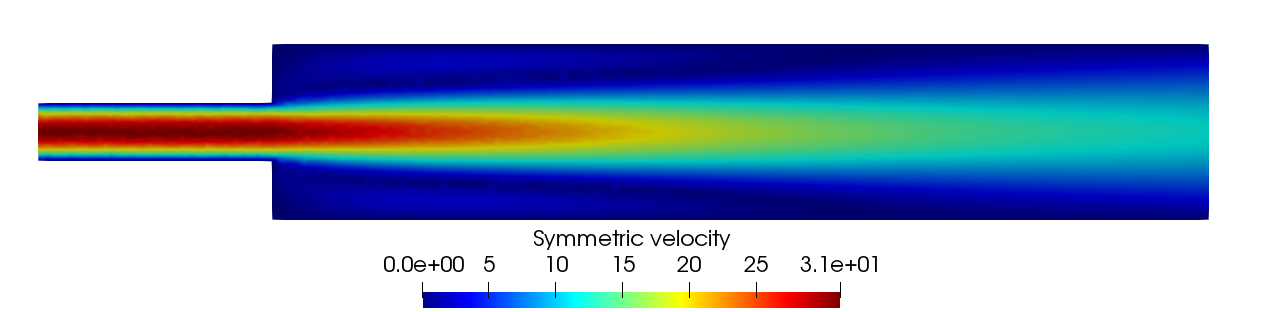}
\caption{}
\label{fig:v_s}
\end{subfigure}
\hfill
\begin{subfigure}[b]{0.49\textwidth}
\centering
\hspace{3mm}\includegraphics[width=0.98\textwidth]{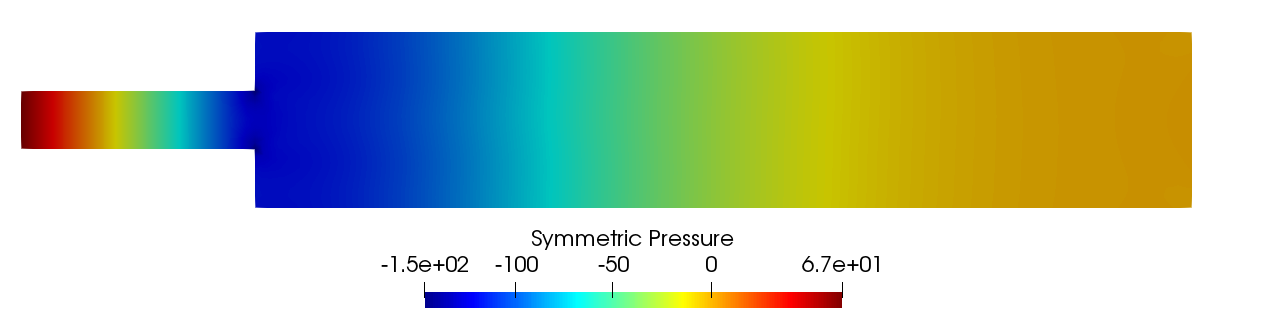}
\caption{}
\label{fig:p_s}
\end{subfigure}
\hfill
\caption{\emph{Neumann Control}: optimal solutions with $\alpha = 0.01$ and $\mu=0.5$, belonging to the natural optimal (panels (a) and (b) for state velocity and pressure, respectively) and the non-natural optimal (panels (c) and (d)) branches.}
\label{fig:Neumann_solution}
\end{figure}

\begin{figure}
\centering
\includegraphics[width=0.489\textwidth]{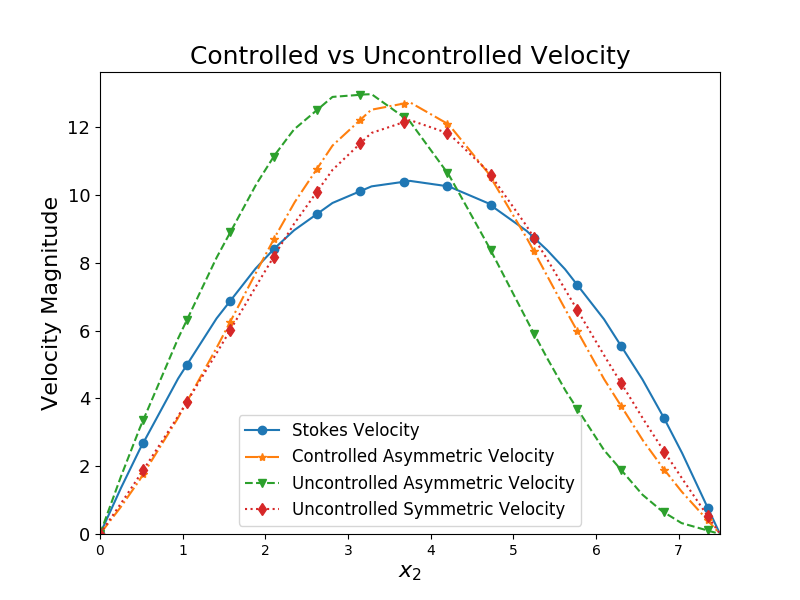}
\includegraphics[width=0.46\textwidth]{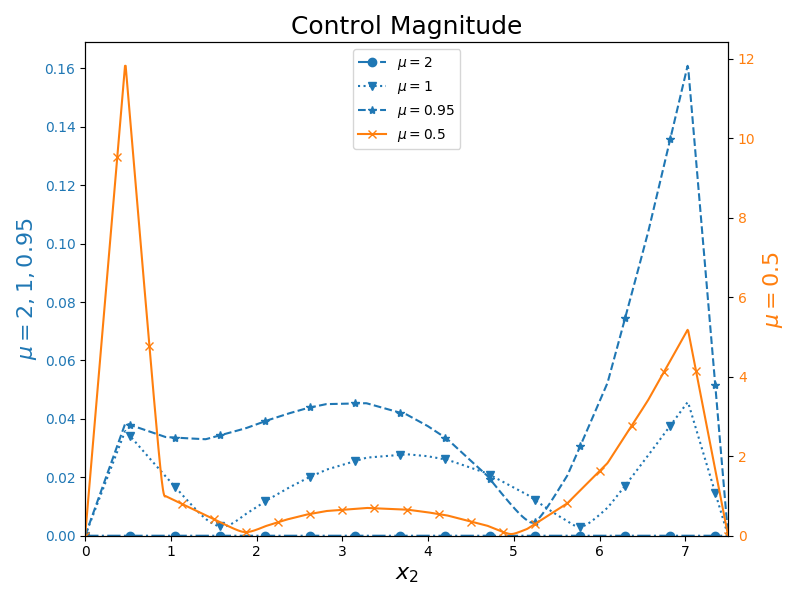}
\caption{\emph{Neumann Control}. \emph{Left}: comparison of velocity profiles in the controlled and uncontrolled cases for $\alpha = 0.01$, $\mu = 0.5$ on $\Gamma_{\text{obs}}$ w.r.t.\ the desired profile when following the natural optimal branch. \emph{Right}: representation of control variable evolution for $\alpha = 0.01$, $\mu=2, 1, 0.95, 0.5$ over $\Gamma_{\text{out}}$ when following the natural optimal branch.}
\label{fig:mag}
\end{figure}

\begin{figure}
\centering
\includegraphics[width=0.489\textwidth]{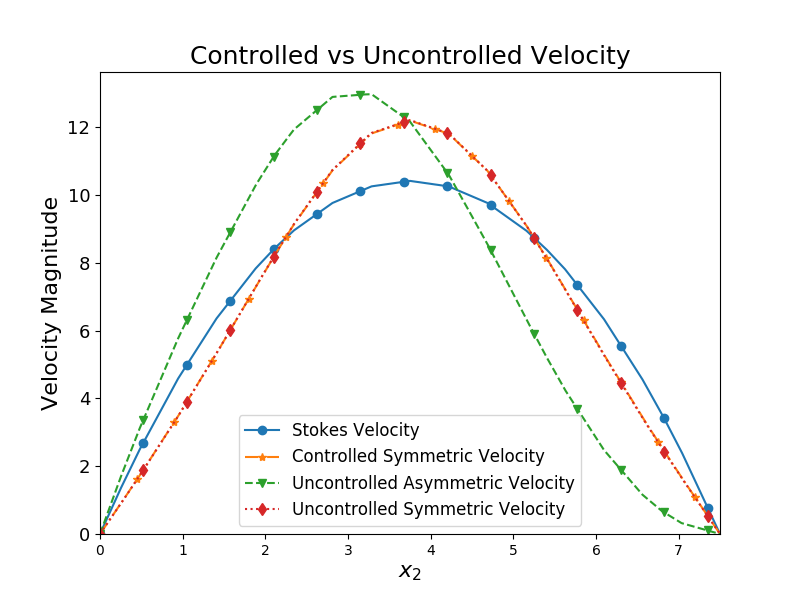}
\includegraphics[width=0.46\textwidth]{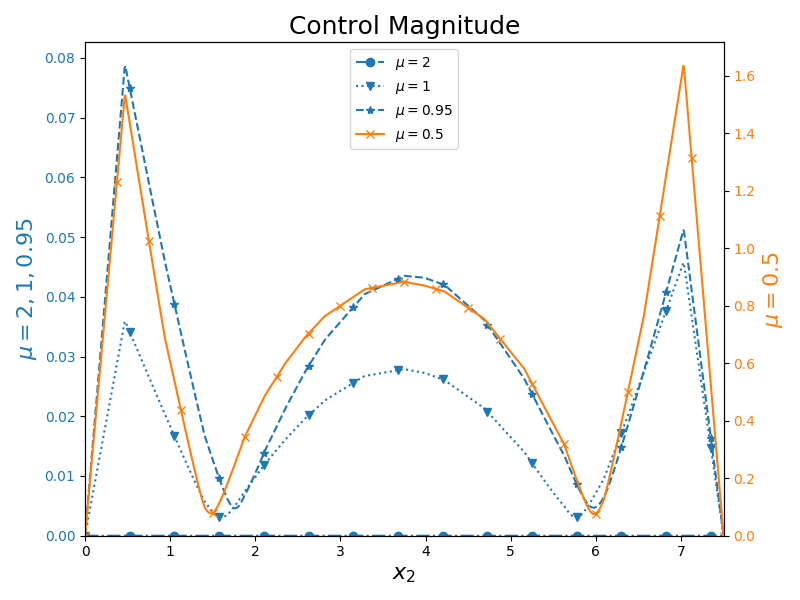}
\caption{\emph{Neumann Control}. \emph{Left}: comparison of velocity profiles in the controlled and uncontrolled cases for $\alpha = 0.01$, $\mu = 0.5$ on $\Gamma_{\text{obs}}$ w.r.t.\ the desired profile when following the non-natural optimal branch. The lines marked by ``Controlled Symmetric Velocity'' and ``Uncontrolled Symmetric Velocity'' overlap. \emph{Right}: representation of control variable evolution for  $\alpha = 0.01$, $\mu=2, 1, 0.95, 0.5$ over $\Gamma_{\text{out}}$ when following the non-natural optimal branch.}
\label{fig:s_mag}
\end{figure}

The first test case we present is a Neumann control over the boundary $\Gamma_{\text{out}}$, where homogeneous Dirichlet conditions are applied to $\Gamma_{\text{wall}} \eqdot \Gamma_0 \cup \Gamma_{\text{D}}$.
More specifically, in this case, the optimality conditions read: given $\mu \in \Cal P$ find $X \in \mathbb X$ such that
\begin{equation}
\label{eq:Neumann_eq}
\begin{cases}
v\mathbb{I}_{\Gamma_{\text{obs}}} -\mu \Delta w - v\cdot\nabla w + (\nabla v)^T w + \nabla q= v_\text{d} \mathbb{I}_{\Gamma_{\text{obs}}} \quad &\text{in} \ \Omega, \\
\nabla \cdot w = 0 \quad &\text{in} \ \Omega, \\
w =0 \quad &\text{on} \ \Gamma_{\text{in}} \cup \Gamma_{\text{wall}}, \\
- qn + (\mu \nabla w) n = 0 \quad &\text{on} \ \Gamma_{\text{out}}, \\
\alpha u \mathbb{I}_{\Gamma_{\text{out}}} = w\mathbb{I}_{\Gamma_{\text{out}}} \quad &\text{in} \ \Omega, \\
-\mu \Delta v + v\cdot\nabla v + \nabla p=0 \quad &\text{in} \ \Omega, \\
\nabla \cdot v = 0 \quad &\text{in} \ \Omega, \\
v = v_{\text{in}} \quad &\text{on} \ \Gamma_{\text{in}}, \\
v = 0 \quad &\text{on} \ \Gamma_{\text{wall}}, \\
- pn + (\mu \nabla v) n = u \quad &\text{on} \ \Gamma_{\text{out}}.\\
\end{cases}
\end{equation}
The desired velocity $v_\text{d}$ will always be of the symmetric type for this specific example. In other words, we are studying which is the best choice for Neumann boundary condition, to reach the exiting symmetric profile shown in Figure \ref{fig:vd_S}. We study the behavior of the controlled solution varying $\alpha = 1, 0.1, 0.001, 0.0001$, where the greater is the value of $\alpha$ the lower is the strength of the control.
In Figure \ref{fig:Neumann_solution} we show some representative solutions for $\alpha = 0.01$ and $\mu = 0.5$, for state velocity and pressure variables. In this case, the natural optimal branch is composed by asymmetric solutions (Figure \ref{fig:Neumann_solution}, top), while there is a further non-natural optimal branch made up by symmetric solutions (Figure \ref{fig:Neumann_solution}, bottom). Results obtained following the natural optimal and non-natural optimal branches are shown in Figures \ref{fig:mag} and \ref{fig:s_mag}, respectively.
Therefore, we conclude that the Neumann control affects \emph{weakly} the system, as it is not able to steer the system towards the desired symmetric configuration after bifurcation has occurred, thus not changing drastically the features already observed for the uncontrolled state equations (see Figure \ref{fig:bifurcation}).\\
The left plot of Figure \ref{fig:mag} depicts the velocity profile magnitude over $\Gamma_{\text{obs}}$ for the highest value of the Reynolds number when following the natural optimal branch. Even though the obtained velocity (marked by an orange line) is indeed different from the desired profile (denoted by a blue line), especially for what concerns peak values, we observe that the Neumann control straightens the flux near the end of the channel (compare the orange line to the green line, which represents the uncontrolled asymmetric profile), even when high Reynolds numbers are considered. The resulting profile is similar to the uncontrolled symmetric velocity (red line), even though full symmetry is not achieved. The action of the control variable is shown in the right plot of Figure \ref{fig:mag} when changing the parameter $\mu$ following the natural optimal branch: the control is stronger for $\mu < \mu^*$ (i.e., when the wall-hugging phenomenon occurs and straightening in necessary), while it remains low in magnitude for $\mu > \mu^*$.\\
Similarly, the left plot of Figure \ref{fig:s_mag} shows the velocity profile magnitude over $\Gamma_{\text{obs}}$ for the highest value of the Reynolds number when following the non-natural optimal branch. In this case, the controlled symmetric profile (orange line) coincides with the uncontrolled symmetric profile (red line). Furthermore, the right plot of Figure \ref{fig:s_mag} shows that the control variable around the critical $\mu^*$ (e.g., $\mu = 1$ and $\mu = 0.95$) is asymmetric to counteract the stable wall-hugging physically driven behavior of the uncontrolled system. We further remark that, compared to the natural optimal branch, the control variable of the non-natural optimal branch is much lower in magnitude.\\
Table \ref{Neumann_J} shows the value of the cost functional \eqref{eq:J_NS} for several values of $\mu$ (rows) and $\alpha$ (columns), following either the natural optimal or non-natural optimal branches. The first column also shows the value of the \emph{uncontrolled functional}, i.e. \eqref{eq:J_NS} evaluated for the uncontrolled velocity $v$ of the equation \eqref{eq:NS_eq} and zero control.
The main observation is that decreasing the value of $\alpha$ results in lower cost functional values, since a lower value of $\alpha$ allows stronger control to take place and drive the velocity to the desired configuration. In all cases, the non-natural branch presents lower values of the functional compared to the natural branch; this has to be expected, as the cost functional measures deviation from a symmetric target, and the non-natural branch is clearly closer to the target being made of symmetric solution (compare e.g.\ for $\mu = 0.05$ and $\alpha = 0.001$ the left panels of Figures \ref{fig:mag}-\ref{fig:s_mag}). However, the natural branch is the one for which the control procedure is influencing the most the cost functional values: for instance, for $\mu = 0.05$ and $\alpha = 0.001$, the cost functional is decreased by $6\%$ on the non-natural branch and by $55\%$ on the natural one w.r.t.\ the corresponding uncontrolled configuration. Again, this has to be expected from the previous discussion of Figures \ref{fig:mag}-\ref{fig:s_mag}, which shows larger impact of the control procedure to straighten the solution on the natural branch. Finally, large values of $\mu$ have negligible cost functionals, as the target velocity almost coincides with the uncontrolled velocity.
From such an analysis we deduce that, when bifurcating phenomena occur, a configuration can perform better than another one, and finding all the solution branches can be of great importance to understand the solution that best recover the desired profile.


\begin{table}
\caption{\emph{Neumann Control}: comparison of the functional value w.r.t. stable and unstable uncontrolled solutions. (Nat.) Natural optimal branch. (n-Nat.) Non-natural optimal branch.}
\label{Neumann_J}
\begin{center}
\tabcolsep=0.11cm
\footnotesize{
\begin{tabular}{|c||c|c||c|c||c|c||c|c||c|c||}
\hline
& Stable & \cellcolor[HTML]{E5E3E3}Unstable & Nat. & \cellcolor[HTML]{E5E3E3}n-Nat. & Nat. & \cellcolor[HTML]{E5E3E3}n-Nat. & Nat. & \cellcolor[HTML]{E5E3E3}n-Nat. & Nat. & \cellcolor[HTML]{E5E3E3}n-Nat. \\ \cline{2-11}
\multirow{-2}{*}{$\mu$} & \multicolumn{2}{c||}{Uncontrolled} & \multicolumn{2}{c||}{$\alpha = 1$} & \multicolumn{2}{c||}{$\alpha = 0.1$} & \multicolumn{2}{c||}{$\alpha = 0.01$} & \multicolumn{2}{c||}{$\alpha = 0.001$} \\ \hline
$2$ & 5.14e--9 & \cellcolor[HTML]{E5E3E3}5.14e--9 & 5.13e--9 & \cellcolor[HTML]{E5E3E3}5.13e--9 & 5.13e--9 & \cellcolor[HTML]{E5E3E3}5.13e--9 & 5.13e--9 & \cellcolor[HTML]{E5E3E3}5.13e--9 & 5.07e--9 & \cellcolor[HTML]{E5E3E3}5.07e--9 \\ \hline
$1.5$ & 4.38e--6 & \cellcolor[HTML]{E5E3E3}4.38e--6 & 4.38e--6 & \cellcolor[HTML]{E5E3E3}4.38e--6 & 4.38e--6 & \cellcolor[HTML]{E5E3E3}4.38e--6 & 4.37e--6 & \cellcolor[HTML]{E5E3E3}4.37e--6 & 4.28e--6 & \cellcolor[HTML]{E5E3E3}4.28e--6 \\ \hline
$1$ & 4.10e--3 & \cellcolor[HTML]{E5E3E3}4.10e--3 & 4.10e--3 & \cellcolor[HTML]{E5E3E3}4.10e--3 & 4.10e--3 & \cellcolor[HTML]{E5E3E3}4.10e--3 & 4.08e--3 & \cellcolor[HTML]{E5E3E3}4.10e--3 & 3.92e--3 & \cellcolor[HTML]{E5E3E3}3.92e--3 \\ \hline
$0.9$ & 3.33e--2 & \cellcolor[HTML]{E5E3E3}1.63e--2 & 3.33e--2 & \cellcolor[HTML]{E5E3E3}1.63e--2 & 3.30e--2 & \cellcolor[HTML]{E5E3E3}1.63e--2 & 3.15e--2 & \cellcolor[HTML]{E5E3E3}1.63e--2 & 2.93e--2 & \cellcolor[HTML]{E5E3E3}1.55e--2 \\ \hline
$0.8$ & 2.08e--1 & \cellcolor[HTML]{E5E3E3}6.52e--2 & 2.07e--1 & \cellcolor[HTML]{E5E3E3}6.52e--2 & 2.04e--1 & \cellcolor[HTML]{E5E3E3}6.51e--2 & 1.88e--1 & \cellcolor[HTML]{E5E3E3}6.51e--2 & 1.70e--1 & \cellcolor[HTML]{E5E3E3}6.15e--2 \\ \hline
$0.7$ & 1.01e+0 & \cellcolor[HTML]{E5E3E3}2.59e--1 & 1.01e+0 & \cellcolor[HTML]{E5E3E3}2.59e--1 & 9.80e--1 & \cellcolor[HTML]{E5E3E3}2.59e--1 & 8.63e--1 & \cellcolor[HTML]{E5E3E3}2.59e--1 & 7.67e--1 & \cellcolor[HTML]{E5E3E3}2.43e--1 \\ \hline
$0.6$ & 4.48e+0 & \cellcolor[HTML]{E5E3E3}1.70e+0 & 4.44e+0 & \cellcolor[HTML]{E5E3E3}1.02e+0 & 4.15e+0 & \cellcolor[HTML]{E5E3E3}1.02e+0 & 3.33e+0 & \cellcolor[HTML]{E5E3E3}1.02e+0 & 2.91e+0 & \cellcolor[HTML]{E5E3E3}9.57e--1 \\ \hline
$0.5$ & 1.88e+1 & \cellcolor[HTML]{E5E3E3}3.92e+0 & 1.83e+1 & \cellcolor[HTML]{E5E3E3}3.92e+0 & 1.50e+1 & \cellcolor[HTML]{E5E3E3}3.92e+0 & 9.61e+0 & \cellcolor[HTML]{E5E3E3}3.92e+0 & 8.54e+0 & \cellcolor[HTML]{E5E3E3}3.68e+0 \\ \hline
\end{tabular}
}
\end{center}
\end{table}

\begin{figure}
\centering
% \begin{subfigure}[b]{0.3\textwidth}
% \centering
% \includegraphics[width=\textwidth]{/ocp/dist_v_sy}
% \caption{}
% \label{fig:v_dist_s}
% \end{subfigure}
% \hfill
\begin{subfigure}[b]{0.49\textwidth}
\centering
\includegraphics[width=0.85\textwidth]{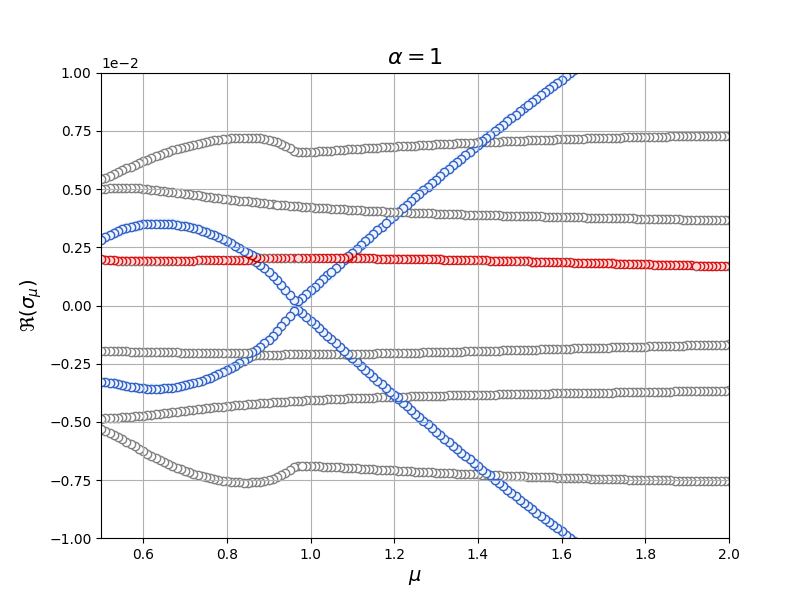}
\caption{}
\label{fig:n_eig1}
\end{subfigure}
\hfill
\begin{subfigure}[b]{0.49\textwidth}
\centering
\includegraphics[width=0.85\textwidth]{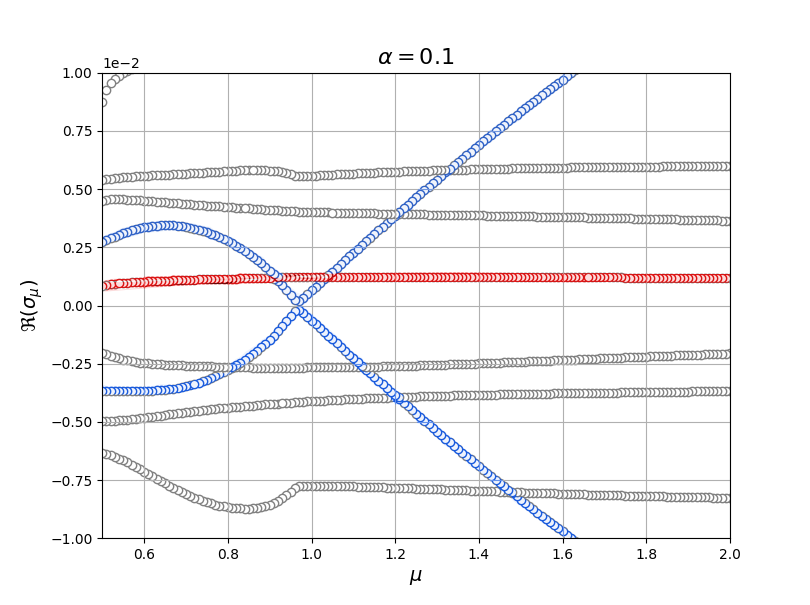}
\caption{}
\label{fig:n_eig_1e1}
\end{subfigure}\\
%\begin{subfigure}[b]{0.3\textwidth}
% \centering
% \includegraphics[width=\textwidth]{/ocp/dist_v_asy}
% \caption{}
% \label{fig:v_dist_as}
% \end{subfigure}
% \hfill
\begin{subfigure}[b]{0.49\textwidth}
\centering
\includegraphics[width=0.85\textwidth]{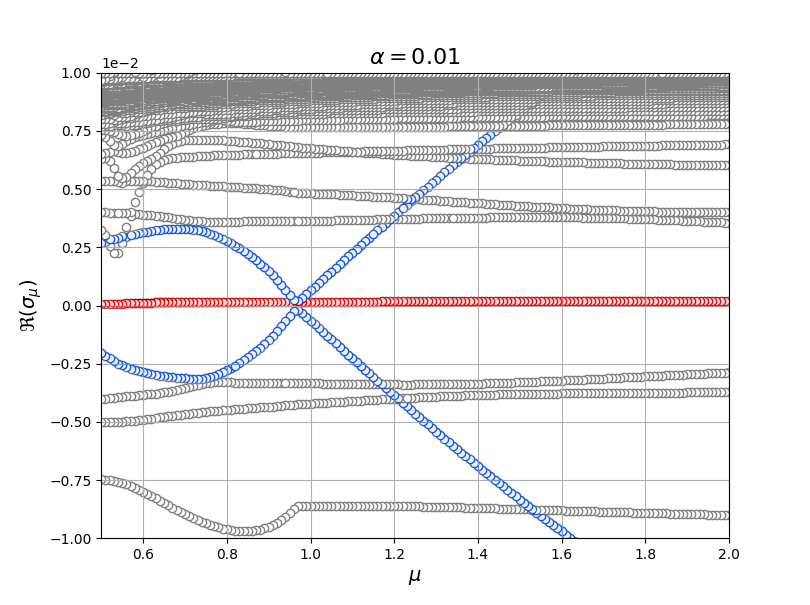}
\caption{}
\label{fig:n_eig_1e2}
\end{subfigure}
\hfill
\begin{subfigure}[b]{0.49\textwidth}
\centering
\includegraphics[width=0.85\textwidth]{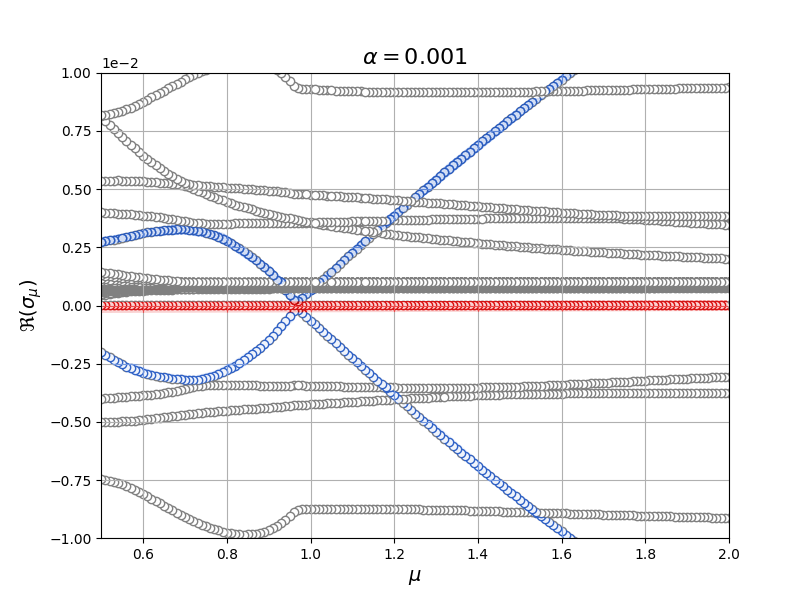}
\caption{}
\label{fig:n_eig_1e3}
\end{subfigure}\\
\caption{\emph{Neumann Control}: spectral analysis with $\alpha = 1, 0.1, 0.01, 0.001$}.
\label{fig:eig_neumann}
\end{figure}

Concerning the stability of the solution, we performed the eigenvalues analysis described in Algorithm \ref{alg:01}. We can derive several insights from the Figure \ref{fig:eig_neumann}, which represents the global eigenvalue problem for the natural branch, against the parameter $\mu$ such that $\Re(\sigma_{\mu}) = [-0.01, 0.01]$.\\
We plot the first $N_{eig} = 100$ eigenvalues of the linearized system \eqref{eq:eigen_ocp} around the global optimal solution, using a Krylov-Schur algorithm.
From the plot, we observe two eigenvalues (highlighted with blue markers) approaching $\Re(\sigma_{\mu}) = 0$: we will refer to this behavior as \emph{shears phenomenon}. Moreover, the number of positive eigenvalues grows inversely with the value of the penalization parameter, and the negative eigenvalues are lowering except for the negative shear eigenvalue.
Furthermore, the positive real eigenvalues
accumulate in the value of $\alpha$: this is very clear in subplots \ref{fig:n_eig_1e2} and \ref{fig:n_eig_1e3}. From the plot, a single eigenvalue (denoted by red markers) approaching zero is visible. \\
One of the conclusions we can obtain from the global eigenvalue analysis is how the concentration of negative eigenvalues is affected by the greater action of the control variable obtained by decreasing $\alpha$: for a fixed range of $\Re(\sigma_{\mu})$, decreasing $\alpha$ (i.e., a more controlled system) results in larger number of positive eigenvalues in $\Re(\sigma_{\mu})$.
Unfortunately, we cannot derive information about the physical stability of the global solution from the performed global eigenvalue analysis, since similar eigenvalue structures are observed for both the natural and non-natural branches (only the former being shown here for the sake of brevity). Therefore, our considerations throughout the work are limited to the numerical stability represented by natural optimality, as discussed above.
Thus, we conclude that the Neumann control is not able to fully steer uncontrolled solutions towards the desired symmetric configuration. However, this will be achieved in the next Section, where a stronger control action will be presented.

\subsection{Distributed Control: strong steering}
\label{distributed}
This Section deals with a distributed control in $\Omega_{u} \equiv \Omega$, thus the control variable $u$ acts as an external forcing term on the whole domain.
Here we consider again $\Gamma_{\text{wall}} = \Gamma_{0} \cup \Gamma_{\text{D}}$. Given $\mu \in \Cal P$, the optimal solution $X \in \mathbb X$ satisfies the following system:
\begin{equation}
\label{eq:Distributed_eq}
\begin{cases}
v\mathbb{I}_{\Gamma_{\text{obs}}} -\mu \Delta w - v\cdot\nabla w + (\nabla v)^T w + \nabla q= v_\text{d} \mathbb{I}_{\Gamma_{\text{obs}}} \quad &\text{in} \ \Omega, \\
\nabla \cdot w = 0 \quad &\text{in} \ \Omega, \\
w =0 \quad &\text{on} \ \Gamma_{\text{in}} \cup \Gamma_{\text{wall}}, \\
- qn + (\mu \nabla w) n = 0 \quad &\text{on} \ \Gamma_{\text{out}}, \\
\alpha u = w \quad &\text{in} \ \Omega, \\
-\mu \Delta v + v\cdot\nabla v + \nabla p=u \quad &\text{in} \ \Omega, \\
\nabla \cdot v = 0 \quad &\text{in} \ \Omega, \\
v = v_{\text{in}} \quad &\text{on} \ \Gamma_{\text{wall}}, \\
v = 0 \quad &\text{on} \ \Gamma_{0}, \\
- pn + (\mu \nabla v) n = 0 \quad &\text{on} \ \Gamma_{\text{out}}, \\
\end{cases}
\end{equation}
First of all, we underline that in distributed \ocp s the action of the control is usually stronger and it affects deeply
the original system.
\begin{figure}
\centering
% \begin{subfigure}[b]{0.3\textwidth}
% \centering
% \includegraphics[width=\textwidth]{/ocp/dist_v_sy}
% \caption{}
% \label{fig:v_dist_s}
% \end{subfigure}
% \hfill
\begin{subfigure}[b]{0.49\textwidth}
\centering
\includegraphics[width=\textwidth]{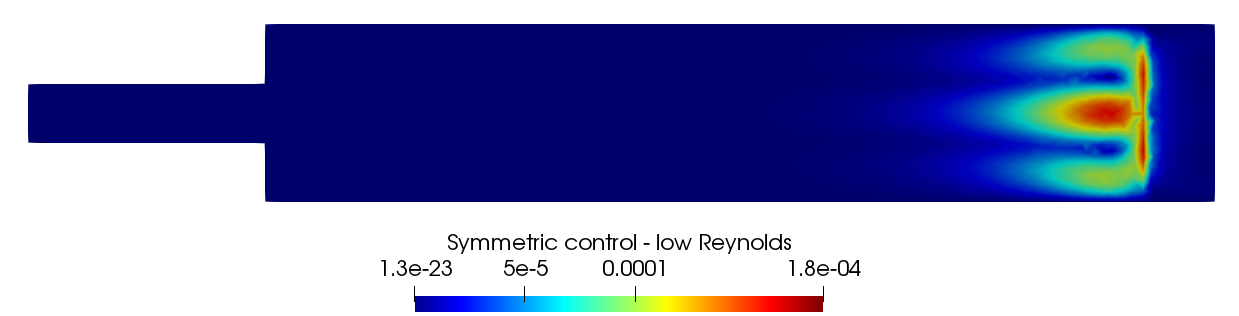}
\caption{}
\label{fig:u_dist_s_low}
\end{subfigure}
\hfill
\begin{subfigure}[b]{0.49\textwidth}
\centering
\includegraphics[width=\textwidth]{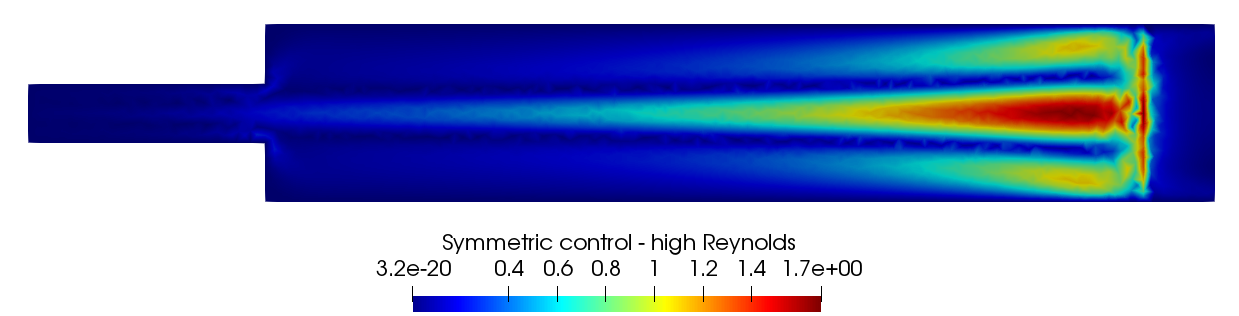}
\caption{}
\label{fig:u_dist_s_high}
\end{subfigure}\\
%\begin{subfigure}[b]{0.3\textwidth}
% \centering
% \includegraphics[width=\textwidth]{/ocp/dist_v_asy}
% \caption{}
% \label{fig:v_dist_as}
% \end{subfigure}
% \hfill
\begin{subfigure}[b]{0.49\textwidth}
\centering
\includegraphics[width=\textwidth]{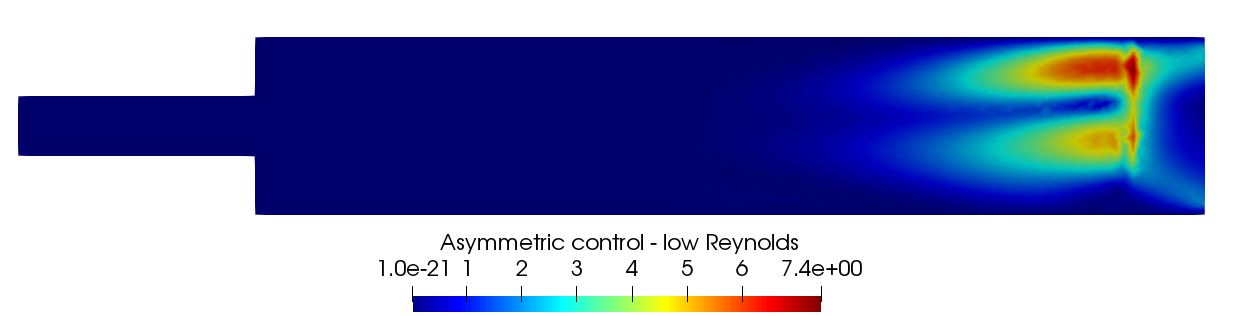}
\caption{}
\label{fig:u_dist_as_low}
\end{subfigure}
\hfill
\begin{subfigure}[b]{0.49\textwidth}
\centering
\includegraphics[width=\textwidth]{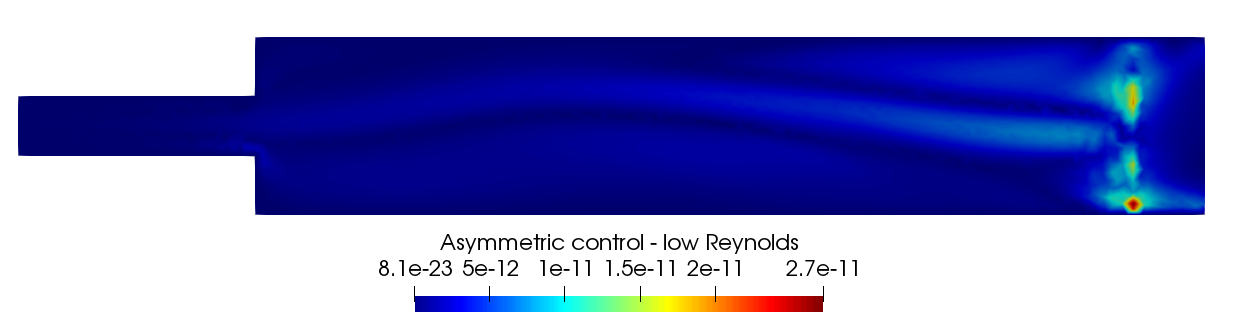}
\caption{}
\label{fig:u_dist_as_high}
\end{subfigure}\\
\caption{\emph{Distributed Control}: optimal control profiles for $\alpha = 0.01$. Left: $\mu = 2$ in (a) and (c); right: $\mu = 0.5$ in (b) and (d). Top: symmetric target in (a) and (b); bottom: asymmetric target in (c) and (d).}
\label{fig:Dist_solution}
\end{figure}

\begin{figure}
\centering
\includegraphics[width=0.489\textwidth]{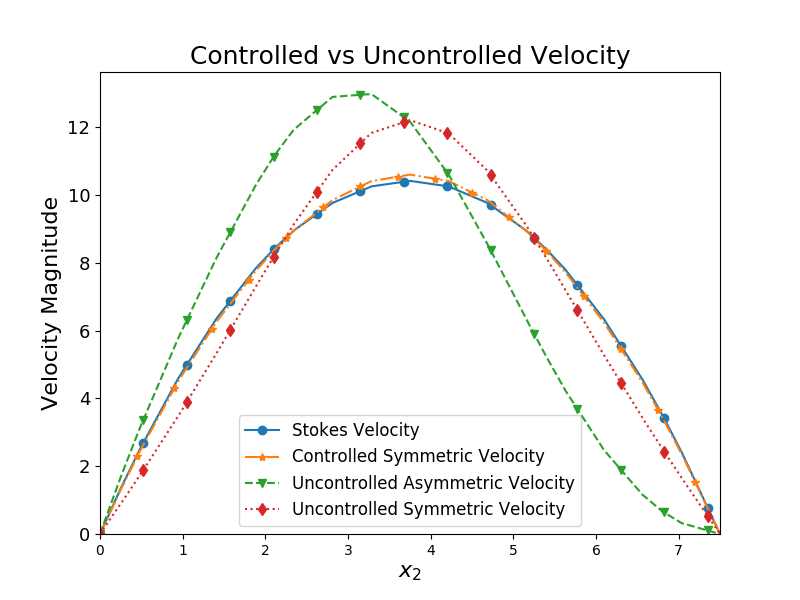}
\includegraphics[width=0.46\textwidth]{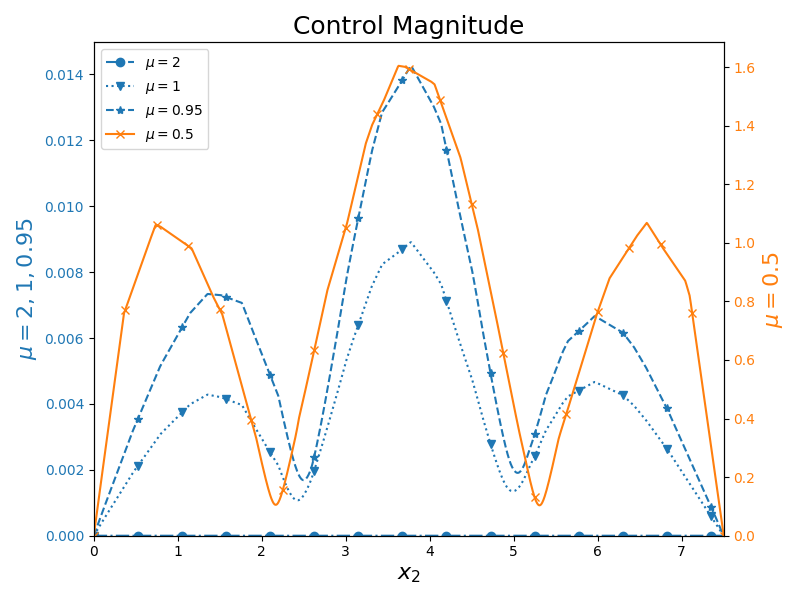}
\caption{\emph{Distributed Control}. \emph{Left}: comparison of velocity profiles in the controlled and uncontrolled cases for $\alpha = 0.01$, $\mu = 0.5$ on $\Gamma_{\text{obs}}$ w.r.t.\ the symmetric desired profile when following the natural optimal branch. \emph{Right}: representation of control variable evolution  for $\alpha = 0.01$, $\mu=2, 1, 0.95, 0.5$ for $x_1 = 45$ when following the natural optimal branch.}
\label{fig:dist_s_mag}
\end{figure}

\begin{figure}
\centering
\includegraphics[width=0.489\textwidth]{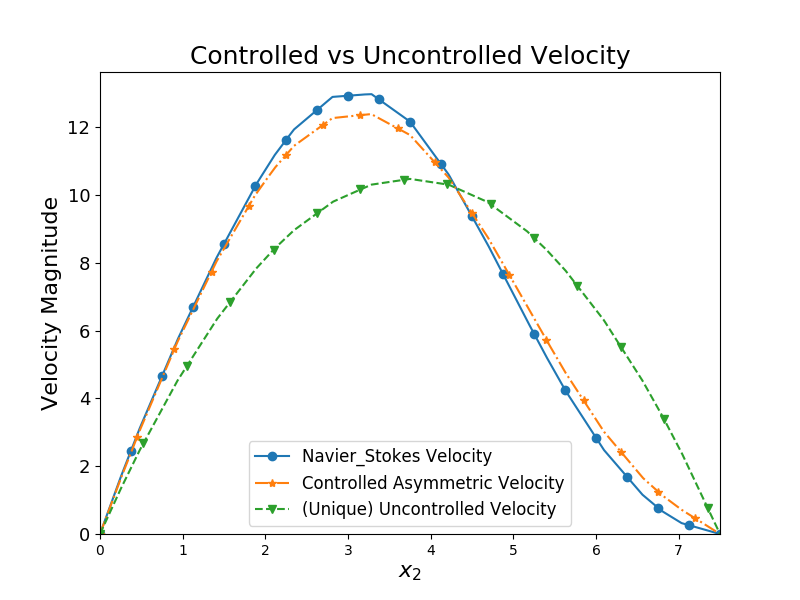}
\includegraphics[width=0.46\textwidth]{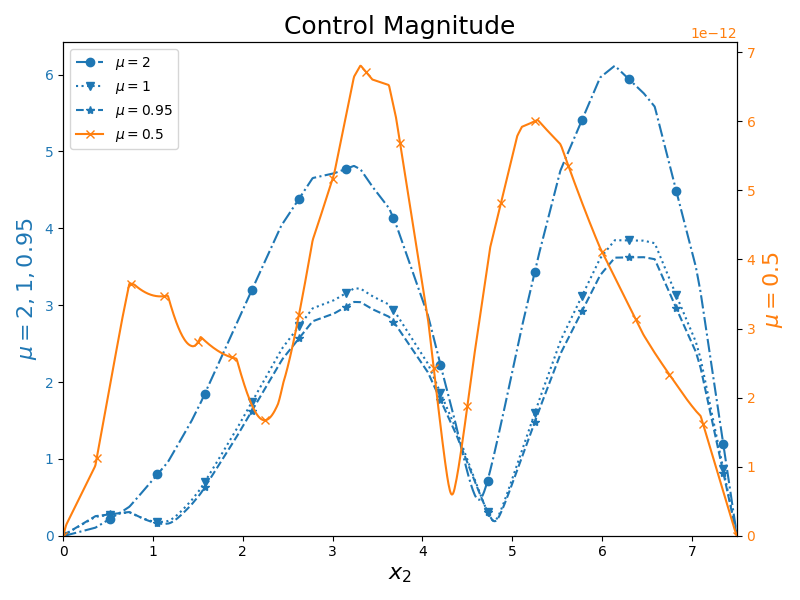}
\caption{\emph{Distributed Control}. \emph{Left}: comparison of velocity profiles in the controlled and uncontrolled cases for $\alpha = 0.01$, $\mu = 2.$ on $\Gamma_{\text{obs}}$ w.r.t.\ the asymmetric desired profile when following the natural optimal branch. \emph{Right}: representation of control variable evolution  for $\alpha = 0.01$, $\mu=2, 1, 0.95, 0.5$ for $x_1 = 45$ when following the natural optimal branch.}
\label{fig:dist_as_mag}
\end{figure}

To show this, we will steer the system towards either symmetric or asymmetric desired profiles $v_\text{d}$:
\begin{itemize}
\item[{$\small{\circ}$}] \emph{Symmetric target}: the aim of this setting is to steer the solution of \eqref{eq:Distributed_eq} to a symmetric profile. We plot two representative control solutions in Figures \ref{fig:u_dist_s_low} and \ref{fig:u_dist_s_high}, obtained for $\mu = 2$ and $\mu = 0.5$ when following the natural optimal branch, which is composed of symmetric solutions.
The stronger action of the control allows the controlled velocity profile to be more diffusive compared to the uncontrolled symmetric profile, as represented in the left plot of Figure \ref{fig:dist_s_mag}, corresponding to the observed slice of the velocity solution for $\mu = 0.5$: in this case the controlled velocity (orange line) and the symmetric target (blue line) almost coincide.
The right plot of Figure \ref{fig:dist_s_mag} shows that a slightly asymmetric control is required only near the critical value $\mu^*$ (also compare to Figures \ref{fig:u_dist_s_low} and \ref{fig:u_dist_s_high} for the cases $\mu = 2$ and $\mu = 0.5$).
Furthermore, the control action is clearly higher when the Re value increases. Indeed, for $\mu = 2$ the control exclusively acts in the proximity of $\Gamma_{\text{obs}}$ with a maximum magnitude of $1.8 \small{\cdot} 10^{-4}$, while for $\mu = 0.5$ it reaches a value of $1.6$ of magnitude.\\
A further non-natural optimal branch exists, and is made of symmetric solutions, but is hardly reachable by numerical continuation methods unless tailored guesses are provided restarting Algorithm \ref{alg:01} in a small neighborhood of $\mu^*$.

\item[{$\small{\circ}$}] \emph{Asymmetric target}: in this case, we desire to recover the asymmetric target for all $\mu \in \mathcal P$.
We plot two representative control solutions in Figures \ref{fig:u_dist_as_low} and \ref{fig:u_dist_as_high}, obtained for $\mu = 2$ and $\mu = 0.5$ when following the natural optimal branch, which is made of asymmetric solutions.
The action of the control is also visible in the left plot of Figure \ref{fig:dist_as_mag}, obtained for $\mu = 2$: indeed, we see how the flux over $\Gamma_{\text{obs}}$ is pushed towards the domain wall (orange line), in contrast to the symmetric profile of the uncontrolled velocity (green line). Namely, also in this case, the distributed control is able to drive the solution towards the desired state.
In order to do so, the control variable has to be large when $\mu > \mu^*$, i.e.\ when the uncontrolled configuration on $\Gamma_{\text{obs}}$ would lead to a symmetric profile. Indeed, in Figure \ref{fig:u_dist_as_low} the maximum control value reaches $7$ for $\mu = 2$ in the upper part of the domain. In contrast Figure \ref{fig:u_dist_as_high} it lowers to $10^{-11}$ for $\mu = 0.5$ when the stable asymmetric velocity solution does not need to be controlled by an external forcing term. This is confirmed in the right plot of Figure \ref{fig:dist_as_mag} for several values of $\mu$.\\
Also in this case a non-natural optimal branch (featuring symmetric solutions) continues to exist, but it is numerically difficult to reach.
\end{itemize}

\begin{table}
\caption{\emph{Distributed Control}: comparison of the functional value. (Sym.) Natural optimal branch for symmetric target. (Asym.) Natural optimal branch for asymmetric target. (Sym.-U.) Unstable uncontrolled solution with symmetric target. (Asym.-S.) Stable uncontrolled solution with asymmetric target. (B.M.E.) Below machine epsilon.}
\label{Distributed_J}
\begin{center}
\tabcolsep=0.11cm
\footnotesize{
\begin{tabular}{|c||c|c||c|c||c|c||c|c||c|c||}
\hline
& Sym.-U. & \cellcolor[HTML]{E5E3E3}Asym.-S. & Sym.& \cellcolor[HTML]{E5E3E3}Asym. & Sym. & \cellcolor[HTML]{E5E3E3}Asym. & Sym. & \cellcolor[HTML]{E5E3E3}Asym. & Sym. & \cellcolor[HTML]{E5E3E3}Asym. \\ \cline{2-11}
\multirow{-2}{*}{$\mu$} & \multicolumn{2}{c||}{Uncontrolled} & \multicolumn{2}{c||}{$\alpha = 1$} & \multicolumn{2}{c||}{$\alpha = 0.1$} & \multicolumn{2}{c||}{$\alpha = 0.01$} & \multicolumn{2}{c||}{$\alpha = 0.001$} \\ \hline
$2$ & 5.14e--9 & \cellcolor[HTML]{E5E3E3}1.88e+1 &5.06e--9 & \cellcolor[HTML]{E5E3E3}1.81e+1& 4.51e--9 & \cellcolor[HTML]{E5E3E3}1.36e+1 & 2.22e--9 & \cellcolor[HTML]{E5E3E3}4.23e+0 & 4.04e--10 & \cellcolor[HTML]{E5E3E3}5.66e--1 \\ \hline
$1.5$ & 4.38e--6 & \cellcolor[HTML]{E5E3E3}1.88e+1 & 4.29e--6 & \cellcolor[HTML]{E5E3E3}1.77e+1 & 3.61e--6 & \cellcolor[HTML]{E5E3E3}1.20e+1 & 1.46e--6 & \cellcolor[HTML]{E5E3E3}3.09e+0 & 2.28e--7& \cellcolor[HTML]{E5E3E3}3.87e--1 \\ \hline
$1$ & 4.10e--3 & \cellcolor[HTML]{E5E3E3}1.86e+1 & 3.95e--3 & \cellcolor[HTML]{E5E3E3}1.67e+1 & 2.99e--3 & \cellcolor[HTML]{E5E3E3}9.15e+0 & 9.14e--4& \cellcolor[HTML]{E5E3E3}1.86e+0& 1.23e--4 & \cellcolor[HTML]{E5E3E3}2.17e--1 \\ \hline
$0.9$ & 1.63e--2 & \cellcolor[HTML]{E5E3E3}1.84e+1 & 1.56e--2 & \cellcolor[HTML]{E5E3E3}1.54e+1 & 1.14e--2 & \cellcolor[HTML]{E5E3E3}7.88e+0 & 3.26e--3 & \cellcolor[HTML]{E5E3E3}1.50e+0 &4.26e--4 & \cellcolor[HTML]{E5E3E3}1.73e--1 \\ \hline
$0.8$ & 6.52e--2 & \cellcolor[HTML]{E5E3E3} 1.54e+1& 6.21e--2 & \cellcolor[HTML]{E5E3E3}1.31e+1 & 4.36e--2 & \cellcolor[HTML]{E5E3E3}6.06e+0 & 1.14e--2 & \cellcolor[HTML]{E5E3E3}1.08e+0& 1.45e--3 & \cellcolor[HTML]{E5E3E3}1.22e--1 \\ \hline
$0.7$ & 2.59e--1 & \cellcolor[HTML]{E5E3E3}1.15e+1 & 2.45e--1 & \cellcolor[HTML]{E5E3E3}9.28e+0 & 1.63e--1 & \cellcolor[HTML]{E5E3E3}3.68e+0 &3.93e--2& \cellcolor[HTML]{E5E3E3}6.16e--1 & 4.81e--3 & \cellcolor[HTML]{E5E3E3}6.90e--2\\ \hline
$0.6$ & 1.70e+0 & \cellcolor[HTML]{E5E3E3}5.34e+0 & 9.54e--1 & \cellcolor[HTML]{E5E3E3}3.76e+0 & 5.94e--1 & \cellcolor[HTML]{E5E3E3}1.24e+0 & 1.28e--1 & \cellcolor[HTML]{E5E3E3}2.00e--1 & 1.70e--2 & \cellcolor[HTML]{E5E3E3}2.22e--2 \\ \hline
$0.5$ & 3.92e+0 & \cellcolor[HTML]{E5E3E3}B.M.E. &3.59e+0 & \cellcolor[HTML]{E5E3E3}B.M.E. & 2.04e+0 & \cellcolor[HTML]{E5E3E3}B.M.E. & 3.92e--1 & \cellcolor[HTML]{E5E3E3}B.M.E. & 4.47e--2 & \cellcolor[HTML]{E5E3E3}B.M.E. \\ \hline
\end{tabular}
}
\end{center}
\end{table}

\begin{figure}
\centering
\includegraphics[scale=0.65]{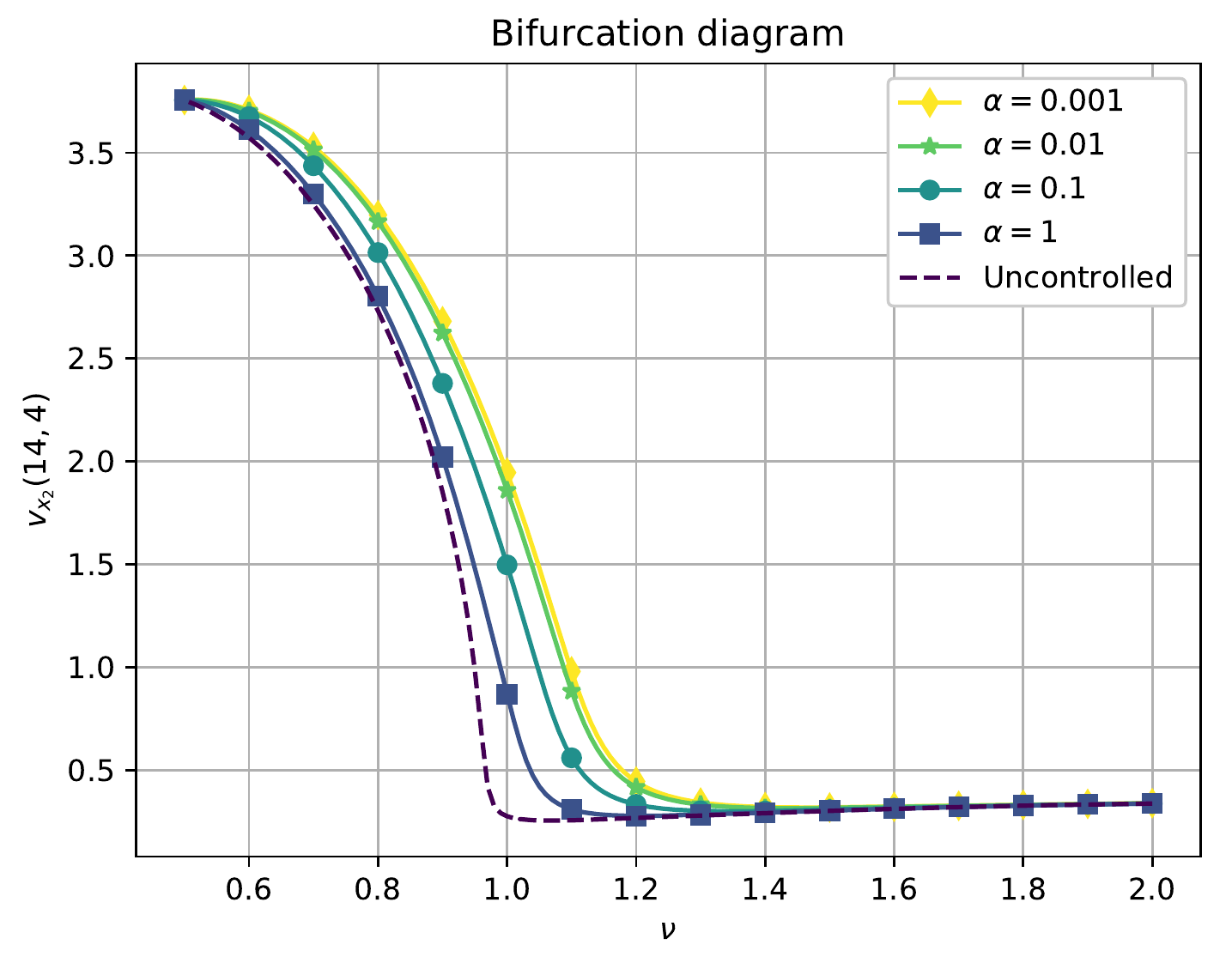}
\caption{\emph{Distributed Control}: bifurcation diagram (upper branch only) for controlled state velocity obtained with $\alpha = 1, 0.1, 0.01, 0.001$ and asymmetric target, compared to the uncontrolled velocity.}
\label{fig:early_bif}
\end{figure}

We show the comparison of the values of the cost functional \eqref{eq:J_NS} in Table \ref{Distributed_J} for the reached natural branch for both symmetric and asymmetric targets. Several values of $\mu$ (rows) and $\alpha$ (columns) have been analyzed and compared to the uncontrolled functional, computed as in the Neumann test case. As expected, we notice that the functional is lower for smaller $\alpha$.
For the symmetric target, the action of the distributed control is indeed able to steer the solution towards the desired symmetric profile. Indeed, for $\mu = 0.05$, the choice $\alpha = 0.01$ shows that the functional is decreased by a $90\%$  w.r.t. its uncontrolled counterpart, while for $\alpha = 0.001$, the cost functional is almost decreased by $99\%$.
Similarly, for the asymmetric target, the maximum action of the control variable is given for low Reynolds and, for $\mu = 2.$ we can observe a decrease of the functional of the $77.5\% $ for $\alpha = 0.01$, up to a $97\%$ for $\alpha = 0.001$. We remark that no control action is needed for $\mu = 0.5 \approx 0.49$, which is the parameter value for which the asymmetric $v_\text{d}$ was computed: this was underlined by very low values of \eqref{eq:J_NS}, which were below the machine precision.

\begin{figure}
\centering
% \begin{subfigure}[b]{0.3\textwidth}
% \centering
% \includegraphics[width=\textwidth]{/ocp/dist_v_sy}
% \caption{}
% \label{fig:v_dist_s}
% \end{subfigure}
% \hfill
\begin{subfigure}[b]{0.49\textwidth}
\centering
\includegraphics[width=0.85\textwidth]{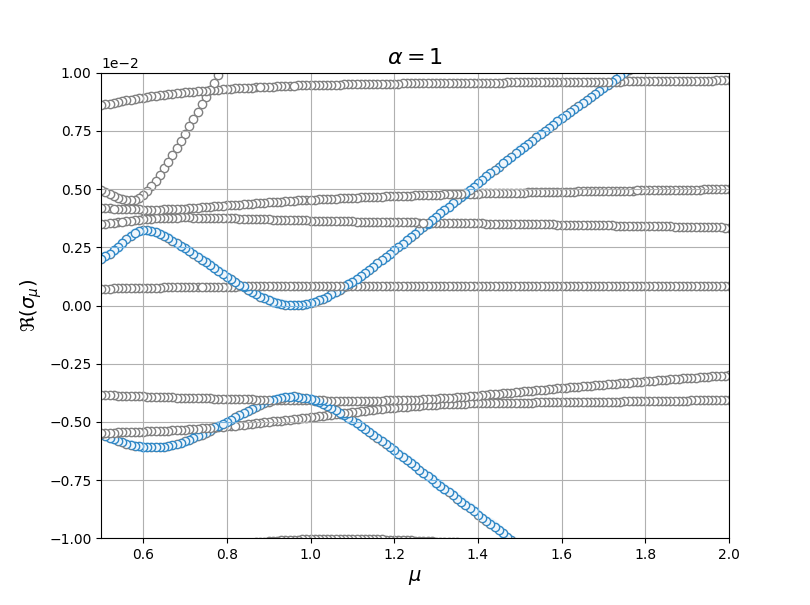}
\caption{}
\label{fig:dist_eig_1}
\end{subfigure}
\hfill
\begin{subfigure}[b]{0.49\textwidth}
\centering
\includegraphics[width=0.85\textwidth]{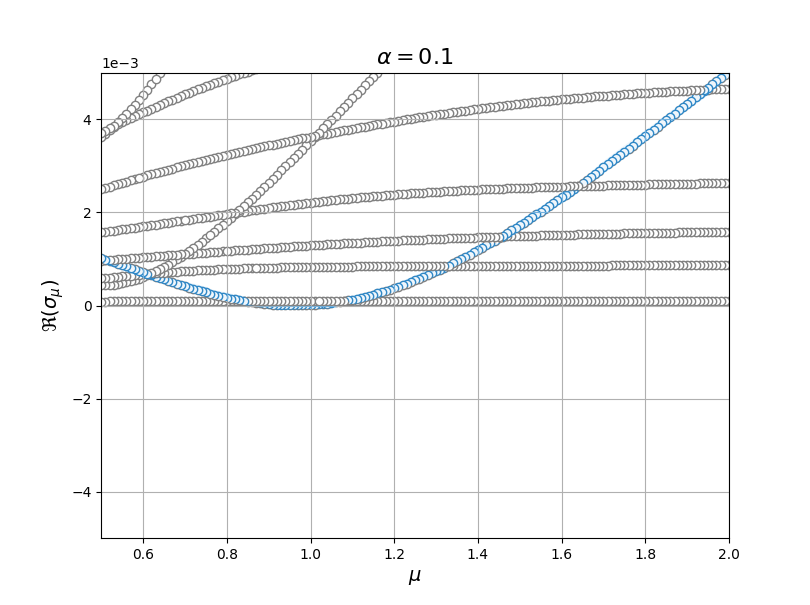}
\caption{}
\label{fig:dist_eig_1e1}
\end{subfigure}\\
%\begin{subfigure}[b]{0.3\textwidth}
% \centering
% \includegraphics[width=\textwidth]{/ocp/dist_v_asy}
% \caption{}
% \label{fig:v_dist_as}
% \end{subfigure}
% \hfill
\begin{subfigure}[b]{0.49\textwidth}
\centering
\includegraphics[width=0.85\textwidth]{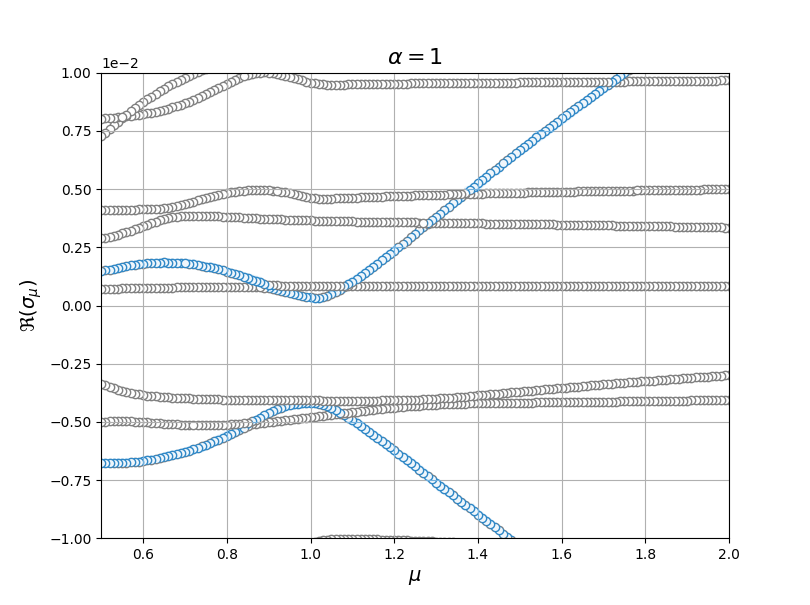}
\caption{}
\label{fig:dist_a_eig_1}
\end{subfigure}
\hfill
\begin{subfigure}[b]{0.49\textwidth}
\centering
\includegraphics[width=0.85\textwidth]{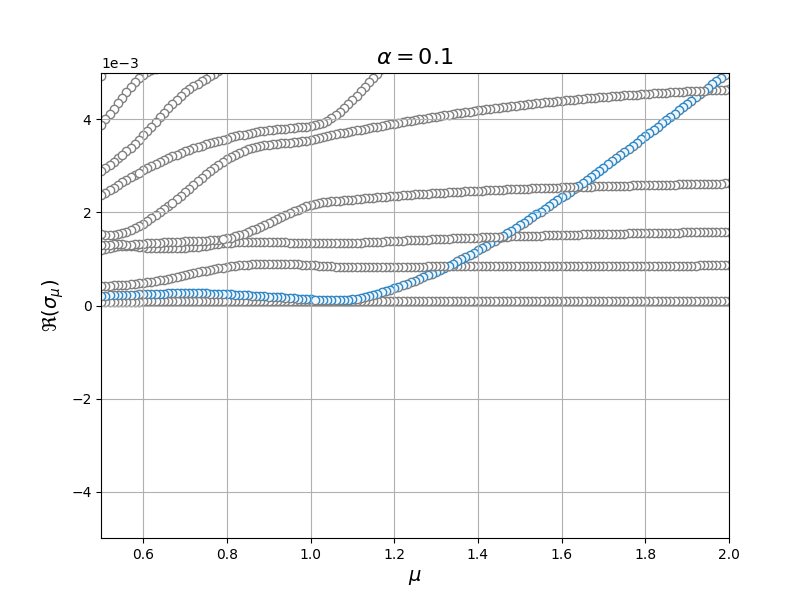}
\caption{}
\label{fig:dist_a_eig_1e1}
\end{subfigure}\\
\caption{\emph{Distributed Control}: spectral analysis with $\alpha = 1, 0.1$ (left to right) for the natural optimal branch with symmetric (top) and asymmetric (bottom) targets.}
\label{fig:dist_eig}
\end{figure}
The plots of Figure \ref{fig:dist_eig} represent the spectral analysis for this optimal control problem: in particular, Figures \ref{fig:dist_eig_1} ($\alpha = 1$) and \ref{fig:dist_eig_1e1} ($\alpha = 0.1$) are associated with the symmetric target when following the corresponding natural optimal branch, while Figures \ref{fig:dist_a_eig_1} ($\alpha = 1$) and \ref{fig:dist_a_eig_1e1} ($\alpha = 0.1$) consider the asymmetric target when following its natural optimal branch. As the behavior between the top and bottom panels of Figure \ref{fig:dist_eig} is comparable, we will only comment in the following on the role of $\alpha$. %Furthermore, as we discussed in the previous test case, we cannot recover stability information about the configuration obtained by the optimal solution by means of the global eigenvalue problem, because the remaining (i.e., non-natural optimal) branches show very similar patterns to the ones in Figure \ref{fig:dist_eig}.

We plot the eigenvalues for $\alpha = 1$ in $\Re(\sigma_{\mu}) = [-0.01, 0.01]$ and for $\alpha = 0.1$ in $\Re(\sigma_{\mu}) = [-0.005, 0.005]$. For this test case, the predominance of positive eigenvalues is visible also for large values of the penalization parameter. The smaller is $\alpha$, the more negative eigenvalues are lowered. For all the $\alpha$ taken into account, the shears phenomenon does not appear: for $\alpha = 1$ a small trace of the shears structure is still visible (highlighted in blue) in Figures \ref{fig:dist_eig_1} and \ref{fig:dist_a_eig_1} where the bottom part of the shear is pushed away from $\Re(\sigma_{\mu}) = 0$. Instead, for $\alpha = 0.1$  the shears structure is completely lost: Figures \ref{fig:dist_a_eig_1} and \ref{fig:dist_a_eig_1e1} show that only one eigenvalue
(representing the top of the shears, and marked in blue) approaches $\Re(\sigma_{\mu}) = 0$ without crossing it.

We finally notice that the point $\mu^{**}$, where the upper shears curve is the closest to the axis $\Re(\sigma_{\mu}) = 0$, allows to obtain further information on the bifurcating phenomenon. From Figure \ref{fig:dist_eig}, $\mu^{**} \approx 0.96$ for the symmetric target, regardless of $\alpha$, while requiring an asymmetric target leads to $\mu^{**} \in [1.0, 1.2]$ with a mild dependence on $\alpha$. 

Thanks to the aid of Figure \ref{fig:early_bif}, which shows the bifurcation diagram for the controlled solution with asymmetric target (a similar plot can be obtained for the symmetric target as well, but is here omitted because the lines almost overlap), we can state that optimal control is not only able to steer the state solution towards a desired branch, but may also affect the location of the bifurcation point.\\ The role of the penalization parameter $\alpha$ will be clarified in the next Section and it will result into a completely new optimal solution behavior in Section \ref{dirichlet}.

\subsection{Channel Control: the $\boldsymbol \alpha$ effect.}
\label{inlet}
This Section aims at describing how the value of the penalization parameter $\alpha$ can affect the natural convergence towards a symmetric target over $\Gamma_{\text{obs}}$. 

\begin{figure}
\centering
% \begin{subfigure}[b]{0.3\textwidth}
% \centering
% \includegraphics[width=\textwidth]{/ocp/dist_v_sy}
% \caption{}
% \label{fig:v_dist_s}
% \end{subfigure}
% \hfill
\begin{subfigure}[b]{0.49\textwidth}
\centering
\includegraphics[width=\textwidth]{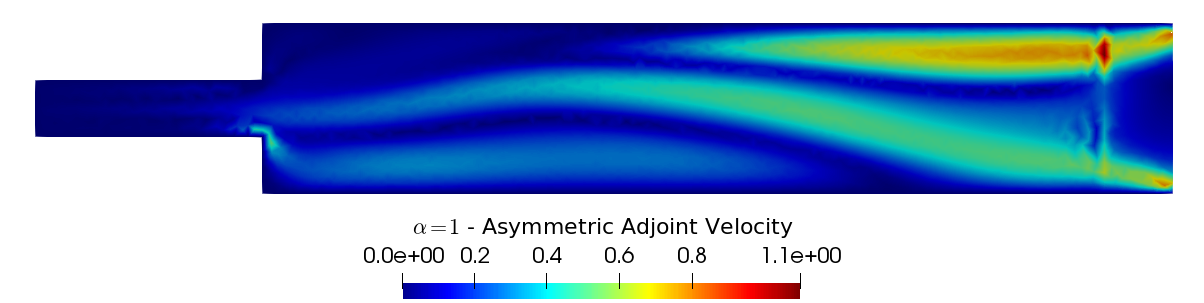}
\caption{}
\label{fig:w_channel_as}
\end{subfigure}
\hfill
\begin{subfigure}[b]{0.49\textwidth}
\centering
\includegraphics[width=\textwidth]{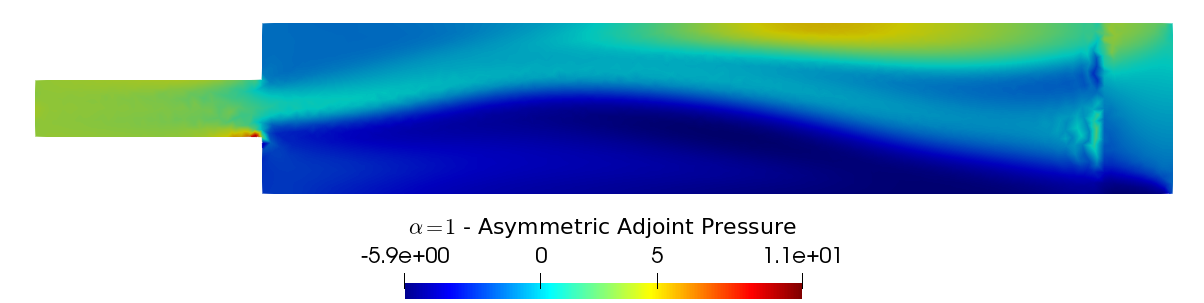}
\caption{}
\label{fig:q_channel_as}
\end{subfigure}\\
%\begin{subfigure}[b]{0.3\textwidth}
% \centering
% \includegraphics[width=\textwidth]{/ocp/dist_v_asy}
% \caption{}
% \label{fig:v_dist_as}
% \end{subfigure}
% \hfill
\begin{subfigure}[b]{0.49\textwidth}
\centering
\includegraphics[width=\textwidth]{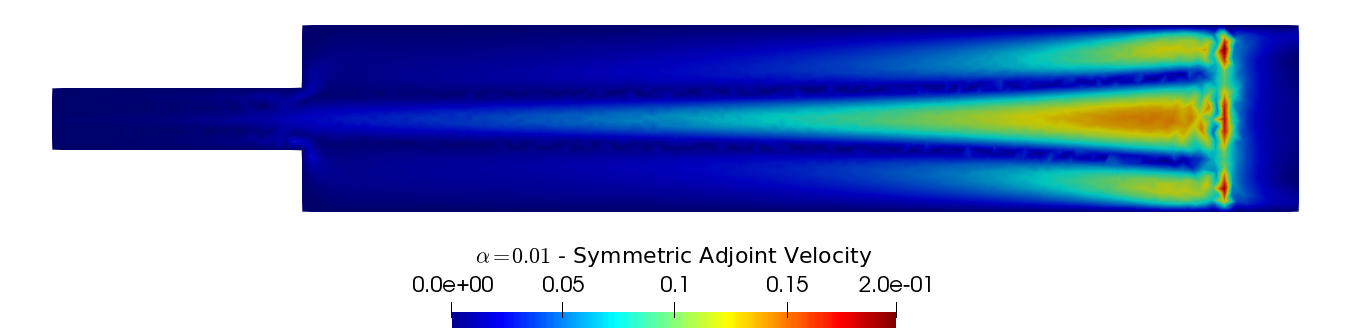}
\caption{}
\label{fig:w_channel_s}
\end{subfigure}
\hfill
\begin{subfigure}[b]{0.49\textwidth}
\centering
\includegraphics[width=\textwidth]{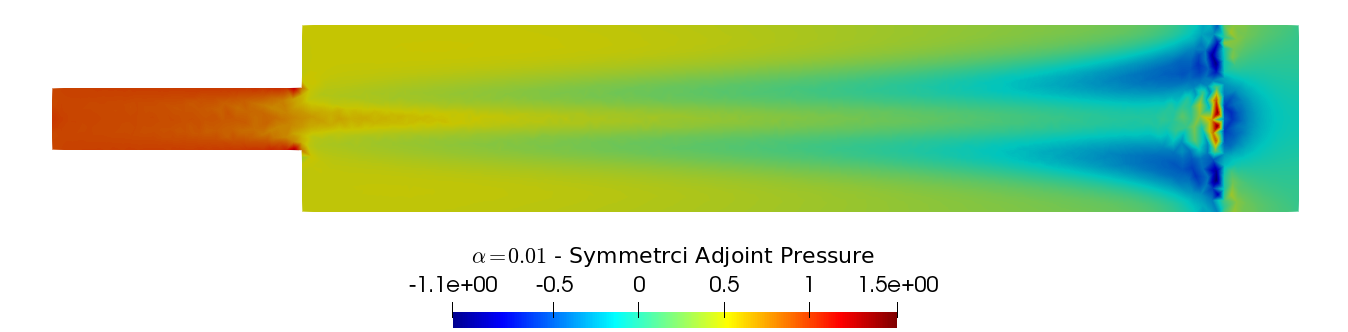}
\caption{}
\label{fig:q_channel_s}
\end{subfigure}\\
\caption{\emph{Channel Control}: two optimal solutions for adjoint velocity and pressure for $\mu = 0.5$: $\alpha = 1$ in (a) and (b), and $\alpha = 0.01$ in (c) and (d), respectively.}
\label{fig:Channel_solution}
\end{figure}

\begin{figure}
\centering
\includegraphics[width=0.489\textwidth]{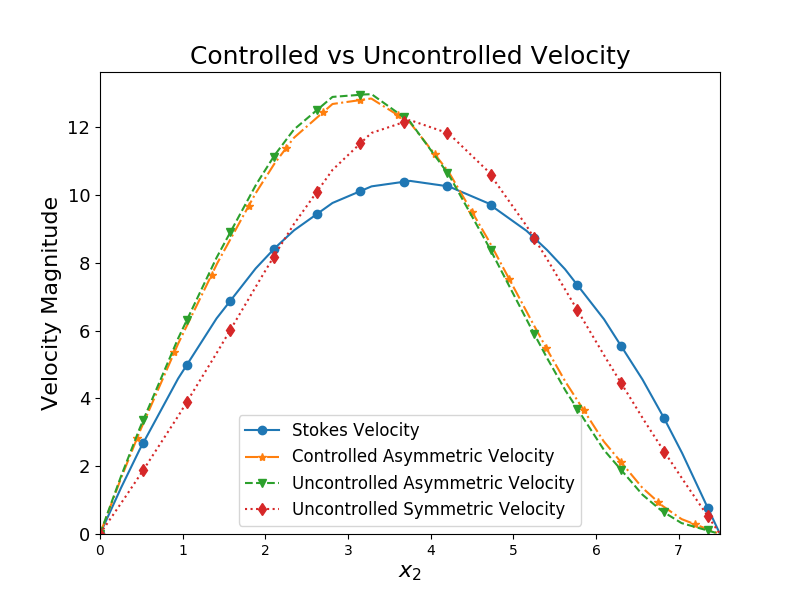}
\includegraphics[width=0.46\textwidth]{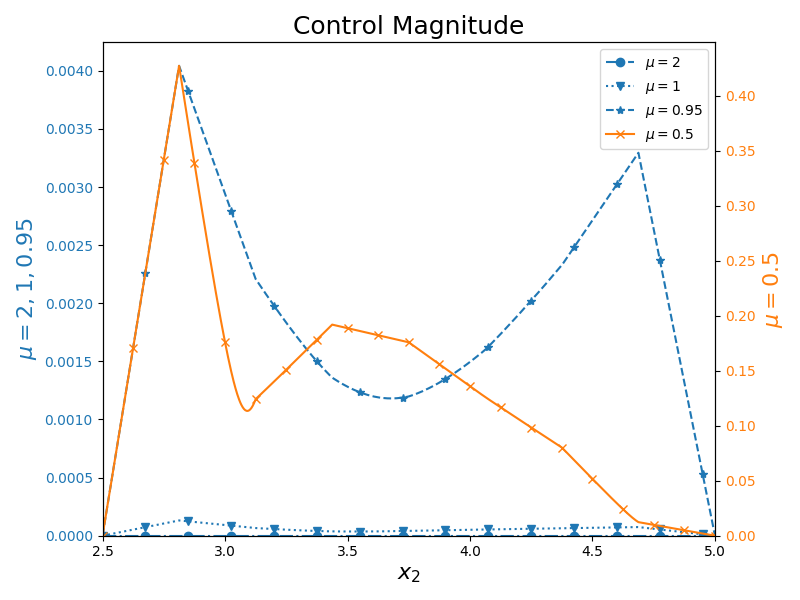}
\caption{\emph{Channel Control}. \emph{Left}: comparison of velocity profiles in the controlled and uncontrolled cases for $\alpha = 1$, $\mu = 0.5$ on $\Gamma_{\text{obs}}$ w.r.t.\ the symmetric desired profile when following the natural optimal branch. \emph{Right}: representation of control variable evolution for $\alpha = 1$, $\mu = 2, 1, 0.95, 0.5$ at $x_1 = 10$ when following the natural optimal branch.}
\label{fig:ch_as_mag}
\end{figure}

\begin{figure}
\centering
\includegraphics[width=0.489\textwidth]{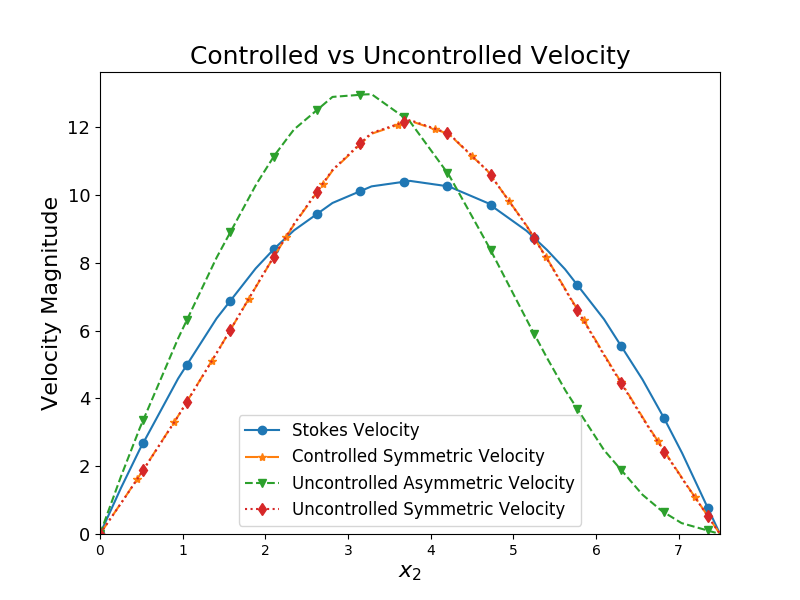}
\includegraphics[width=0.46\textwidth]{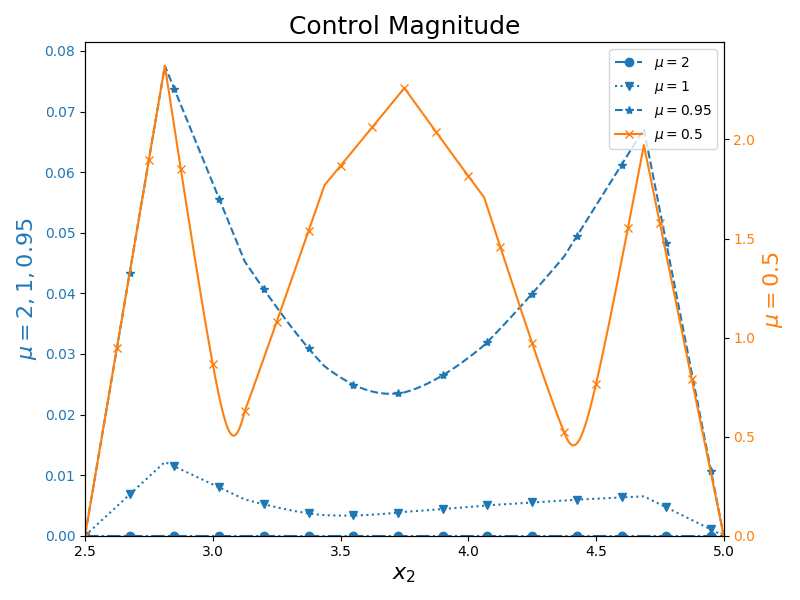}
\caption{\emph{Channel Control}. \emph{Left}: comparison of velocity profiles in the controlled and uncontrolled cases for $\alpha = 0.01$, $\mu = 0.5$ on $\Gamma_{\text{obs}}$ w.r.t.\ the symmetric desired profile when following the natural optimal branch. \emph{Right}: representation of control variable evolution for $\alpha = 0.01$, $\mu=2, 1, 0.95, 0.5$ for $x_1 = 10$ when following the natural optimal branch.}
\label{fig:ch_s_mag}
\end{figure}

Towards this goal, we analyzed the action of a control variable defined at the end of inlet channel, i.e. $\Omega_u = \Gamma_{\text{ch}}$, as depicted in Figure \ref{fig:channel}. The boundary $\Gamma_{\text{wall}}$ is, once again, $\Gamma_{0} \cup \Gamma_{\text{D}}$. 
Within this setting, the problem reads: given $\mu \in \Cal P$ find the optimal solution $X \in \mathbb X$ such that the following holds

\begin{equation}
\label{eq:Channel_eq}
\begin{cases}
v\mathbb{I}_{\Gamma_{\text{obs}}} -\mu \Delta w - v\cdot\nabla w + (\nabla v)^T w + \nabla q= v_\text{d} \mathbb{I}_{\Gamma_{\text{obs}}} \quad &\text{in} \ \Omega, \\
\nabla \cdot w = 0 \quad &\text{in} \ \Omega, \\
w =0 \quad &\text{on} \ \Gamma_{\text{in}} \cup \Gamma_{\text{wall}}, \\
- qn + (\mu \nabla w) n = 0 \quad &\text{on} \ \Gamma_{\text{out}}, \\
\alpha u \mathbb{I}_{\Gamma_{\text{ch}}}= w \mathbb{I}_{\Gamma_{\text{ch}}} \quad &\text{in} \ \Omega, \\
-\mu \Delta v + v\cdot\nabla v + \nabla p=u \mathbb{I}_{\Gamma_{\text{ch}}} \quad &\text{in} \ \Omega, \\
\nabla \cdot v = 0 \quad &\text{in} \ \Omega, \\
v = v_{\text{in}} \quad &\text{on} \ \Gamma_{\text{in}}, \\
v = 0 \quad &\text{on} \ \Gamma_{\text{wall}}, \\
- pn + (\mu \nabla v) n = 0 \quad &\text{on} \ \Gamma_{\text{out}}.\\
\end{cases}
\end{equation}
The optimal control acts as a forcing term capable to change the way the flow enters in the expansion channel. In Figure \ref{fig:Channel_solution} we show the adjoint velocity and pressure profiles obtained for $\mu = 0.5$ and two different penalization values, namely $\alpha = 1$ and $\alpha = 0.01$.
In the first case, following Algorithm \ref{alg:01}, the natural optimal branch presents a wall-hugging behavior, while for smaller values of $\alpha$ the control variable is able to drive the velocity towards a straight flux (see the left panels of Figures \ref{fig:ch_as_mag} and \ref{fig:ch_s_mag}). Therefore, for large values of $\alpha$ the natural optimal branch is composed by asymmetric solutions (i.e., far away from the target), while for smaller values of $\alpha$ the natural optimal branch is made of symmetric solutions.

From the right plots of Figures \ref{fig:ch_as_mag} and \ref{fig:ch_s_mag}, the control is very sensitive close to $\mu^*$ and this is shown by its asymmetric configuration both for the wall-hugging solution and the straight one. For $\alpha = 1, 0.1, 0.01$, we were able to detect two solutions using different initial guesses in the continuation method, showing symmetric and asymmetric features coexisting for some values of $\mu < \mu^*$. The smaller was $\alpha$, the more difficult was to recover the non-natural branch. For example, when $\alpha = 0.001$, the action of the control variable drives the wall-hugging phenomenon towards a straight flux so strongly that we were not able to really reconstruct the whole optimal non-natural branch. Indeed, either the Newton's solver did not converge (this happens also for $\alpha = 0.1$ and $\mu = 0.5$, compare Table \ref{Inlet_J}) or converged to the natural branch consisting in symmetric features.

As usual, the role of $\alpha$ is highlighted in reducing objective functional,
as Table \ref{Inlet_J} shows. As already specified in Sections \ref{neumann} and \ref{distributed}, the straight configuration is lowering much more the functional than the asymmetric solution, due to its similarity with the symmetric $v_\text{d}$, which is the fixed target for this test case. In this case, the role of $\alpha$ is crucial in order to reach a solution which represents better the desired state. Indeed, the control was able to steer the solution to the symmetric profile for $\alpha = 0.1, 0.01, 0.001$. From the functional point of view, we do not have a notable decrease, as can be observed in Table \ref{Inlet_J}, where the value of \eqref{eq:J_NS} is presented for different values of $\mu$ and $\alpha$ w.r.t. the uncontrolled problem solution. Yet, acting at the end of the inlet channel still allows to drive the optimal solution towards a natural convergence to the symmetric $v_\text{d}$, but the parabolic profile on $\Gamma_{\text{obs}}$ is not reached (the functional decreases only of a $10\%$ for $\mu = 0.5$ and $\alpha = 0.001$ w.r.t.\ the uncontrolled symmetric solution).

\begin{table}
\caption{\emph{Channel Control}: comparison of the functional value w.r.t. stable and unstable uncontrolled solutions. \emph{Headers:} (Nat.) Natural optimal branch. (n-Nat.) Non-natural optimal branch.
%(Sym.) the solution has symmetric profile. (Asym) the soluton has asymmetric profile.
\emph{Trailing cell characters}: (s) the solution has symmetric profile. (a) The solution has asymmetric profile. (nat-C.) Converging to natural branch despite tailored guess. (non-C.) Non-converging Newton's solver for tailored guess.}
\label{Inlet_J}
\begin{center}
\tabcolsep=0.09cm
\footnotesize{
\begin{tabular}{|c||c|c||c|c||c|c||c|c||c|c||}
\hline
%& Asym. & \cellcolor[HTML]{E5E3E3}Sym. & Asym. & \cellcolor[HTML]{E5E3E3}Sym. & Sym. & \cellcolor[HTML]{E5E3E3}Asym. & Sym. & \cellcolor[HTML]{E5E3E3}Asym. & Sym. & \cellcolor[HTML]{E5E3E3}Asym. \\ %\hline{2-11} %\cline{2-11}
& Stable & \cellcolor[HTML]{E5E3E3}Unstable & Nat. & \cellcolor[HTML]{E5E3E3}n-Nat. & Nat. & \cellcolor[HTML]{E5E3E3}n-Nat. & Nat. & \cellcolor[HTML]{E5E3E3}n-Nat. & Nat. & \cellcolor[HTML]{E5E3E3}n-Nat. \\ \cline{2-11}
\multirow{-2}{*}{$\mu$} & \multicolumn{2}{c||}{Uncontrolled} & \multicolumn{2}{c||}{$\alpha = 1$} & \multicolumn{2}{c||}{$\alpha = 0.1$} & \multicolumn{2}{c||}{$\alpha = 0.01$} & \multicolumn{2}{c||}{$\alpha = 0.001$} \\ \hline
$2$ & 5.14e--9 & \cellcolor[HTML]{E5E3E3}5.14e--9 & 5.14e--9s& \cellcolor[HTML]{E5E3E3}5.14e--9s& 5.14e--9s & \cellcolor[HTML]{E5E3E3}5.14e--9s & 5.14e--9s & \cellcolor[HTML]{E5E3E3}5.14e--9s & 5.07e--9s& \cellcolor[HTML]{E5E3E3}5.14e--9s \\ \hline
$1.5$ & 4.38e--6 & \cellcolor[HTML]{E5E3E3}4.38e--6 & 4.38e--6s& \cellcolor[HTML]{E5E3E3}4.38e--6s & 4.38e--6s & \cellcolor[HTML]{E5E3E3}4.38e--6s & 4.38e--6s & \cellcolor[HTML]{E5E3E3}4.38e--6s & 4.28e--6s & \cellcolor[HTML]{E5E3E3}4.38e--6s \\ \hline
$1$ & 4.10e--3 & \cellcolor[HTML]{E5E3E3}4.10e--3 & 4.10e--3s& \cellcolor[HTML]{E5E3E3}4.10e--3s & 4.10e--3s & \cellcolor[HTML]{E5E3E3}4.10e--3s & 4.08e--3s & \cellcolor[HTML]{E5E3E3}4.10e--3s & 3.92e--3s & \cellcolor[HTML]{E5E3E3}4.10e--3s \\ \hline
$0.9$ & 3.33e--2 & \cellcolor[HTML]{E5E3E3}1.63e--2 & 3.33e--2a& \cellcolor[HTML]{E5E3E3}1.63e--2s & 1.63e--1s & \cellcolor[HTML]{E5E3E3}3.33e--2a & 1.63e--1s & \cellcolor[HTML]{E5E3E3} nat-C. & 2.93e--2s & \cellcolor[HTML]{E5E3E3}non-C. \\ \hline
$0.8$ & 2.08e--1 & \cellcolor[HTML]{E5E3E3}6.52e--2 & 2.08e--1a& \cellcolor[HTML]{E5E3E3}6.52e--2s & 6.52e--2s& \cellcolor[HTML]{E5E3E3}2.07e--1a& 6.52e--2s & \cellcolor[HTML]{E5E3E3}2.04e--1a & 6.51e--2s & \cellcolor[HTML]{E5E3E3}nat-C. \\ \hline
$0.7$ & 1.01e+0 & \cellcolor[HTML]{E5E3E3}2.59e--1 & 1.01e+0a& \cellcolor[HTML]{E5E3E3}2.59e--1s & 2.59e--1s & \cellcolor[HTML]{E5E3E3}1.01e+0a& 2.59e--1s & \cellcolor[HTML]{E5E3E3}9.76e--1a & 2.24e--1s & \cellcolor[HTML]{E5E3E3}nat-C. \\ \hline
$0.6$ & 4.48e+0 & \cellcolor[HTML]{E5E3E3}1.70e+0 & 4.48e+0a& \cellcolor[HTML]{E5E3E3}1.02e+0s& 1.02e+0s & \cellcolor[HTML]{E5E3E3}4.43e+0a& 1.02e+0s & \cellcolor[HTML]{E5E3E3}4.03e+0a& 9.90e--1s & \cellcolor[HTML]{E5E3E3}nat-C.\\ \hline
$0.5$ & 1.88e+1 & \cellcolor[HTML]{E5E3E3}3.92e+0 & 1.87e+1a& \cellcolor[HTML]{E5E3E3}3.92e+0s & 3.92e+1s & \cellcolor[HTML]{E5E3E3}non-C. & 3.87e+0s & \cellcolor[HTML]{E5E3E3}non-C. &3.50e+0s & \cellcolor[HTML]{E5E3E3}nat-C. \\ \hline
\end{tabular}
}

\end{center}
\end{table}

\begin{figure}
\centering
\begin{subfigure}[b]{0.49\textwidth}
\centering
\includegraphics[width=0.85\textwidth]{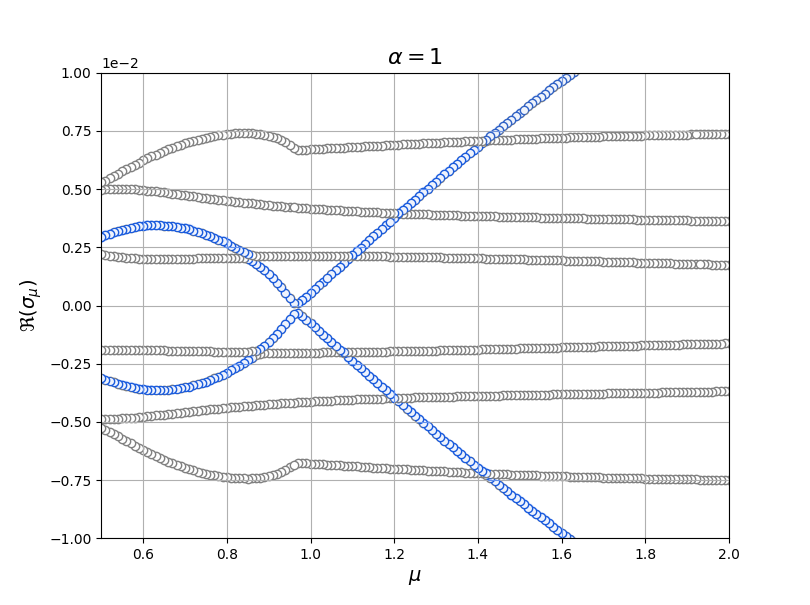}
\caption{}
\label{fig:in_eig1}
\end{subfigure}
\hfill
\begin{subfigure}[b]{0.49\textwidth}
\centering
\includegraphics[width=0.85\textwidth]{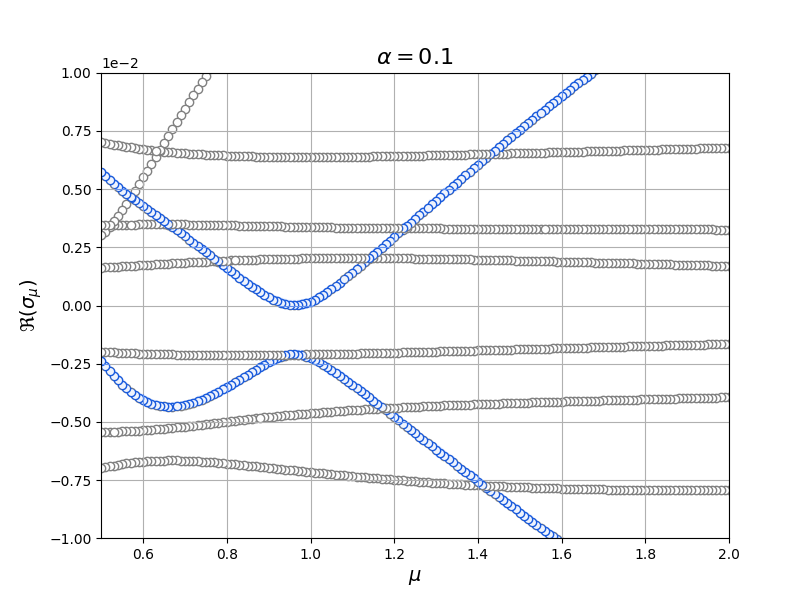}
\caption{}
\label{fig:in_eig_1e1}
\end{subfigure}\\
\begin{subfigure}[b]{0.49\textwidth}
\centering
\includegraphics[width=0.85\textwidth]{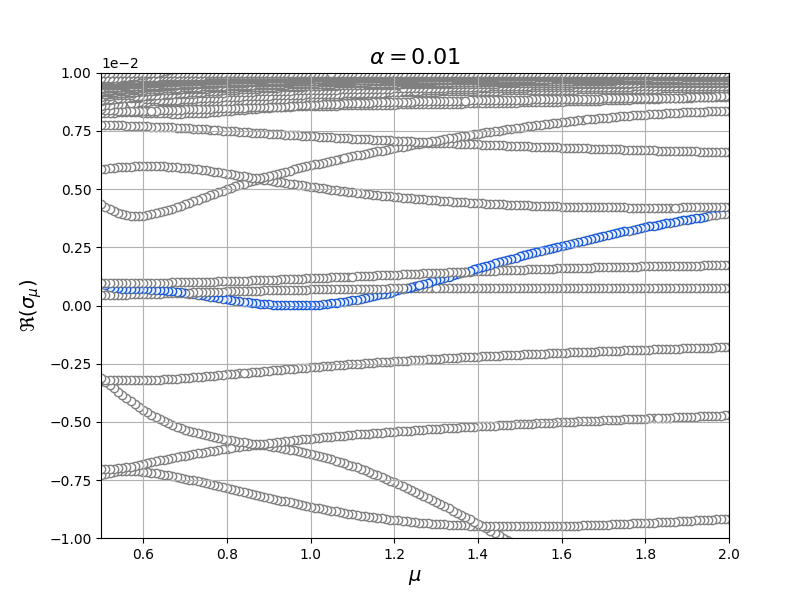}
\caption{}
\label{fig:in_eig_1e2}
\end{subfigure}
\hfill
\begin{subfigure}[b]{0.49\textwidth}
\centering
\includegraphics[width=0.85\textwidth]{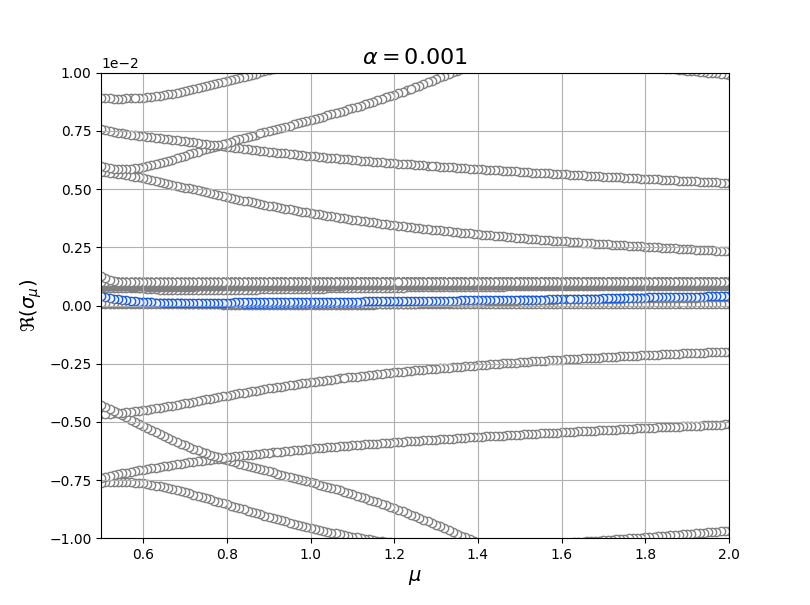}
\caption{}
\label{fig:in_eig_1e3}
\end{subfigure}\\
\caption{\emph{Channel Control}: spectral analysis with $\alpha = 1, 0.1, 0.01, 0.001$ following the natural branch.}
\label{fig:eig_inlet}
\end{figure}

Figure \ref{fig:eig_inlet} shows the eigenvalues of the global eigenproblem
in the range $\Re(\sigma_{\mu}) = [-0.01, 0.01]$ when following the natural optimal branch. %We remark that the eigenvalues behavior is preserved also for the other non-natural optimal branch. 
For $\alpha = 1$, we can observe the shears phenomenon, which disappears for other values of the penalization parameter. Lowering the value of $\alpha$, leads to a positive-dominated eigenvalues ensemble. 
Furthermore, a clustering around the value of $\alpha$ can be observed in plots \ref{fig:in_eig_1e2} and \ref{fig:in_eig_1e3}. In the next Section, very peculiar features have been observed as well, while changing the value of the penalization parameter $\alpha$ in a Dirichlet control.

\subsection{Dirichlet Control: flux action.}
\label{dirichlet}
In this final example, we propose a Dirichlet control over the boundary $\Omega_u \equiv \Gamma_{\text{D}}$. We fix the symmetric configuration as desired state $v_\text{d}$ over the line $\Gamma_{\text{obs}}$, while we set $\Gamma_{\text{wall}} = \Gamma_0$. In other words, we are trying to control a Dirichlet boundary condition in order to lead the controlled solution towards the symmetric profile. The problem to be solved reads: given $\mu \in \Cal P$ find $X \in \mathbb X$ such that
\begin{equation}
\label{eq:Dirichlet_eq}
\begin{cases}
v\mathbb{I}_{\Gamma_{\text{obs}}} -\mu \Delta w - v\cdot\nabla w + (\nabla v)^T w + \nabla q= v_\text{d} \mathbb{I}_{\Gamma_{\text{obs}}} \quad &\text{in} \ \Omega, \\
\nabla \cdot w = 0 \quad &\text{in} \ \Omega, \\
w =0 \quad &\text{on} \ \Gamma_{\text{in}} \cup \Gamma_{\text{D}} \cup \Gamma_{\text{wall}}, \\
- qn + (\mu \nabla w) n = 0 \quad &\text{on} \ \Gamma_{\text{out}}, \\
\alpha u  = w \quad &\text{on} \ \Gamma_D, \\
-\mu \Delta v + v\cdot\nabla v + \nabla p=0 \quad &\text{in} \ \Omega, \\
\nabla \cdot v = 0 \quad &\text{in} \ \Omega, \\
v = v_{\text{in}} \quad &\text{on} \ \Gamma_{\text{in}}, \\
v = u \quad &\text{on} \ \Gamma_{\text{D}}, \\
v = 0 \quad &\text{on} \ \Gamma_{\text{wall}}, \\
- pn + (\mu \nabla v) n = 0 \quad &\text{on} \ \Gamma_{\text{out}}.\\
\end{cases}
\end{equation}

\begin{figure}
\centering
\includegraphics[width=0.489\textwidth]{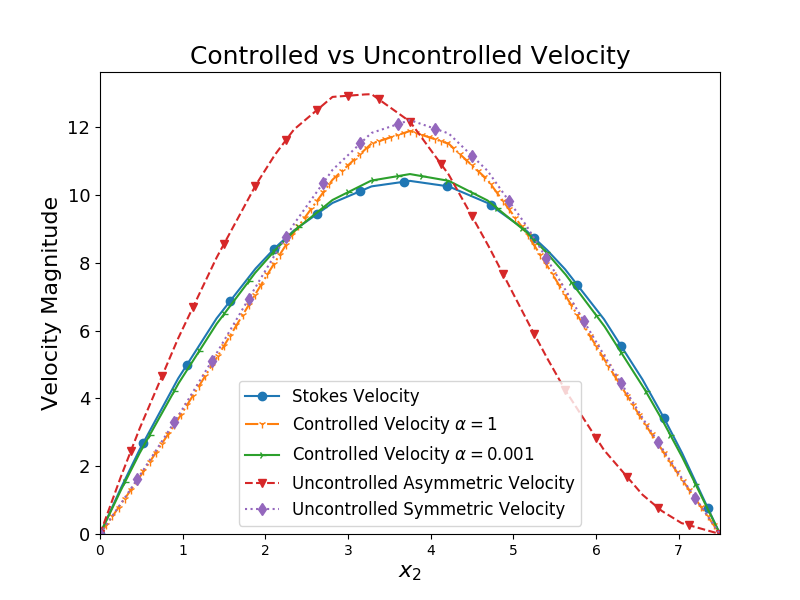}
\includegraphics[width=0.46\textwidth]{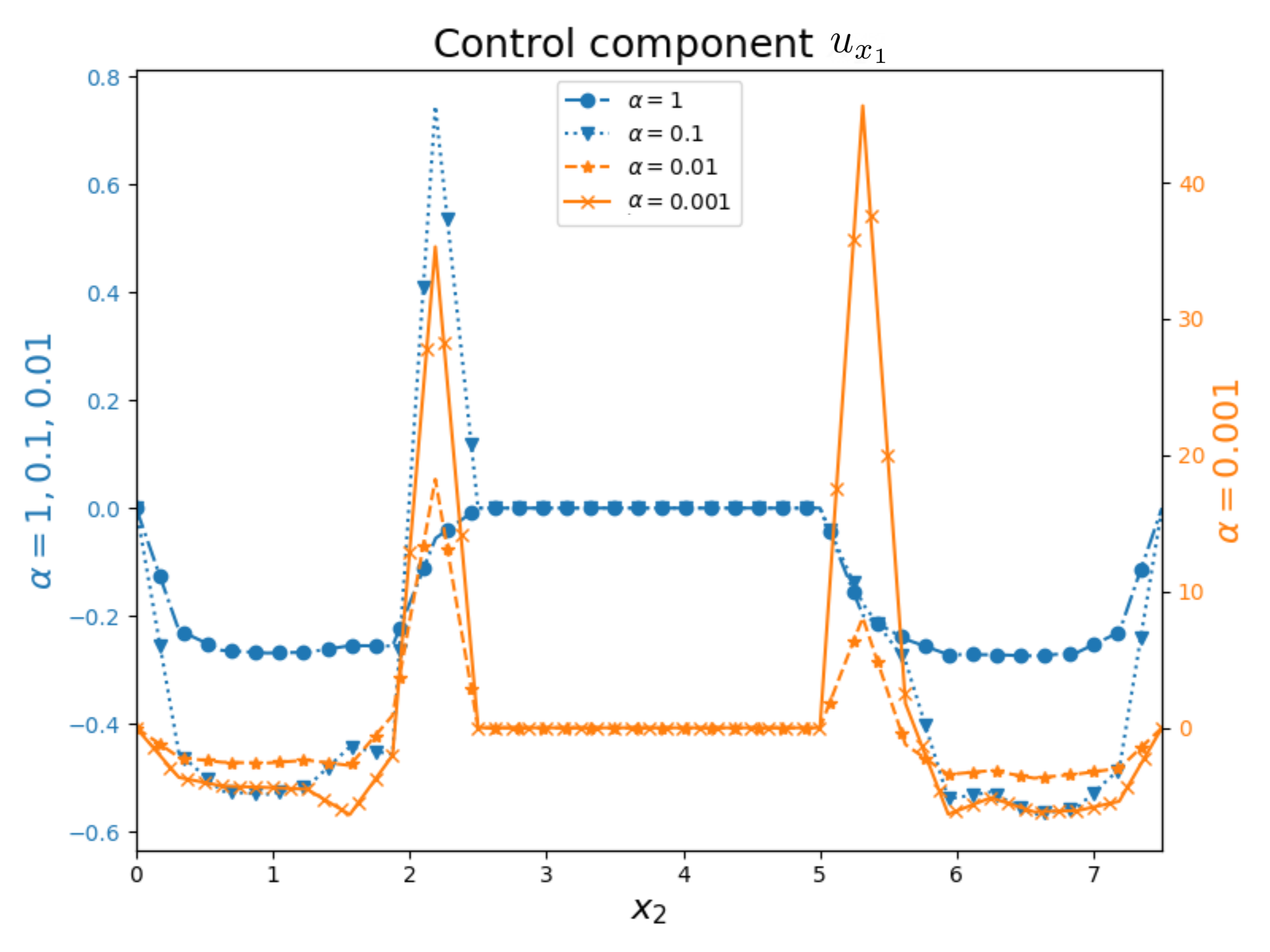}
\caption{\emph{Dirichlet Control}. \emph{Left}: comparison of velocity profiles in the controlled and uncontrolled cases for $\alpha = 1, 0.01$, $\mu = 0.5$ on $\Gamma_{\text{obs}}$ w.r.t. the symmetric desired profile when following the natural optimal branch. \emph{Right}: representation of control variable evolution for $\alpha=0.1, 0.01, 0.001, 0.001$ and $\mu = 0.5$ at $x_1 = 10$ when following the natural optimal branch.}
\label{fig:diri_mag}
\end{figure}

\begin{figure}
\centering
\begin{subfigure}[b]{0.49\textwidth}
\centering
\includegraphics[width=\textwidth]{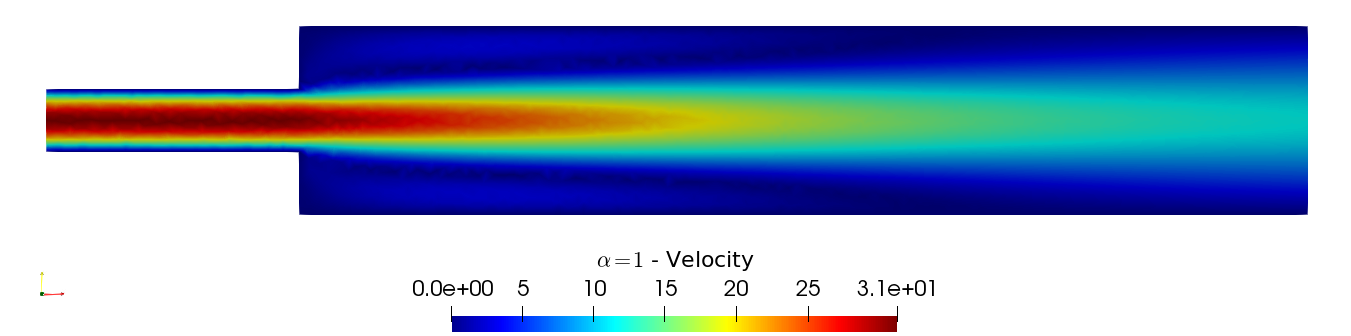}
\caption{}
\label{fig:diri_v_1}
\end{subfigure}
\hfill
\begin{subfigure}[b]{0.49\textwidth}
\centering
\includegraphics[width=\textwidth]{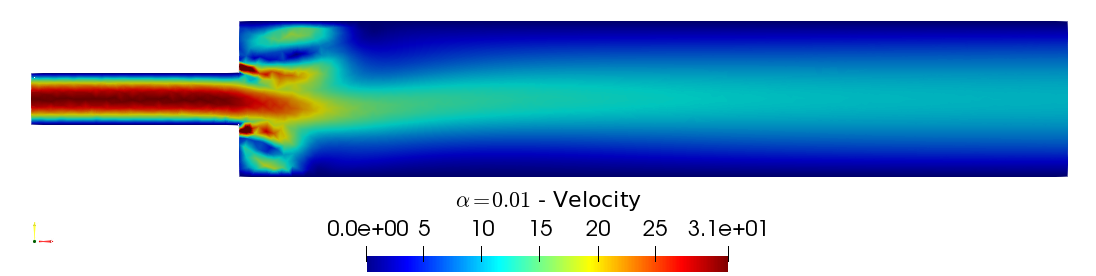}
\caption{}
\label{fig:diri_v_1e3}
\end{subfigure}
\caption{\emph{Dirichlet Control}: two optimal velocity solutions for $\mu=0.5$, with $\alpha = 1$ and $\alpha = 0.001$, left and right, respectively.}
\label{fig:Dirichlet_solution}
\end{figure}

Allowing the flux to freely enter or exit from the boundary $\Gamma_{\text{D}}$ drastically changes the optimal solution behavior. Since we are asking for a symmetric desired profile, the main action of the control is to straighten the flow: this behavior can be observed from Figure \ref{fig:diri_v_1} and the left plot of Figure \ref{fig:diri_mag}. Indeed, even for large values of $\alpha$, the velocity profile reaches the symmetric configuration, while for lower values of the penalization parameter, the velocity on $\Gamma_{\text{obs}}$ is parabolic.
This feature is highlighted also from the functional values in Table \ref{Diri_J}, where the functional \eqref{eq:J_NS} is shown for several $\mu$ (rows) and $\alpha$ (columns) w.r.t. the uncontrolled stable and unstable solutions. The cost functional largely decreases for smaller values of $\alpha$, e.g. $\alpha = 0.001$: for example, focusing on $\mu = 0.5$, the functional only lowers of $18\%$ for $\alpha = 0.01$, in contrast to almost $82\%$ for $\alpha = 0.001$. Within the setting of $\alpha = 0.001$, the system manifests an interesting and unexpected profile, shown in Figure \ref{fig:diri_v_1e3}. The flux presents an asymmetric configuration for low value of $\mu$. Namely, for low $\alpha$ a bifurcating solution appears as depicted in Figure \ref{fig:diri_v_1e3}. The asymmetric behavior is due to the control variable which not only allows the flow to exit from $\Gamma_{\text{D}}$ (in order to avoid the asymmetric recirculation of the wall-hugging solution), but it is adding flux near the channel, in order to achieve the straight configuration and the parabolic velocity profile given by the symmetric target velocity over the observation domain, as it is represented in the right plot of Figure \ref{fig:diri_mag}.

\begin{figure}
\centering
\begin{subfigure}[b]{0.49\textwidth}
\centering
\includegraphics[width=0.85\textwidth]{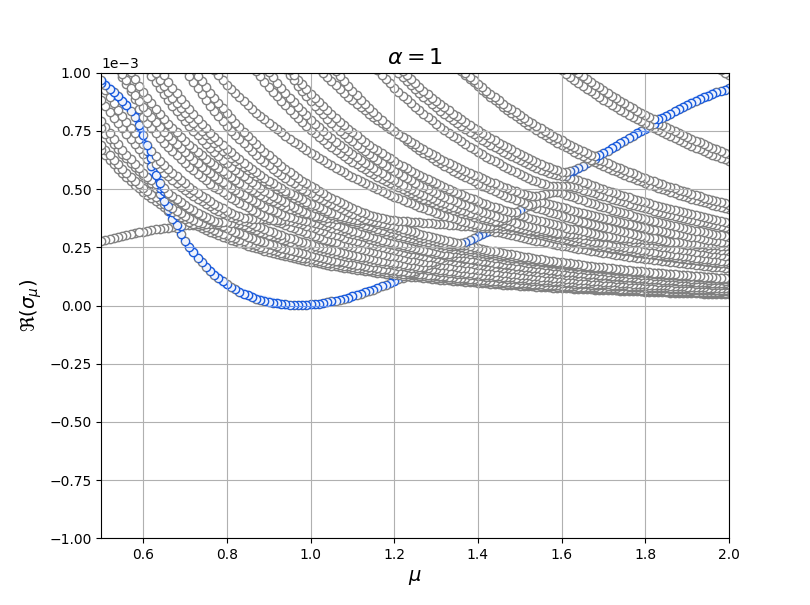}
\caption{}
\label{fig:diri_eig1}
\end{subfigure}
\hfill
\begin{subfigure}[b]{0.49\textwidth}
\centering
\includegraphics[width=0.85\textwidth]{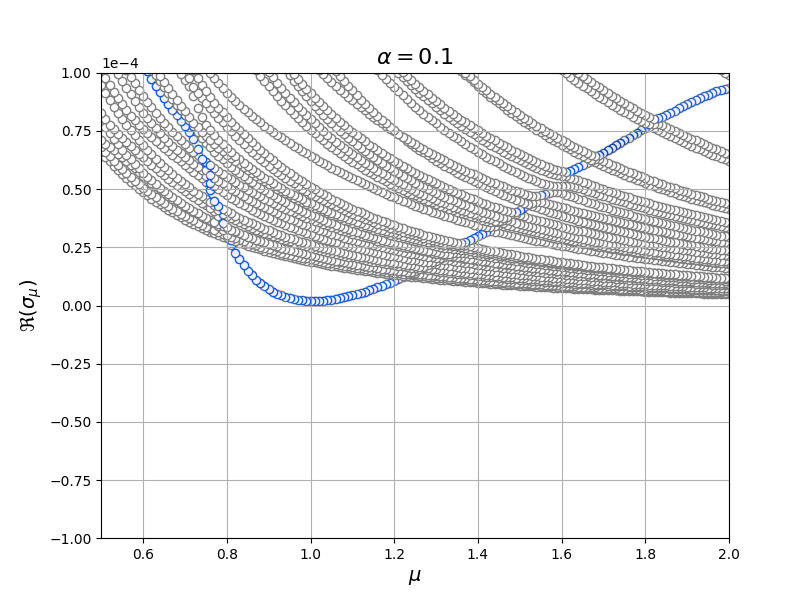}
\caption{}
\label{fig:diri_eig_1e1}
\end{subfigure}\\
\begin{subfigure}[b]{0.49\textwidth}
\centering
\includegraphics[width=0.85\textwidth]{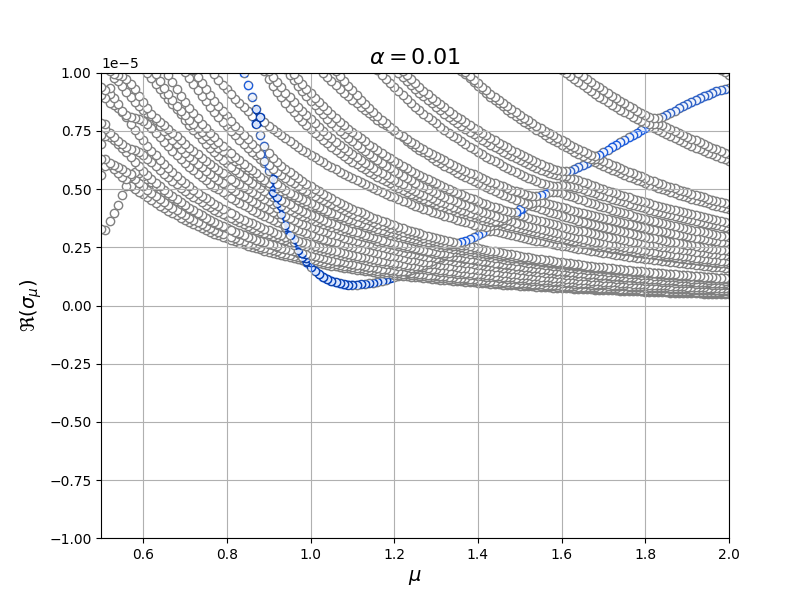}
\caption{}
\label{fig:diri_eig_1e2}
\end{subfigure}
\hfill
\begin{subfigure}[b]{0.49\textwidth}
\centering
\includegraphics[width=0.85\textwidth]{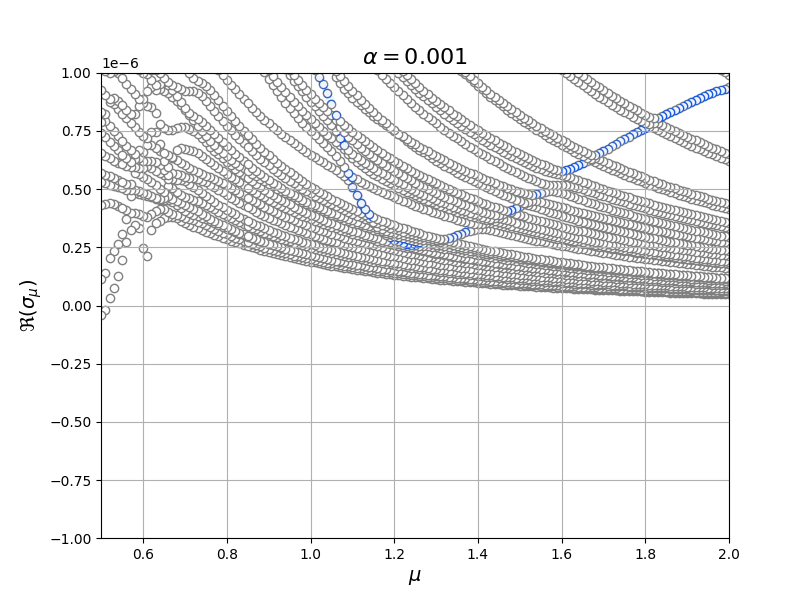}
\caption{}
\label{fig:diri_eig_1e3}
\end{subfigure}\\
\caption{\emph{Dirichlet Control}: spectral analysis with $\alpha = 1, 0.1, 0.01, 0.001$.}
\label{fig:eig_diri}
\end{figure}

\begin{table}[b]
\caption{\emph{Dirichlet Control}: comparison of the functional value w.r.t. the stable and unstable uncontrolled solutions.}
\label{Diri_J}
\begin{center}
\footnotesize{
\begin{tabular}{|c||c|c||c|c||c|c||c|c||c|c||}
\hline
\multirow{2}{*}{$\mu$} &Stable & \cellcolor[HTML]{E5E3E3} Unstable & \multicolumn{8}{c||}{Controlled Solution} \\ \cline{2-11}
& \multicolumn{2}{c||}{Uncontrolled} & \multicolumn{2}{c||}{$\alpha = 1$} & \multicolumn{2}{c||}{$\alpha = 0.1$} & \multicolumn{2}{c||}{$\alpha = 0.01$} & \multicolumn{2}{c||}{$\alpha = 0.001$} \\ \hline
$2$ & { 5.14e--9 } & \cellcolor[HTML]{E5E3E3} { 5.14e--9 } & \multicolumn{2}{c||}{ 4.98e--9 } & \multicolumn{2}{c||}{ 4.83e--9 } & \multicolumn{2}{c||}{ 4.79e--9 } & \multicolumn{2}{c||}{ 4.79e--9 } \\ \hline
$1.5$ & { 4.38e--6 } & \cellcolor[HTML]{E5E3E3} { 4.38e--6 } & \multicolumn{2}{c||}{ 4.24e--6 } & \multicolumn{2}{c||}{ 4.10e--6} & \multicolumn{2}{c||}{ 4.07e--6} & \multicolumn{2}{c||}{ 4.06e--6} \\ \hline
$1$ & { 4.10e--3 } & \cellcolor[HTML]{E5E3E3} { 4.10e--3 } & \multicolumn{2}{c||}{ 3.94e--3 } & \multicolumn{2}{c||}{ 3.78e--3} & \multicolumn{2}{c||}{ 3.74e--3 } & \multicolumn{2}{c||}{ 3.72e--3} \\ \hline
$0.9$ & { 3.33e--2 } & \cellcolor[HTML]{E5E3E3} { 1.63e--2 } & \multicolumn{2}{c||}{ 1.56e--2 } & \multicolumn{2}{c||}{ 1.49e--2} & \multicolumn{2}{c||}{ 1.47e--2 } & \multicolumn{2}{c||}{ 1.45e--2 } \\ \hline
$0.8$ & { 2.08e--1 } & \cellcolor[HTML]{E5E3E3} { 6.52e--2} & \multicolumn{2}{c||}{ 6.20e--2 } & \multicolumn{2}{c||}{ 5.88e--2} & \multicolumn{2}{c||}{ 5.78e--2} & \multicolumn{2}{c||}{ 5.46e--2} \\ \hline
$0.7$ & { 1.01e+0} & \cellcolor[HTML]{E5E3E3} { 2.69e--1 } & \multicolumn{2}{c||}{ 2.44e--1 } & \multicolumn{2}{c||}{ 2.29e--1} & \multicolumn{2}{c||}{ 2.21e--1 } & \multicolumn{2}{c||}{ 1.82e--1 } \\ \hline
$0.6$ & { 4.48e+0} & \cellcolor[HTML]{E5E3E3} { 1.70e+0} & \multicolumn{2}{c||}{ 9.49e--1 } & \multicolumn{2}{c||}{ 8.73e--1} & \multicolumn{2}{c||}{ 8.09e--1 } & \multicolumn{2}{c||}{ 3.57e--1 } \\ \hline
$0.5$ & { 1.88e+1 } & \cellcolor[HTML]{E5E3E3} { 3.92e+0 } & \multicolumn{2}{c||}{ 3.58e+0 } & \multicolumn{2}{c||}{ 3.21e+0} & \multicolumn{2}{c||}{ 2.41e+0 } & \multicolumn{2}{c||}{ 4.73e--1 } \\ \hline
\end{tabular}
}
\end{center}
\end{table}

\begin{figure}
\centering
\includegraphics[width=0.489\textwidth]{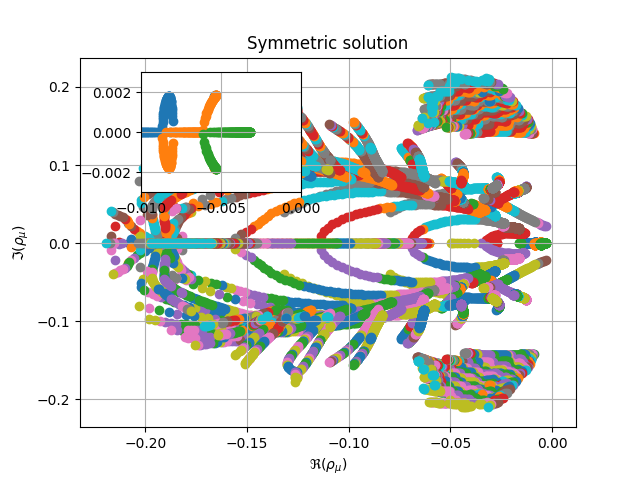}
\includegraphics[width=0.489\textwidth]{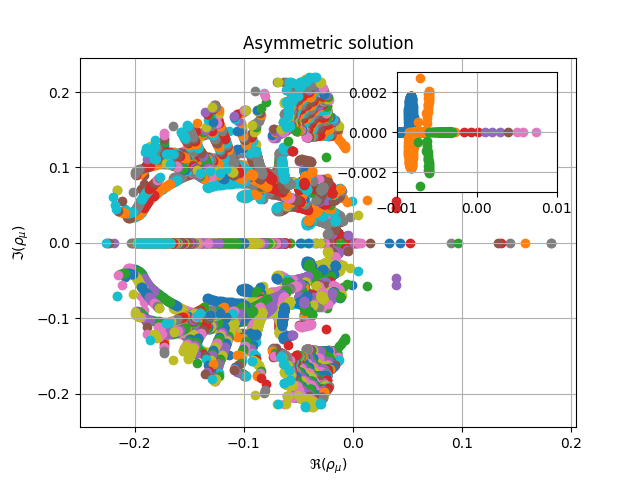}
\caption{\emph{Dirichlet Control}. Eigenvalues of the state eigenproblem in the complex plane: asymmetric and symmetric solutions, left and right, respectively.}
\label{fig:dirichlet_eig_state}
\end{figure}
Figure \ref{fig:eig_diri} provides the eigenvalue analysis, where we show some close-ups starting with $\Re(\sigma_{\mu}) = [-0.001, 0.001]$ for $\alpha = 1$ in the top-right image, and restricting the vertical interval due to the order of the lowered $\alpha$ in the remaining images. As already noticed in Section \ref{distributed}, the strongest is the control, the more the eigenvalues are positive.
Furthermore, as it already happened in the distributed control case with asymmetric target in Section \ref{distributed}, the value $\mu^*$ seems to be not relevant anymore as an indication of the bifurcation point. Recalling the definition of $\mu^{**}$ in Section \ref{distributed} as the value of the parameter for which the top curve of the shears (marked in blue) is approaching $\Re(\sigma_\mu) = 0$ from above, Figure \ref{fig:eig_diri} shows that such a curve is moving away from $\Re(\sigma_\mu) = 0$ as $\alpha$ decreases, and thus no such point $\mu^{**}$ exists. 
The results of the previous Sections have shown that the top shear structure is typically associated to the wall-hugging bifurcation, and that $\mu^{**}$ provides an indication of the bifurcation: we are thus lead to believe that the standard bifurcating configuration observed in the uncontrolled case, consisting of a branch of symmetric solutions and a branch of wall-hugging ones, is not present here, with the latter branch disappearing. However, the system seems to be featuring a \emph{different bifurcation}, presented in Figure \ref{fig:diri_v_1e3}. Indeed, we can see an eigenvalue crossing the line $\Re(\sigma_\mu) = 0$ for the global eigenproblem in Figure \ref{fig:diri_eig_1e3} for $\alpha = 0.001$.
Therefore, if we plot the eigenvalues of the state eigenproblem of Algorithm \ref{alg:01} (Figure \ref{fig:dirichlet_eig_state}), we see how the symmetric profile does never cross the origin, while the asymmetric solution in Figure \ref{fig:diri_v_1e3} for $\alpha=0.001$ does. 
In the setting with the modified BCs, the physical stable solution behavior is a feature of straight profile.
Moreover, from Figure \ref{fig:dirichlet_eig_state}, we can clearly observe a couple of complex and conjugate eigenvalue crossing the imaginary axis. This is, in fact, a paradigm for Hopf bifurcation \cite{AQpreprint,seydel2009practical} and represents another evidence of how deeply the system changed its inner features.

\begin{remark}[Lagrange multipliers]
From a numerical point of view, we employed Lagrange multipliers to solve the optimality system \eqref{eq:Dirichlet_eq}. The condition $v = u$ on $\Gamma_{\text{D}}$ has been weakly imposed in integral form.
% \begin{equation}
% \label{eq:multiplier1}
% \int_{\Gamma_{\text{D}}} v \lambda \; ds = \int_{\Gamma_{\text{D}}} u \lambda \; ds \quad \forall \lambda \in \V.
% \end{equation}
% This equation reflects in system \eqref{eq:gal_ocp_ns2}, since the term
% \begin{equation}
% \label{eq:neumann_gamma_D}
% \int_{\Gamma_{\text{D}}} \phi (\nabla \cdot v n + pn) \; ds \quad \forall \phi \in \V,
% \end{equation}
% appears.
This will result in extra terms in the adjoint equations. Furthermore, we weakly impose also the boundary condition $w = 0$ on $\Gamma_{\text{D}}$ with another multiplier. The reason of this latter decision will be explained in Section \ref{sec_ROM}.
\end{remark}

\subsection{Comparative Eigenvalue Analysis}
\label{comparison}
In this Section, we sum up all the observations and results derived from the global eigenvalue analysis over the four test cases. Therefore, we now list the similarities between these, especially for what concerns the variation against different values of the penalization parameter $\alpha$:
\begin{itemize}
\item[{$\small{\circ}$}] the \emph{ eigenvalues cluster} around the value of $\alpha$. This behavior is well represented in Figures \ref{fig:n_eig_1e2}, \ref{fig:n_eig_1e3}, \ref{fig:in_eig_1e2} and \ref{fig:in_eig_1e3}. These eigenvalues come from the optimality equation;
\item[{$\small{\circ}$}] the \emph{predominance of positive eigenvalues} over the negative ones. In all the performed spectral analysis we have observed that the control action lowers the negative eigenvalues. The stronger is the control, the greater is the number of positive eigenvalues, as it is represented in Figures \ref{fig:dist_eig_1e1} and \ref{fig:eig_diri};
\item[{$\small{\circ}$}] the \emph{shears effect} for low controlled system. The shears configuration is characteristic of control problems which do not highly change the uncontrolled system solution. It is the case of Neumann control in Figure \ref{fig:n_eig1} and of the channel control for $\alpha = 1$ as shown in Figure \ref{fig:in_eig1}. For the other cases, the smaller is $\alpha$ the least visible is this eigenvalue configuration: in some cases, the structure is completely broken;
\item[{$\small{\circ}$}] the $\mu^{**}$ \emph{identification}. It is clear that the shears (or their top curve, if shears are broken) approach to the real line in the point $\mu^{**}$ for which the bifurcating phenomenon of the controlled system occurs. This is the same situation we found for the uncontrolled problem, in which the path of the eigenvalue identifies the value of bifurcation parameter $\mu^*$. Moreover, this situation is often preserved regardless of $\alpha$. Indeed, the positive shears eigenvalue is still present in Figures \ref{fig:dist_eig_1e1}, \ref{fig:in_eig_1e1}, \ref{fig:in_eig_1e2} and even in the Dirichlet optimal control, as shown in Figure \ref{fig:eig_diri}. In some cases, a shift of the $\mu^{**}$ compared to $\mu^*$ has been observed.
\end{itemize}

Since the structure of the spectral analysis is highly influenced by the control strength, we tried to perform an eigenvalue analysis dealing with only state and adjoint equations. For all the test cases, shears occurs.
The shears structure is symmetric when the solution
shows the wall-hugging property, while it is slightly asymmetric when the state flow is straight.
We guess
that this behavior is due to the different reaction of state and adjoint blocks to the bifurcating phenomena. Indeed, for the symmetric flux, the behavior of the state equation has to be preserved for all $\mu$, while the adjoint problem, which is strictly linked to the control variable, puts more effort in rebalancing the flux, resulting in an asymmetric contribution that causes the shears to be slightly asymmetric.

The spectral analysis of a nonlinear system is a indispensable tool to understand bifurcation phenomena, eventually. Under the point of view of computational costs, it is a very tough task, most of all for nonlinear \ocp s. Indeed, FE discretization leads to huge systems to be solved for a wide sample of parameter $\mu \in \Cal P$. In the next Section we propose ROMs as a suitable approach to overcome this issue.

\begin{figure}[b]
\centering
\includegraphics[width=0.489\textwidth]{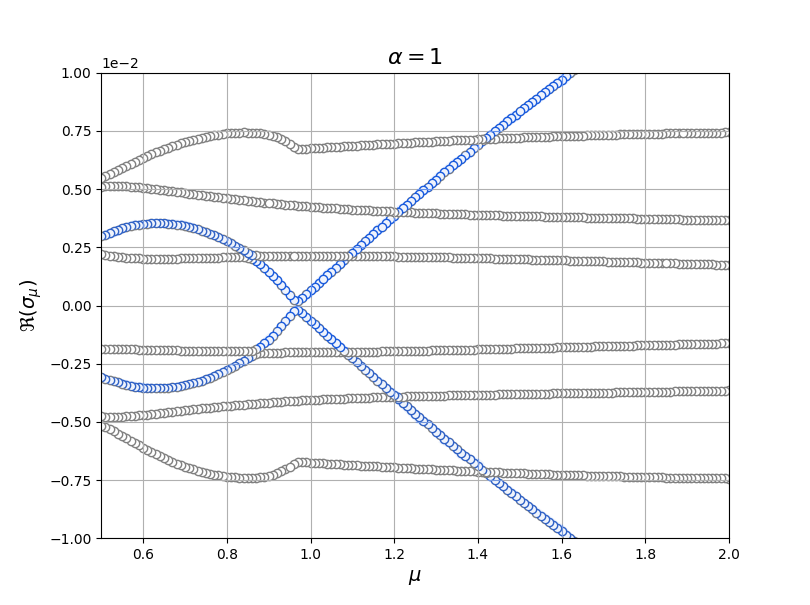}
\includegraphics[width=0.489\textwidth]{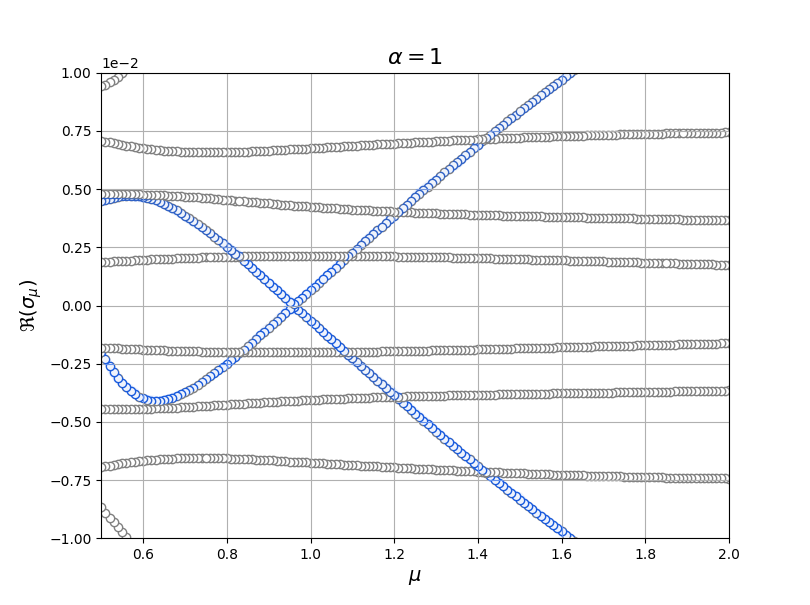}
\caption{\emph{Comparative Analysis}. \emph{Left}: asymmetric velocity with Neumann control for $\alpha = 1$. \emph{Right}: symmetric velocity with distributed control for $\alpha = 1$. }
\label{fig:tenaglie}
\end{figure}

\section{ROMs for Nonlinear \ocp s}
\label{sec_ROM}

This Section introduces ROM approximation techniques for nonlinear \ocp s. The proposed reduced strategy is independent from the governing state equation. Indeed, we refer to \cite {Strazzullo1, Strazzullo3, ZakiaMaria, Zakia} for previous contributions to ROM for nonlinear \ocp s, extending them to standard techniques proposed in \cite{HESS2019379, Hess2019, pichiquaini, pichirozza,pintore2019efficient, PR15}, for bifurcating systems. Section \ref{sec_rom_gen} introduces the basics ideas of ROM approach and the standard assumptions to guarantee its efficiency and applicability. Moreover, in Section \ref{POD}, we will describe the used reduction strategy, relying on POD-Galerkin basis construction, see \cite{ballarin2015supremizer, burkardt2006pod, Chapelle2013, hesthaven2015certified} as introductory references, combined with aggregated spaces techniques, following the linear \ocp s fashion, as already presented in \cite{bader2016certified,bader2015certified,dede2010reduced,gerner2012certified,karcher2014certified,karcher2018certified, negri2015reduced,negri2013reduced,quarteroni2007reduced}. Finally, numerical results are shown in Section \ref{rom_results}.

\subsection{General Reduction Strategy}
\label{sec_rom_gen}
In Section \ref{NS_ocp}, we analyzed how optimal flow control can modify an expected behavior in bifurcating systems. In such problems, a parametric study of the state solution is necessary in order to understand the solution properties. The FE method cannot be exploitable when several instances of parametrized problem have to be studied. Indeed, the computational time required can be too costly, especially in an optimal control setting, where the state equation is associated to an optimization system. ROM techniques aim at building a \emph{reduced} surrogate of the FE approximation. The main idea of ROM is to spend computational resources to build a new (smaller) model starting from FE solutions. Such a \emph{reduced system}, although having a considerably lower dimension, is built in a way that it does not lose accuracy and can be used to analyze several parametric configurations in a versatile low-dimensional framework.
We now briefly introduce ROM ideas in the \ocp s setting. At the continuous level, we have already defined the solution branch $\Cal X_i$, which represents how an optimal solution
$X(\boldsymbol \mu) = (y(\boldsymbol \mu), u(\boldsymbol \mu),z(\boldsymbol \mu))$ of \eqref{ocp} changes w.r.t.\ the parameter $\bmu \in \Cal P$. For the sake of clarity, in this Section we will underline the parameter dependence of our solution branch, since it is of primary importance to understand the ROM basics and notions. From now on we will refer to FE as \emph{high-fidelity} approximation. Indeed, we suppose that the FE discretization reliably represents the solution branch $\Cal X_i$ as its discretized counterpart $\Cal X_i^{\Cal N}$. \\
The ROM aims at representing the high-fidelity $\Cal X_i^{\Cal N}$ through the construction of basis derived from \emph{snapshots}, i.e. properly chosen solutions of \eqref{FE_ocp}. The reduced function spaces are contained in the FE spaces: a standard Galerkin projection is performed, resulting in an efficient low-dimensional solution which is \A{still accurate w.r.t.\ }the high-fidelity model. We exploited a \emph{branch-wise reduction}, namely, for every bifurcating solution branch $\Cal X_i$, we build a different ROM. \\
Of course this is the best approach form the accuracy standpoint, but in \cite{pichirozza,pichiquaini} a different context a global approach was pursued, where the loss in accuracy is balanced with the construction of a single ROM.
We suppose to have already constructed the reduced spaces $\state_N \subset \state {\discy} \subset \state$ and $\control{_N} \subset \control \discu \subset \control$, the former for state and adjoint variables, the latter for control, respectively (description of the construction of reduced spaces is postponed to Section \ref{POD}). The \emph{function space dimension} $N$ is usually much lower than $\Cal N$. We will refer to this stage as \emph{offline phase}: here, apart from the basis construction, all the parameter independent quantities are assembled and stored. \\
After this reduced spaces building process, one can solve the following low-dimensional problem in an \emph{online phase}
given $\boldsymbol \mu \in \Cal P$, find
$X_N (\boldsymbol \mu) = (y_N(\boldsymbol \mu), u_N(\boldsymbol \mu),z_N(\boldsymbol \mu)) \in
\mathbb X_N \eqdot
\Cal \state_N \times \Cal \control_N \times \Cal \state_N$ such that it holds:
\begin{equation}
\label{ROM_ocp}
\begin{cases}
D_{y}\Lg(X_N; y_\text{d}, \bmu)[\omega] = 0 & \forall \omega \in \state{_N},\\
D_u\Lg(X_N; y_\text{d}, \bmu)[\kappa] = 0 & \forall \kappa \in \control{_N},\\
D_z\Lg(X_N; y_\text{d}, \bmu)[\zeta] = 0 & \forall \zeta \in \state{_N}.\\
\end{cases}
\end{equation}
Namely, at each new parametric instance $\bmu \in \Cal P$, the system \eqref{ROM_ocp}, which inherits the features from the high-fidelity dynamics \eqref{FE_ocp}, is assembled and solved. Also in this case, we can deal with the nonlinearity applying a Newton's method, as in the FE setting. Moreover, the stability analysis can be performed as described in Algorithm \ref{alg:01}, employing the Galerkin projection in the reduced space. It is clear that, it is convenient to exploit the ROM only if one does not have to build from scratch the reduced model for any parametric instance.
For this reason, the system \eqref{ROM_ocp} is assumed to be affinely decomposed, i.e.\ all the variational forms involved can be written as the product of $\boldsymbol \mu -$independent forms and
$\boldsymbol \mu -$dependent functions \cite{hesthaven2015certified}. %Focusing on the system \eqref{ROM_ocp}.
%this translates in the following finite sums:
%\begin{comment}
%\begin{equation}
%\label{affinity}
%\begin{aligned}
%D_{y}\Lg(X_N, y_\text{d}, \bmu)[\omega] =
%\displaystyle \sum_{\mathsf q=1}^{Q_y} \Theta_{{y}}^\mathsf q(\boldsymbol{\mu})D_{y}\Lg^\mathsf q(X_N; y_\text{d})[\omega],\\
%D_{u}\Lg(X_N; y_\text{d}, \bmu)[\kappa], =
%\displaystyle \sum_{\mathsf q=1}^{Q_u} \Theta_{{u}}^\mathsf q(\boldsymbol{\mu})D_{u}\Lg^\mathsf q(X_N; y_\text{d})[\kappa],\\
%D_{z}\Lg(X_N; y_\text{d}, \bmu)[\zeta] =
%\displaystyle \sum_{\mathsf q=1}^{Q_z} %\Theta_{{z}}^\mathsf %q(\boldsymbol{\mu})D_{z}\Lg^\mathsf q(X_N; %y_\text{d})[\zeta],
%\end{aligned}
%\end{equation}
%where $\Theta_{{y}}^\mathsf q(\boldsymbol{\mu}), \Theta_{u}^\mathsf q(\boldsymbol{\mu})$ and $\Theta_{z}^\mathsf q(\boldsymbol{\mu})$ are $\boldsymbol{\mu}-$dependent smooth functions, while
%$D_{y}\Lg^\mathsf q(X_N, y_\text{d})[\omega],$ %\\$D_{u}\Lg^\mathsf q(X_N; y_\text{d})[\kappa]$ and $D_{z}\Lg^\mathsf q(X_N; y_\text{d})[\zeta]$ are $\boldsymbol{\mu} -$independent bilinear forms describing the optimality system. \end{comment}
When the affine dependency assumption is verified, the online phase does not depend on $\Cal N$, and can be performed in a small amount of time. Conversely, the offline process is performed only once and can take advantage of High Performance Computing (HPC) resources.
\begin{remark}
Our test cases deal with Navier-Stokes governing equations \eqref{eq:NS_eq}, then, it has at most quadratically
nonlinear terms and the affine decomposition is not fulfilled. One can employ hyper-reduction techniques as the Empirical Interpolation Method (EIM) to recover it. We refer the interested reader to \cite{barrault2004empirical} or \cite[Chapter 5]{hesthaven2015certified}.
\end{remark}

In the next Section, we will focus on the ROM offline and online phase, showing the strategy employed to build the reduced function spaces.

\subsection{Offline and online stages}
\label{POD}
In this Section, we present how to build a reduced space for \ocp s. We exploit here the POD algorithm. Thanks to this algorithm, $N_{\text{\text{max}}}$ snapshots are sampled and then compressed in order to generate function spaces of dimension $N < N_{\text{\text{max}}}$.

It is well known that optimization governed by PDEs($\bmu$) constraints leads to the solution of a saddle point system \cite{Benzi, bochev2009least, hinze2008optimization, Stoll}, as already specified in Section \ref{FE}. In order to guarantee the well-posedness of such a structure, the matrix $\mathsf B$ of system \eqref{J_saddle} must verify the inf-sup stability condition \eqref{FE_lbb} for every $\boldsymbol \mu \in \Cal P$.
In the FE approximation, the above-mentioned relation holds since state and adjoint spaces are equally discretized. However, the inf-sup stability must hold at the reduced level too, since the relation is provable if the reduced spaces for state and adjoint variables coincide. The standard POD construction process leads to the reduced function spaces for state and adjoint which may be different even under the assumption of the same starting FE spaces. To overcome this issue, the basis are usually manipulated in order to stabilize the system. Indeed, we apply \emph{aggregated spaces} technique, as already done in several papers about ROM for \ocp s, see \cite{bader2016certified,bader2015certified,dede2010reduced,gerner2012certified,karcher2014certified,karcher2018certified, negri2015reduced,negri2013reduced,quarteroni2007reduced} as references. The strategy aims at building a common space for state and adjoint which is able to describe both state and adjoint variables. \\Let us suppose to have applied a standard POD for all the involved variables and to have defined the following spaces $\label{state_r} {\state}_{N}= \text{span }\{\chi^{y}_n, \chi^{z}_n, \; n = 1, \dots, N\}$ and $\label{control_r}  {\control}_{N} = \text{span}\{\chi^{u}_n, \; n = 1, \dots, N\}$, with
$$
\mathsf Z =
\begin{bmatrix}
\mathsf Z_{\mathsf x} \\
\mathsf Z_{\mathsf z}
\end{bmatrix},
\spazio \text{and} \spazio
\mathsf Z_{\mathsf x} =
\begin{bmatrix}
\mathsf Z_{\mathsf y} \\
\mathsf Z_{\mathsf u}
\end{bmatrix}
$$
where
$
\mathsf Z_{\mathsf y} \equiv \mathsf Z_{\mathsf z} = [\chi_{1}^{y} | \cdots | \chi_{N}^{y}| \chi_{1}^{z} | \cdots | \chi_{N}^{z}] \in \mathbb R^{\Cal N_{y} \times 2N}
$ and $
\mathsf Z_{\mathsf u} = [\chi_{1}^{u} | \cdots | \chi_{N}^{u}] \in \mathbb R^{\Cal N_{u} \times N}
$ are the reduced basis matrices for each variable and $\mathsf Z $ spans the global space $\mathbb X_N$.
We want to solve the optimality system in a low dimensional framework at each parametric instance.
To this end, we employ a Galerkin projection into the reduced spaces and the system \eqref{ocp} will be
\begin{equation}
\label{G_compact_ROM}
\mathsf G_{N}(\mathsf X_N; \boldsymbol \mu) \mathsf X_N = \mathsf F_N,
\end{equation}
where
$$\mathsf G_{N}(\mathsf X_N; \boldsymbol \mu) \eqdot \mathsf Z^T \mathsf G(\mathsf Z \mathsf X_N; \boldsymbol \mu),
\spazio \text{and} \spazio \mathsf F_N \eqdot \mathsf Z^T \mathsf F.$$
The system \eqref{G_compact_ROM} is nonlinear, thus we apply Newton's method and we iteratively obtain
\begin{equation}
\mathsf {X}_N^{j + 1} = \mathsf {X}_N^j+ \mathsf{Jac}_N(\mathsf X_N^{j}; \boldsymbol \mu)^{-1}(\mathsf F_N - \mathsf G(\mathsf X_N^j; \boldsymbol \mu)\mathsf X_N^j), \spazio j \in \mathbb N.
\end{equation}
from the FE approximation, the Fr\'echet derivative inherits the saddle point structure, i.e.
\begin{equation}
\label{Frechet_ROM}
\mathsf{Jac}_N (\mathsf X_N; \boldsymbol \mu) \mathsf X_N =
\begin{bmatrix}
\mathsf A_N & \mathsf B_N^T \\
\mathsf B_N & 0 \\
\end{bmatrix}
\begin{bmatrix}
\mathsf x_N \\
\mathsf z_N
\end{bmatrix},
\end{equation}
with $\mathsf{Jac}_N (\mathsf X_N; \boldsymbol \mu) = \mathsf Z^T \mathsf{Jac}(\mathsf Z \mathsf X_N; \bmu)\mathsf Z$, $\mathsf A_N = \mathsf Z_{\mathsf x}^T \mathsf A \mathsf Z_{\mathsf x}$ and $\mathsf B_N = \mathsf Z_{\mathsf z}^T \mathsf B \mathsf Z_{\mathsf x}$. \\
We now have all the ingredients to define a \emph{reduced Brezzi inf-sup condition} as follows
\begin{equation}
\label{ROM_infsup}
\beta_{\text{Br},N}(\bmu) \eqdot \adjustlimits \inf_{0 \neq \mathsf z_N} \sup_{0 \neq \mathsf x_N} \frac{\mathsf z_N^T\mathsf B_N \mathsf x_N}{\norm{\mathsf x_N}_{\mathbb Y \times \mathbb U}\norm{\mathsf z_N}_{\Cal \state}} \geq \overline{\beta}_{\text{Br},N} > 0.
\end{equation}

If $\bmu \neq \bmu^{\ast}$, relation \eqref{ROM_infsup} is verified thanks to the aggregated space definition.
We remark that this technique is increasing the dimension of the global reduced system. However, it is usually still much smaller then $\Cal N$. For the sake of simplicity and for a consistent construction of state and adjoint space, we always choose $N_{\text{max}}$ and $N$ equal for all the involved variables.

\begin{remark}[Supremizer Stabilization] Dealing with Navier-Stokes governing equations,
one has to take care not only of the global inf-sup condition \eqref{ROM_infsup}, but also with the state equation inf-sup condition. Indeed, Navier-Stokes problem can be recast as saddle point itself, which results in a nested saddle point structure when a Navier-Stokes problem is used as state equation of an optimal control problem. The aggregated space techniques has to be accompanied by \emph{supremizer stabilization} for the reduced velocity space.
This approach \cite{rozza2007stability} consists in defining a supremizer operator
$T^{\boldsymbol \mu}: \mathbb P^{\Cal N_p} \rightarrow{{\mathbb V}}^{\Cal N_v}$ as %follows:
% \begin{equation}

$$(T^{\boldsymbol \mu} s, \phi)_{\mathbb {V}} = b(\phi, s; \boldsymbol \mu),  \quad \forall \phi \in {{\mathbb V}}^{\Cal N_v},$$
where $b\cd$ is the bilinear form representing the continuity equation defined in Section \ref{sec:NS}.
Then, we enrich the reduced velocity space through the pressure supremizers as follows:
$$
{{\mathbb V}}_N = \text{span}\{ \chi^{v}_n, \; \chi^{{T_p}}_n, \chi^{w}_n, \; \chi^{{T_q}}_n, \; n = 1,\dots,N\},
$$
where $\chi^{T_p}_n$ and $\chi^{T_q}_n$ are the basis supremizers obtained by state and adjoint pressure snapshots, respectively. Enlarging in this way the reduced space for velocity will guarantee inf-sup stability for the Navier-Stokes equation. This approach, i.e. supremizer stabilization combined with aggregated spaces, is the key for the well-posedness of the whole optimality system \eqref{NS_ocp}. This will lead to a reduced system of dimension $13N$, which is still convenient compared to the global FE approximation dimension.
\end{remark}

\subsection{Numerical Results}
\label{rom_results}
\begin{figure}
\centering
\begin{subfigure}[b]{0.49\textwidth}
\centering
\includegraphics[width=0.75\textwidth]{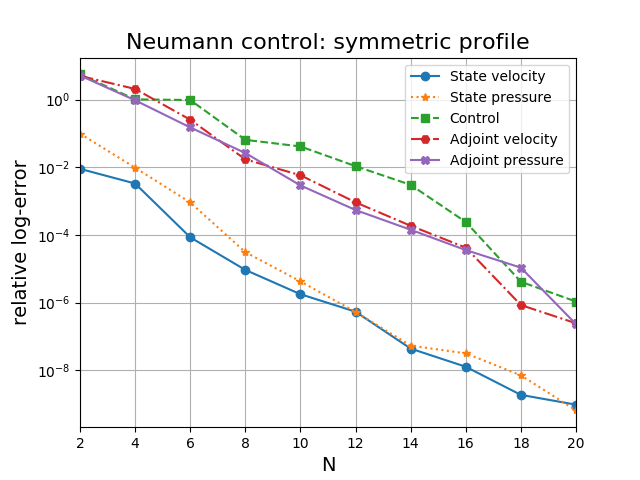}
\caption{}
\label{fig:neumann_mean_s}
\end{subfigure}
\hfill
\begin{subfigure}[b]{0.49\textwidth}
\centering
\includegraphics[width=0.75\textwidth]{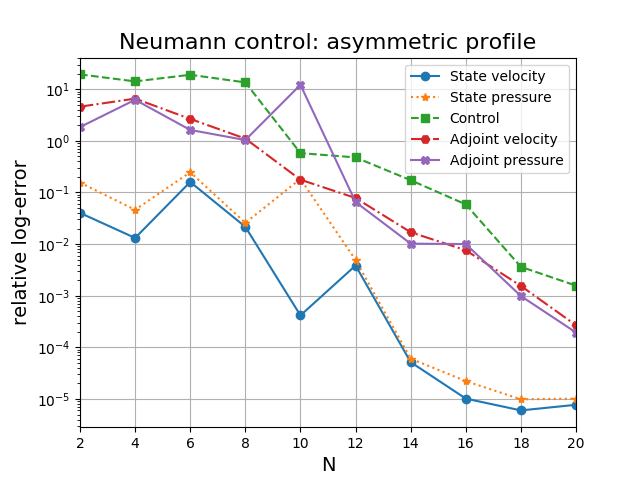}
\caption{}
\label{fig:neumann_mean_as}
\end{subfigure}\\
\caption{Average error over $\mu$ with $\overline N = 20$ and $\alpha = 0.01$ for symmetric and asymmetric profile in (a) and (b) for Neumann control, respectively.}
\label{fig:rom_neumann_av}
\end{figure}

\begin{figure}
\centering
\begin{subfigure}[b]{0.49\textwidth}
\centering
\includegraphics[width=0.75\textwidth]{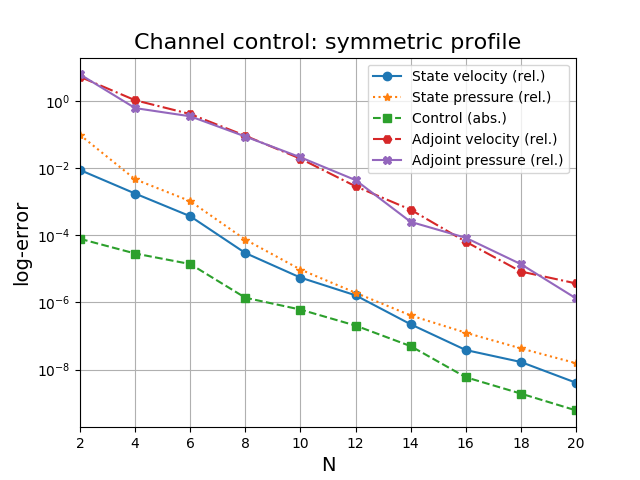}
\caption{}
\label{fig:channel_mean_s}
\end{subfigure}
\hfill
\begin{subfigure}[b]{0.49\textwidth}
\centering
\includegraphics[width=0.75\textwidth]{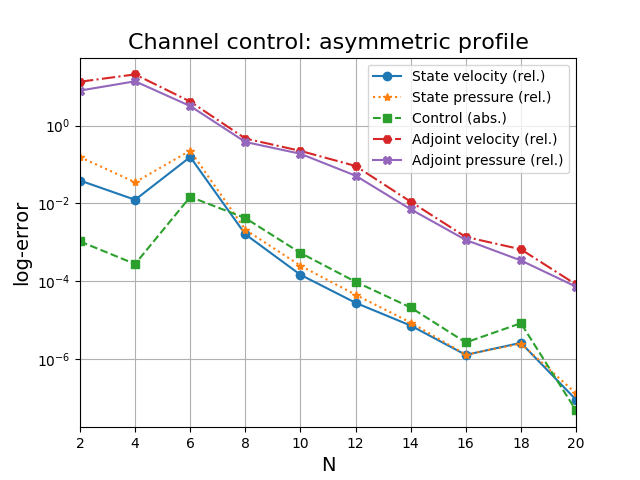}
\caption{}
\label{fig:channel_mean_as}
\end{subfigure}\\
\caption{Average error with $\overline N= 20$ over $\mu$ for symmetric ($\alpha = 1$) and asymmetric ($\alpha = 0.01$) profile in (a) and (b) for Channel control, respectively.}
\label{fig:rom_channel_av}
\end{figure}

In this Section we present the numerical results deriving from the reduction of the four controlled test cases described in Section \ref{NS_ocp}. For each numerical test case, the offline setting is given by $N_\text{max} = 51$ snapshots evaluated for equidistant parameters in the range of $\mathcal P = [0.5, 2]$ and the POD algorithm is chosen for the ROM construction. Let us define the \emph{basis number} $\overline N$ as the maximum value of $N$, i.e. $N \in \{1, \hdots, \overline N\}$. For Dirichlet test case we chose $\overline N = 12$ basis functions, while for the other test cases the basis number is $\overline N = 20$. \C{Such value is chosen in analogy with reduction results obtained in the uncontrolled scenario, see e.g.\ \cite{pichi2021artificial,pintore2019efficient,khamlich2021model}}. \A{The former choice is due to the presence of two multiplier variables, which increase the global ROM dimension (from $13 \overline N$ to $15 \overline N$), and jeopardize the robustness of the reduced nonlinear solver. We remark that the final reduced systems are still much smaller than their high-fidelity counterparts, which involve from 50 to 70 thousands of degrees of freedom, depending on the type of control imposed, boundary or distributed, respectively.} We
perform an online phase solving \eqref{ROM_ocp} for $151$ equidistant value of $\mu$ in the same parameter space $\Cal P$. The performance has been tested through separate error analysis for each variable. The reliability of the ROM approach has been evaluated through
\begin{itemize}
\item[{$\small{\circ}$}] an average error over the parameter space against an increasing value of the reduced spaces dimension $N$ from one up to $\overline N$;
\item[{$\small{\circ}$}] a $\mu$-dependent error computed for the value $\overline N$.
\end{itemize}
The two error analyses highlight different features of the reduced system that we are going to discuss in the following.
Indeed, the average error gives us information about how the reliability changes as the behavior of the solution changes. The straight profile appears to be always the best approximated, due to its Stokes-like (symmetric) nature for all the value $\mu \in \Cal P$. This is the case of Neumann and Channel control, which average error is depicted in Figures \ref{fig:neumann_mean_s} and \ref{fig:channel_mean_s}. Their asymmetric counterparts, Figures \ref{fig:neumann_mean_as} and \ref{fig:channel_mean_as}, show how representing the two different features of the solution, a Stokes-like one for $\mu > \mu^{\ast}$ and a wall-hugging for the lower values of $\mu$, using the same value of $\overline N$ is more difficult than the symmetric one. Nonetheless, the provided accuracy for basis size $\overline N$ is satisfactory for many practical applications in both target cases.
This argument applies to the control and adjoint variables, yet the state ones are the best described by ROM for all the test cases. Because of the optimality equation, the adjoint variables feel the direct influence of the control, which is the most challenging one to be approximated by the reduced model due to its high variability in $\mu$. Indeed, the control variable presents a sort of \emph{on-off} behavior which drastically affects the efficiency of reduced representation. 
\begin{figure}
\centering
\begin{subfigure}[b]{0.49\textwidth}
\centering
\includegraphics[width=0.75\textwidth]{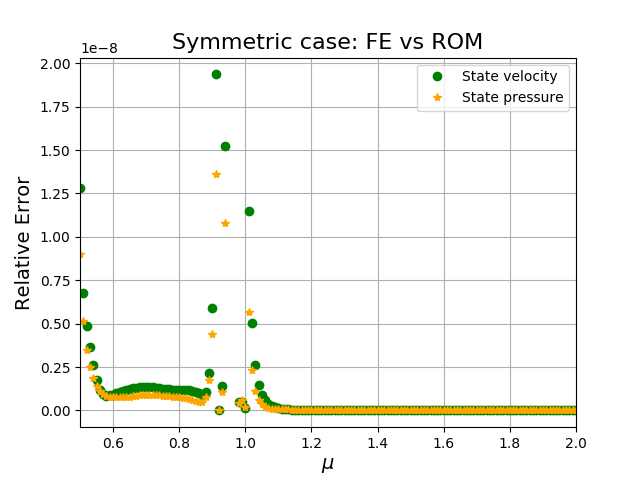}
\caption{}
\label{fig:neumann_mu_s_state}
\end{subfigure}
\hfill
\begin{subfigure}[b]{0.49\textwidth}
\centering
\includegraphics[width=0.75\textwidth]{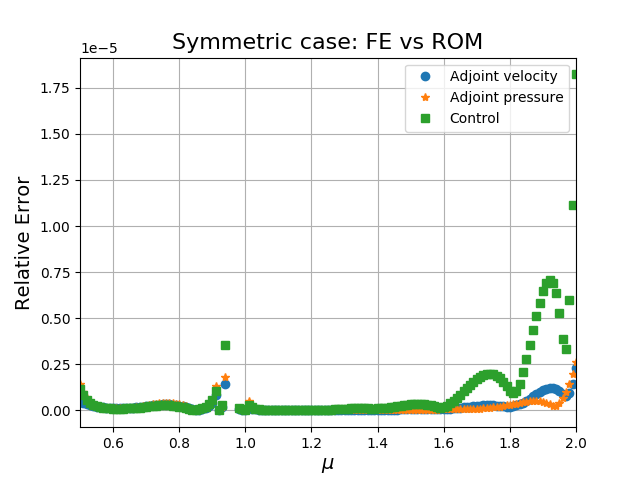}
\caption{}
\label{fig:neumann_mu_s_adj}
\end{subfigure}\\
\caption{The $\mu$-dependent error with $\overline N = 20$ and $\alpha = 0.01$ for symmetric profile for state variable in (a) and adjoint and control variables in (b) for Neumann control, respectively.}
\label{fig:rom_neumann_mu}
\end{figure}

\begin{figure}
\centering
\begin{subfigure}[b]{0.49\textwidth}
\centering
\includegraphics[width=0.75\textwidth]{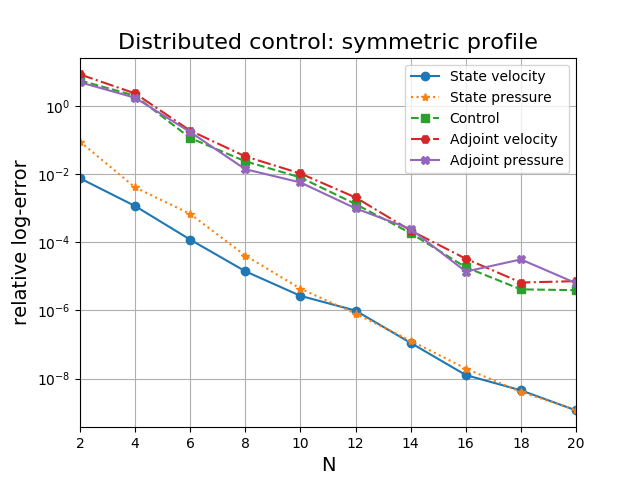}
\caption{}
\label{fig:dist_mean_s}
\end{subfigure}
\hfill
\begin{subfigure}[b]{0.49\textwidth}
\centering
\includegraphics[width=0.75\textwidth]{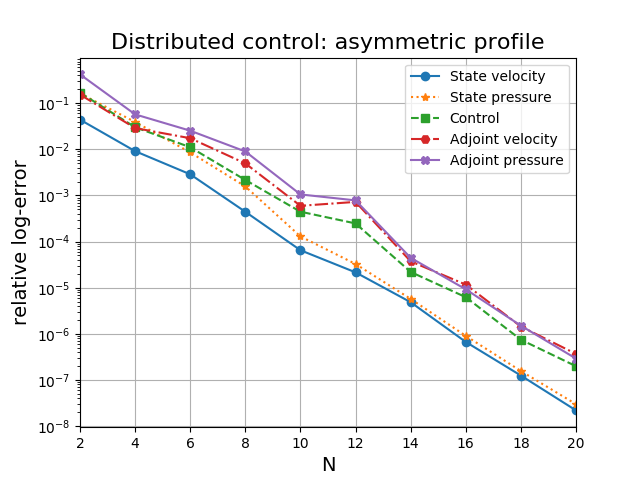}
\caption{}
\label{fig:dist_mean_as}
\end{subfigure}\\
\caption{Average error over $\mu$ with $\overline N= 20$ and $\alpha = 0.01$ for symmetric and asymmetric profile in (a) and (b) for Distributed control, respectively.}
\label{fig:rom_dist_av}
\end{figure}

\begin{figure}
\centering
\begin{subfigure}[b]{0.49\textwidth}
\centering
\includegraphics[width=0.75\textwidth]{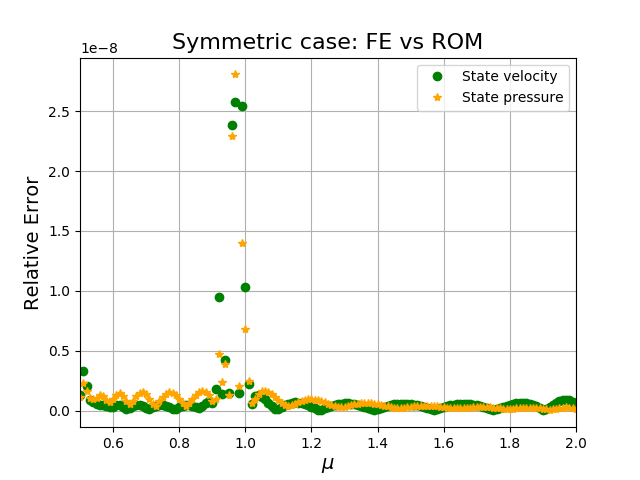}
\caption{}
\label{fig:dist_mu_s_state}
\end{subfigure}
\hfill
\begin{subfigure}[b]{0.49\textwidth}
\centering
\includegraphics[width=0.75\textwidth]{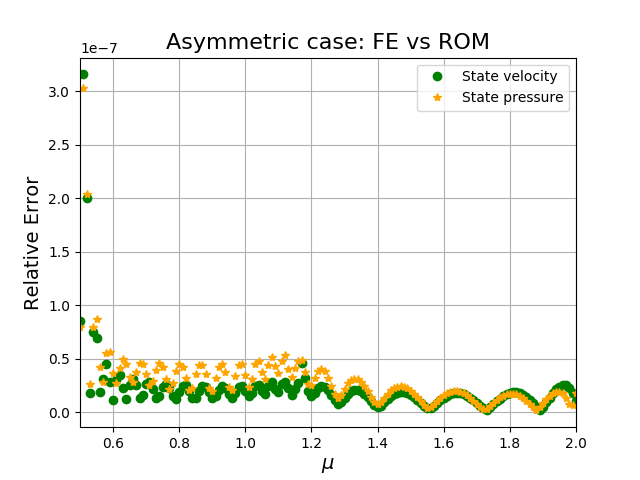}
\caption{}
\label{fig:dist_mu_s_adj}
\end{subfigure}\\
\caption{The $\mu$-dependent error with $\overline N= 20$ and $\alpha = 0.01$ for symmetric and asymmetric profile of the state variable in (a) and (b) for Distributed control, respectively.}
\label{fig:rom_dist_mu}
\end{figure}

\begin{figure}
\centering

\begin{subfigure}[b]{0.49\textwidth}
\centering
\includegraphics[width=0.75\textwidth]{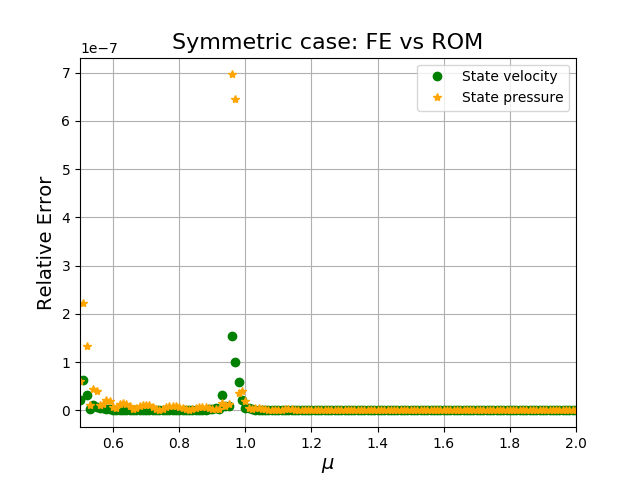}
\caption{}
\label{fig:channel_mu_s_state}
\end{subfigure}
\hfill
\begin{subfigure}[b]{0.49\textwidth}
\centering
\includegraphics[width=0.75\textwidth]{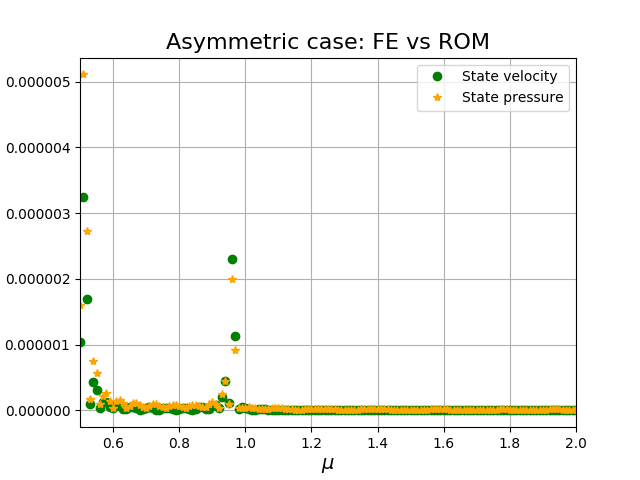}
\caption{}
\label{fig:channel_mu_s_adj}
\end{subfigure}\\
\caption{The $\mu$-dependent error for $\overline N= 20$ for symmetric and asymmetric profile of the state variable in (a) and (b) for Channel control, respectively.}
\label{fig:rom_channel_mu}
\end{figure}

\begin{figure}
\centering
\begin{subfigure}[b]{0.49\textwidth}
\centering
\includegraphics[width=0.75\textwidth]{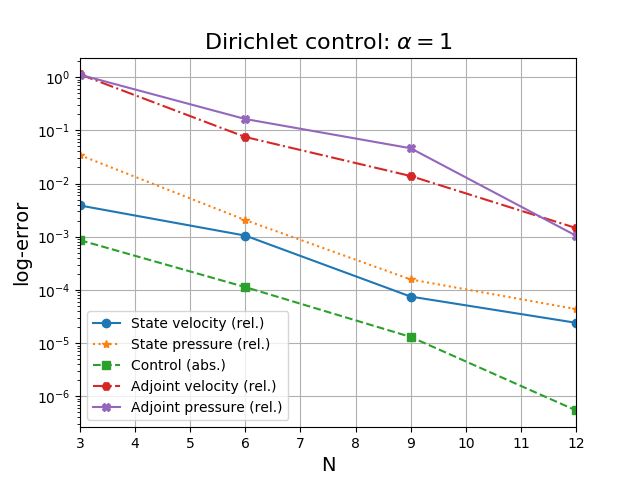}
\caption{}
\label{fig:diri_mean_1}
\end{subfigure}
\hfill
\begin{subfigure}[b]{0.49\textwidth}
\centering
\includegraphics[width=0.75\textwidth]{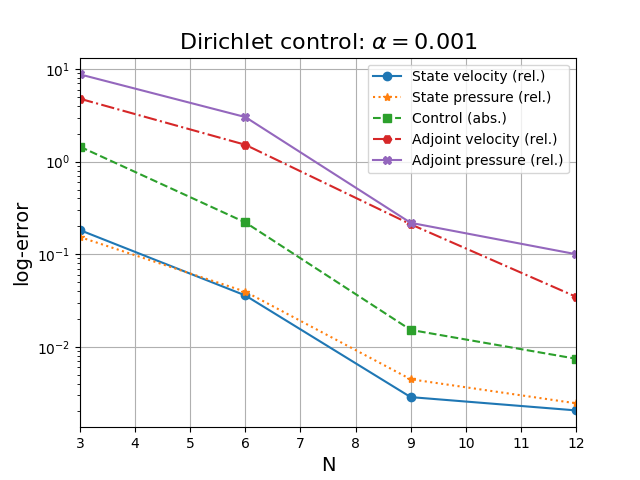}
\caption{}
\label{fig:diri_mean_1e3_av}
\end{subfigure}\\
\caption{Average error over $\mu$ with $\overline N= 12$ for $\alpha = 1$ and $\alpha = 0.001$ in (a) and (b) for Dirichlet control, respectively.}
\label{fig:rom_diri}
\end{figure}

\begin{figure}
\centering
\begin{subfigure}[b]{0.49\textwidth}
\centering
\includegraphics[width=0.75\textwidth]{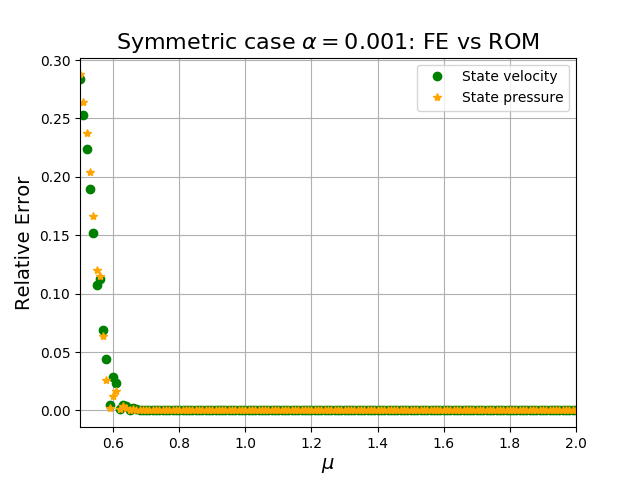}
\caption{}
\label{fig:diri_mu_1e3_state}
\end{subfigure}
\hfill
\begin{subfigure}[b]{0.49\textwidth}
\centering
\includegraphics[width=0.75\textwidth]{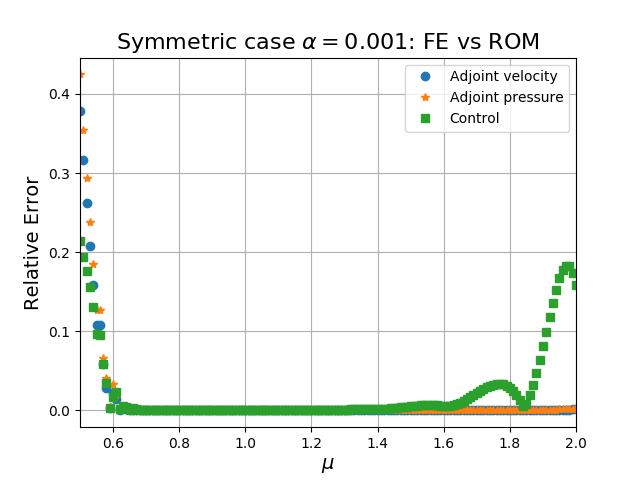}
\caption{}
\label{fig:diri_mu_1e3_adj}
\end{subfigure}\\
\caption{The $\mu$-dependent error for $\alpha = 0.001$ with $\overline N = 12$ for state variable in (a) and adjoint and control variables in (b) for Dirichlet control, respectively.}
\label{fig:rom_diri_mu}
\end{figure}

For example, if we deal with Stokes target $v_{\text{d}}$, the control is $\emph{off}$ for high values of viscosity, but, when $\mu \sim \mu^{\ast}$, it starts to grow in magnitude and to change drastically its features. This is represented in Figures \ref{fig:neumann_mu_s_adj} and \ref{fig:diri_mu_1e3_adj}, where the higher values for the control error and, thus, for the adjoint variable, is for higher values of $\mu$.
Since the control magnitude is essentially zero, for Channel and Dirichlet test cases with low Reynolds number (high viscosity), we chose to plot the absolute errors instead of the relative ones, in order to prevent division by zero.
This is not the case of Distributed control, see for example Figure \ref{fig:dist_mean_as}, which presents good errors decay for all the variables, since its strong action causes the control magnitude to always be a meaningful normalization factor.\\
The most challenging case to be approximated is the Dirichlet one for $\alpha = 0.001$. It is the most complex dynamics, where a new bifurcation appears. Indeed, it can reach an error of almost $10^{-3}$ for the controlled state. Even though this result is worse than the ones obtained for the other test cases, where average error are ranging between $10^{-5}$ and $10^{-8}$, the accuracy provided by the ROM is still acceptable for many practical purposes. We remark that this performance is strictly correlated to the more complex features of the Dirichlet problem. This is evident also from the $\mu$-dependent errors in Figure \ref{fig:diri_mu_1e3_state} and \ref{fig:diri_mu_1e3_adj}. In both the pictures, we see a great increment of the error for high Re. Although the phenomenon appears also for the other test cases, see Figures \ref{fig:dist_mu_s_adj}, \ref{fig:channel_mu_s_state} and \ref{fig:channel_mu_s_adj}, it is not as strong as the Dirichlet case.

Furthermore, the $\mu$-dependent error gives an \emph{a posteriori} information about the bifurcation point. Indeed, in order to have good accuracy property, ROM approach requires regularity on the parametric dependence of the solution instead of the spatial one. This means that reduced errors will generally exhibit a peak at $\mu^*$.
In fact, an increasing value of the error can be seen around $\mu^{\ast} \sim 0.96$ for Neumann, Distributed and Channel control \C{as one can observe from Figures \ref{fig:neumann_mu_s_state}, \ref{fig:dist_mu_s_state} and \ref{fig:channel_mu_s_state} for the state variable, respectively}. 

This feature can be very useful when there is no previous knowledge about bifurcating behaviors. In this sense, ROM is not only an useful approach to solve in a fast way very complicated time consuming systems, but also to detect parameters which can be related to the bifurcating nature of the problem at hand, since their instances will be the worst approximated. \C{Namely, the ROMs confirm ``a posteriori" the location of the bifurcation points.} We conclude this analysis by noticing that the same considerations hold for the Dirichlet control, but this time at the left end of the parametric domain $\mathcal P$ where such phenomenon is clearly influenced by the new configuration observed in Figure \ref{fig:diri_v_1e3}. 

\section{Conclusions}
\label{conclusions}
In this work we proposed a first attempt to steer bifurcation phenomena arising from nonlinear PDE($\bmu$) through an optimal control formulation. First of all, we built a general framework which can be applied to general nonlinear \ocp s and we proposed a global stability analysis trough the solution of the eigenvalue problem associated to the optimization system. We tested the stability analysis over four control problems governed by bifurcating Navier-Stokes equations in a sudden-expansion channel. We studied how the control can affect the classical stable wall-hugging solution of the uncontrolled state equation and we proposed some observations on different configurations and features of the optimal solution based on the eigenvalue systems. \\
Furthermore, we employed ROM for all the test cases, proposing it as a strategy to solve the parametrized analysis of the optimization system in a low-dimensional setting, while confirming the applicability of reduction strategies to complex nonlinear models. \\
We believe that this work is a first step towards a better comprehension of the very complex action of optimal control over a nonlinear bifurcating system. We believe that the content of this work could pave the way to many improvements on the topic. Among them, one could be the deeper analysis of the Dirichlet test case, which gives the more unexpected and, consequently, difficult to interpret, results.

\section*{Acknowledgements}
We acknowledge the support by European Union Funding for Research and Innovation -- Horizon 2020 Program -- in the framework of European Research Council Executive Agency: Consolidator Grant H2020 ERC CoG 2015 AROMA-CFD project 681447 ``Advanced Reduced Order Methods with Applications in Computational Fluid Dynamics''. We also acknowledge the PRIN 2017  ``Numerical Analysis for Full and Reduced Order Methods for the efficient and accurate solution of complex systems governed by Partial Differential Equations'' (NA-FROM-PDEs) and the INDAM-GNCS project ``Tecniche Numeriche Avanzate per Applicazioni Industriali''.
The computations in this work have been performed with RBniCS \cite{rbnics} library, developed at SISSA mathLab, which is an implementation in FEniCS \cite{fenics} of several reduced order modelling techniques; we acknowledge developers and contributors to both libraries.

 \bibliographystyle{abbrv}
\bibliography{maria,fede}

\end{document}